\documentclass[11pt,a4paper,reqno]{article}
\usepackage{amsmath,amssymb,exscale,mathrsfs}
\usepackage{graphicx,color,epsfig}%les deux lignes a rajouter pour le graphique
\usepackage{subfigure}
\usepackage[subfigure]{ccaption}
\usepackage{float}

\usepackage{enumerate}

\usepackage{color}
\definecolor{marin}{rgb}   {0.,   0.3,   0.7} 
\definecolor{rouge}{rgb}   {0.8,   0.,   0.} 
\definecolor{sepia}{rgb}   {0.8,   0.5,   0.} 
\usepackage[colorlinks,citecolor=marin,linkcolor=rouge,
            bookmarksopen,
            bookmarksnumbered
           ]{hyperref}

\newtheorem{lemma}{Lemma}[section]

\newtheorem{theorem}[lemma]{Theorem}
\newtheorem{proposition}[lemma]{Proposition}

\newtheorem{remark}[lemma]{Remark}
\newtheorem{example}[lemma]{Example}

\newtheorem{notation}[lemma]{Notation}
\newtheorem{definition}[lemma]{Definition}
\newtheorem{conclusion}[lemma]{Conclusion}
\newtheorem{assumption}{Assumption}

\numberwithin{equation}{section}
\newcommand{\QED}{\mbox{}\hfill \raisebox{-0.2pt}{\rule{5.6pt}{6pt}\rule{0pt}{0pt}} 
          \medskip\par}             
\newenvironment{Proof}{\noindent
    \abovedisplayskip = 0.5\abovedisplayskip
    \belowdisplayskip=\abovedisplayskip{\bfseries Proof. }}{\QED}

\newcommand{\N}{\mathbb{N}}

\newcommand{\R}{\mathbb{R}}

\newcommand{\C}{\mathbb{C}}

\newcommand{\T}{\mathbb{T}}
\newcommand{\Z}{\mathbb{Z}}

\newcommand{\cF}{{\cal F}}
\newcommand{\cG}{{\cal G}}

\newcommand{\be}{\begin{equation}}
\newcommand{\ee}{\end{equation}}
\newcommand{\bea}{\begin{eqnarray}}
\newcommand{\eea}{\end{eqnarray}}
\newcommand{\bee}{\begin{eqnarray*}}
\newcommand{\eee}{\end{eqnarray*}}

\def\ni{\noindent}
\def\bs{\bigskip}
\def\ms{\medskip}
\def\ss{\smallskip}
\def\eps{\varepsilon}

\def\pref#1{{\rm \ref{#1}}}

\def\I{\mbox{\rm I}}
\def\pa{\partial}

\def\ub{\underbar}
\def\elle{\ell \hspace*{-0.7pt}e}
\def\te{t \hspace*{-0.7pt}e}

\author{Philippe Chartier \thanks{INRIA-Rennes Bretagne Atlantique, IPSO Project} 
\and 
Nicolas Crouseilles \thanks{INRIA-Rennes Bretagne Atlantique, IPSO Project} 
  \and 
Mohammed Lemou
 \thanks{CNRS and IRMAR, Universit\'e de Rennes 1 and INRIA-Rennes Bretagne Atlantique, IPSO Project}
  \and 
Florian M\'ehats 
 \thanks{IRMAR, Universit\'e de Rennes 1 and INRIA-Rennes Bretagne Atlantique, IPSO Project}}
\title{Uniformly accurate numerical schemes for highly oscillatory Klein-Gordon and  nonlinear Schr\"odinger equations}        
\begin{document}
\maketitle

\begin{abstract}
This work is devoted to the numerical simulation of nonlinear Schr\"odinger and Klein-Gordon equations.
We present a general strategy to construct numerical schemes which are {\em uniformly accurate} with respect to
the oscillation frequency.  This is a stronger feature than the usual so called ``Asymptotic preserving"  property, the last being also satisfied by our scheme
in the highly oscillatory limit. Our strategy enables to simulate the oscillatory problem without
using any mesh or time step refinement, and the orders of our schemes are preserved uniformly in all regimes.  In other words, since our numerical method
is not based on the derivation and the simulation of asymptotic models, it works in
the regime where the solution does not oscillate rapidly,  in the highly oscillatory limit
regime,  and in the intermediate regime with the same order of accuracy.  In the same spirit as in  \cite{clm}, the method is based on two main ingredients. First, we embed our problem in a {\em suitable ``two-scale" reformulation}  with the introduction of an additional variable. Then a link is made with classical strategies based on Chapman-Enskog expansions
in kinetic theory despite the dispersive  context of the targeted equations, allowing to separate the fast time  scale from the slow one.  Uniformly accurate (UA)  schemes are  eventually derived from this new formulation and their properties and performances are assessed both theoretically and numerically.
 
\end{abstract}

\tableofcontents
%%%%%%%%%%%%
%%% INTRODUCTION
%%%%%%%%%%%%

\section{Introduction} 
This work is concerned with the numerical solution of highly-oscillatory differential equations in an infinite dimensional setting. Our main two applications here are the nonlinear Schr\" odinger equation and the nonlinear Klein-Gordon equation, although, prior to addressing them specifically, we envisage the more general situation of an abstract differential equation in a Hilbert space. To be a bit more specific, we shall consider equations of the form 
\begin{equation} \label{eq:IVP}
\frac{d}{dt} u^\eps(t) = \cF(t,t/\eps,u^\eps(t)),\qquad t \in [0,T],  \quad u^\eps(0) = u_0, 
\end{equation}
where the vector field $(t,\tau,u) \mapsto \cF(t,\tau,u)$ is supposed to be periodic of period $P$ with respect to the variable $\tau$ (we shall denote $\T \equiv \R/(P\Z)$). The parameter $\eps$ is supposed to have a positive real value in an interval of the form $]0,\eps_0]$ for some $\eps_0 >0$. However, $\eps$ is not necessarily vanishing and may be as well thought of as being close to $1$: this means  we can consider
equation \eqref{eq:IVP} simultaneously in different regimes, namely highly-oscillatory for small values of $\eps$ or smooth for larger values of $\eps$, and our aim is to design a versatile numerical method, capable of handling these two extreme regimes as much as all intermediate ones. 

Generally speaking, standard numerical methods for equation \eqref{eq:IVP} exhibit errors of the form $\Delta t^p/\eps^q$ for some positive $p$ and $q$. The user of such methods is thus  forced to restrict the step-size $\Delta t$ to values less than $\eps^{q/p}$ in order to obtain some accuracy. This becomes an unacceptable constraint for vanishing values of $\eps$. Whenever equation \eqref{eq:IVP} admits a limit  model, {\em Asymptotic-Preserving} (AP) schemes \cite{jin} have been designed to overcome this restriction: the methods we construct obey the corresponding requirement, i.e. they degenerate into a consistent numerical scheme for the limit model whenever $\eps$ tends to zero. 

As favorable as this property may seem, the error behavior of an AP scheme may deteriorate for ``intermediate" regimes where $\eps$ is neither very small nor large. The derivation of asymptotic models  for \eqref{eq:IVP} has been the subject of many works  --\,see e.g. \cite{strobo,cmss,perko,SV} for time-averaging techniques and \cite{allaire,frenod-raviart,nguetseng} for homogenization techniques\,-- and a hierarchy of averaged vector fields and models at different orders of $\eps$ can be classically written from  asymptotic expansions of the solution. However, these asymptotic models are valid only when $\varepsilon$ is small enough and any numerical methods based on the direct approximation of such averaged vector fields introduce a truncation index $n$ and a corresponding incompressible error $\eps^n$. 

In sharp contrast, our strategy in this paper consists in developing numerical schemes that solve directly \eqref{eq:IVP} for a wide range of $\eps$-values with uniform accuracy.
%In particular, we construct a numerical method in the so-called {\em Asymptotic Preserving} (AP) class 
%\cite{jin}: such schemes are consistent with the  model \eqref{eq:IVP} for all positive value of $\eps$, and 
%degenerate into consistent schemes with the asymptotic  (or time-averaged) model when $\eps \to 0$. 
%%%%%%%
%In this article, we are interested in the construction of  efficient numerical schemes for nonlinear 
%Schr\"odinger equations which involve rapid oscillations in time.  We consider
% the solution $u^\eps(t, x)$ of the following 
%nonlinear Schr\"odinger equation 
%\be
%\label{schro}
%i\partial_t u^\eps = -\frac{1}{\varepsilon}\Delta u^\eps + f(|u^\eps|^2) u^\eps, \;\; u^\eps(t=0, x)=u_0(x),  
%\ee
%with $\eps>0$, $t\geq 0$ corresponds to the time variable, $x\in [0, 2\pi]^d$, 
%$d$ being the dimension. Throughout  this paper, the function $f$  is assumed to be smooth, say $\cal C^\infty $, real-valued function. 
% We finally note that  in our setting, the set of the eigenvalues of the Laplacian operator is the set of all the integers. 
The main output of our work are numerical methods for highly 
oscillatory equations of type \eqref{eq:IVP}, which are {\em uniformly accurate} (UA) with respect to the parameter 
$\varepsilon \in ]0,\eps_0], \ \eps_0>0$. These methods, as we shall demonstrate,  are able to capture the various 
scales occurring in the system, while keeping numerical parameters (for instance $\Delta t$) independent of the degree of 
stiffness ($\eps$).

%To prove that our numerical schemes enjoy the above described, and numerically observed, nice  properties, we perform a convergence analysis of the numerical schemes and prove that indeed the
%accuracy of  our schemes is independent of $\eps\in [0,\eps_0]$.  Numerical schemes of order $2$ in time are constructed  and the uniform preservation of this order with respect to $\eps$  is proved: the error in time between the numerical solution and the exact solution is uniformly 
%bounded (with respect to $\eps$) by $C \Delta t^2$, where $\Delta t$ is the time step.
The main idea underlying our strategy (see also \cite{clm}) consists in separating the two time scales naturally present in \eqref{eq:IVP}, namely the slow time $t$ and the fast time $t/\eps$. To this aim, we embed the solution $u^\eps$ into a two-variable function $(t,\tau) \in[0,T] \times \T \mapsto U^\eps(t,\tau)$ while imposing that $U^\eps$ coincides with $u^\eps$ on the {\em diagonal} $\tau=t/\eps$. Clearly, this implies that 
$$
\frac{d}{dt} u^\eps(t) = \pa_t U^\eps(t,t/\eps) + \frac{1}{\eps} \pa_\tau U^\eps(t,t/\eps) = \cF(t,t/\eps,U^\eps(t,t/\eps)).
$$
By virtue of the ``separation" principle, we then consider the equation over the whole $(t,\tau)$-domain, i.e.
\begin{equation} \label{eq:TR}
\pa_t U^\eps(t,\tau) + \frac{1}{\eps} \pa_\tau U^\eps(t,\tau) = \cF(t,\tau,U^\eps(t,\tau)).
\end{equation}
An observation of paramount importance is that no initial condition for \eqref{eq:TR} is evident, since only the value $U^\eps(0,0)=u_0$ is prescribed: consequently, as such, the transport equation \eqref{eq:TR} is not a Cauchy problem and may have many solutions. This apparent obstacle is in fact the way out to our numerical difficulties: given that for any smooth initial condition $U^\eps(0,\tau)=U^\eps_0(\tau)$ satisfying $U^\eps_0(0)=u_0$, we can recover the solution $u^\eps$ from the values of $U^\eps$ on the diagonal $\tau=t/\eps$, the missing Cauchy condition should be regarded as an additional degree of freedom. 

Now, it turns out that for some specific choice of $U^\eps_0$, it is possible to prove that $U^\eps$ and its time-derivatives are bounded on $[0,T] \times \T$ uniformly w.r.t. $\eps$. The point is, in this two-scale 
formulation \eqref{eq:TR} of \eqref{eq:IVP}, that stiffness is confined in the sole term 
$\frac{1}{\eps}\pa_\tau U^\eps$. Interpreting this singularly 
perturbed term as a ``collision" operator, we can derive
the asymptotic behavior of $U^\eps$ through a Chapman-Enskog expansion (see for instance \cite{degond}) from which 
averaged models (first and second order) can  easily  be obtained.  The initial datum $U^\eps_0$ is then chosen so as to satisfy  this expansion at $t=0$, a requirement compatible with $U^\eps_0(0)=u_0$. Two numerical schemes
are then proposed for this augmented problem, following the strategy in \cite{clm}. In the present work, 
these schemes are proved to be uniformly accurate with respect to $\eps$: they have respectively  orders $1$ and $2$ {\it uniformly in $\eps$}. These properties are assessed by numerical experiments on the nonlinear Klein-Gordon and Schr\"odinger equations.

\bs

This paper is organized as follows. 
In Section \ref{sect:setting}, we present the two-scale 
formulation in a general framework and perform in Subsection \ref{chapman} the 
Chapman-Enskog expansion of $U^\eps$. 
The question of the choice of the initial datum  $U^\eps(0,\tau)$ for this augmented 
equation \eqref{eq:TR} is addressed. In Section \ref{sect:fo}, 
a first-order numerical scheme is introduced and analyzed while a second-order one is similarly studied in Section \ref{sect:so}. Finally, Section \ref{sect:num} is 
devoted to a series of numerical tests which confirm the theoretical properties of our schemes when applied to the Schr\"odinger and Klein-Gordon equations and demonstrate the relevance of our strategy.

%%%%%%%%%
%%% SECTION 2
%%%%%%%%%

\section{Two-scale formulation of the oscillatory equation}
In this section, we formulate and analyze the equation obtained by decoupling the slow variable $t$ and the fast one $\tau=t/\eps$. 

\subsection{Setting of the problem}\label{sect:setting}
Given $\eps_0>0$, we  consider the following highly-oscillatory evolution problem
\begin{equation}
\label{eqftildegene}
\partial_t u^\eps = \cF(t,t/\eps,  u^\eps), \qquad t \geq 0, \quad \eps \in ]0,\eps_0], \quad u^\eps(0)=u_0,
\end{equation}
where the unknown $t  \mapsto  u^\eps(t)$ is a smooth map onto a Sobolev  space $H^s$ (either $H^s(\T_x^d)$ or $H^s(\R^d)$ with $d\geq 1$) and the vector-field $(t,\tau, u) \mapsto \cF(t,\tau,u) \in H^s$ is a smooth map, $P$-periodic w.r.t.  $\tau\in\T$ ($\T \equiv \R/P\Z$). Let us emphasize that $\cF$ may also depend on $\eps$, although we shall not reflect specifically this dependence: whenever this is the case, all bounds on $\cF$ and its derivatives then implicitly hold {\em uniformly} in $\eps$. In order to work in Banach algebras, we require that $s>d/2+8$, a condition whose necessity will become apparent for the numerical schemes. 

As described in the Introduction section, we envisage $u^\eps(t)$ as the diagonal solution of the following transport equation which constitutes our starting point:
\be
\label{eqF}
\pa_t U^\eps + \frac{1}{\eps}\pa_\tau U^\eps = \cF(t,\tau, U^\eps),\qquad U^\eps(0,\tau)=U_0^\eps(\tau),
\ee
where the unknown is now the function $(t,\tau)\mapsto U^\eps(t,\tau) \in H^s$. The choice of the Cauchy condition $U_0^\eps(\tau)$ is discussed below, but it is already clear that $u^\eps(t)$ and $U^\eps(t,t/\eps)$ coincide provided that $U_0^\eps(0)=u_0$. 

For our purpose, we shall need that the vector field $\cF$ obeys the following assumption, where each derivation w.r.t. $t$ or $\tau$ typically costs 2 derivatives in the space variable. Indeed, for applications to nonlinear Klein-Gordon or Schr\"odinger equation --\,see \eqref{cF1} and \eqref{cF2}\,--, one has in mind vector fields of the form 
$$\cF(t,\tau, U^\eps)=f\left(e^{it\Delta}U^\eps,e^{i\tau\Delta}U^\eps\right),$$
where $f$ is a smooth function.

\ms
\ni
\begin{assumption}
\label{mainassump} For all $\alpha \in \{0 \cdots 3\}$, $\beta \in  \{0,1\}$ and $\gamma \in \{0 \cdots 3\}$, for all $s,\sigma$ such that $s\geq \sigma>2(\alpha +\beta)+d/2$, the functional $\pa_t^\alpha \pa_\tau^\beta \pa_u^\gamma \cF$ is continuous and locally bounded from  $\R_+ \times \T \times H^{s}$ to
$${\cal L}( \underbrace{H^\sigma\times\ldots\times H^\sigma}_{\gamma \;\;  \mbox{\small{times}}},H^{\sigma-2(\alpha+\beta)}).$$
\end{assumption}
 \subsection{Bounds  in $H^\sigma$ of the solution of the transport equation}
Similar related transport equations will occur  in our analysis, with possibly other functions than $\cF$ and initial conditions with various regularities. A somehow preliminary result thus concerns the existence and uniqueness of the solution of the  {\em general} Cauchy problem 
\be \label{eq:cp}
\pa_t \Phi^\eps + \frac{1}{\eps}\pa_\tau \Phi^\eps = G(t,\tau, \Phi^\eps),  \quad \Phi^\eps(0,\tau) = \Phi^\eps_0(\tau) \in H^\sigma,
\ee
where $\Phi^\eps_0$ (possibly) depends on $\eps \in ]0,\eps_0]$ and is assumed to be not identically zero.
\begin{proposition} \label{prop:exist}
Let $T>0$, let $\sigma> d/2$ and suppose that $G$ is a locally Lipschitz continuous map from $[0,T] \times \T \times H^\sigma$ into $H^\sigma$, that it admits derivatives $\pa_\tau G$ and $\pa_u G$ which are continuous and locally bounded from $[0,T] \times \T \times H^\sigma$ into, respectively, $H^{\sigma-2}$ and  ${\cal L}(H^{\sigma-2},H^{\sigma-2})$. If $\Phi^\eps_0 \in C^0(\T;H^\sigma)\cap C^1(\T;H^{\sigma-2})$ is uniformly bounded in $\eps \in ]0,\eps_0]$ with respect to the $L^\infty_\tau(H^\sigma)$ norm then, for any $ \kappa > 1$, there exists $0<T_\kappa\leq T$ such that for all $\eps\in ]0, \eps_0]$, equation \eqref{eq:cp}  has a unique solution $\Phi^\eps \in C^0([0,T_\kappa]\times \T;H^\sigma)$
and we have  
\begin{equation} \label{eq:est}
 \forall t\in [0,T_\kappa],\qquad \sup_{\eps \in ]0,\eps_0]} \|\Phi^\eps(t,\cdot)\|_{L^\infty_\tau(H^\sigma)} \leq \kappa \;  \sup_{\eps \in ]0,\eps_0]} \|\Phi^\eps_0(\cdot)\|_{L^\infty_\tau(H^\sigma)}.
\end{equation}
%$$
%T_\kappa = \min\Big(T,\frac{\log( \kappa)}{C_G(c \, \kappa \,   \sup_{\eps \in ]0,\eps_0]} \|\Phi^\eps_0(\cdot)\|_{L^\infty_\tau(H^\sigma)})}\Big),
%$$ 
%where $c$ is the norm of the continuous embedding of $H^\sigma$ into $L^\infty$. 
Moreover, $\Phi^\eps$ has first derivatives w.r.t. both $t$ and $\tau$ which are functions of $C^0( [0,T_\kappa] \times \T;H^{\sigma-2})$. If in addition, $G$ satisfies the estimate
$$
\forall (t,\tau) \in [0,T] \times \T, \; \forall v \in H^\sigma,\quad 
\|G(t,\tau,v)\|_{H^\sigma} \leq C_G \|v\|_{H^\sigma} + D_G
$$
for some positive constants $C_G$ and $D_G$, then equation \eqref{eq:cp}  has a unique solution in $C^0([0,T]\times \T;H^\sigma)$
satisfying 
$$
 \forall t\in [0,T],\qquad \sup_{\eps \in ]0,\eps_0]} \|\Phi^\eps(t,\cdot)\|_{L^\infty_\tau(H^\sigma)} \leq  (\sup_{\eps \in ]0,\eps_0]} \|\Phi^\eps_0(\cdot)\|_{L^\infty_\tau(H^\sigma)} + D_G t) e^{t C_G}.
$$
\end{proposition} 
\begin{Proof}
Considering a smooth solution $\Phi^\eps(t,\tau)$ of \eqref{eq:cp} and denoting $\varphi^\eps(t,\tau) = \Phi^\eps(t,\tau+t/\eps)$, it is easy to check that 
$$
\pa_t \varphi^{\eps}(t,\tau) = \pa_t \Phi^\eps(t,\tau+t/\eps) + \frac{1}{\eps} \pa_\tau \Phi^\eps(t,\tau+t/\eps) = G(t,\tau+t/\eps,\varphi^\eps(t,\tau)),
$$
so that the smooth function $t  \mapsto \varphi^\eps(t,\tau)$, parametrized by $(\tau,\eps) \in \T \times ]0,\eps_0]$, is then solution of the ordinary differential equation 
\be \label{eq:cpb}
\pa_t \varphi^\eps(t,\tau)  = G(t,\tau+t/\eps, \varphi^\eps(t,\tau)), \quad \varphi^\eps(0,\tau) = \Phi^\eps_0(\tau).
\ee
According to Cauchy-Lipschitz theorem in $H^\sigma$ (a Banach space), equation \eqref{eq:cpb} has a unique maximal solution on an interval of the form $[0,T^{\eps}_{max}[$ (when $0<T_{max}^\eps<T$) or a solution 
on $[0,T]$, which furthermore satisfies the following inequality
$$
\|\varphi^\eps(t,\tau)\|_{H^\sigma} \leq \|\Phi^\eps_0(\tau)\|_{H^\sigma} + \int_0^t \|G(\theta,\tau+\theta/\eps,\varphi^\eps(\theta,\tau)\|_{H^\sigma} d\theta.
$$
Denote $$R=  \sup_{\eps \in ]0,\eps_0]} \|\Phi^\eps_0(\tau)\|_{H^\sigma}$$
and
$$M_\kappa=\sup\{\|G(t,\tau,u)\|_{H^\sigma},\quad 0\leq t\leq T,\quad \tau\in \T,\quad\|u\|_{H^{\sigma}}\leq \kappa R\}.$$
Now, as long as $\|\varphi^\eps(t,\tau)\|_{H^\sigma} \leq R$, we have
$$
\|\varphi^\eps(t,\tau)\|_{H^\sigma} \leq \|\Phi^\eps_0(\tau)\|_{H^\sigma} + tM_\kappa,
$$
so that 
$$T^{\eps}_{max} \geq T_\kappa:=\min\Big(T,\frac{(\kappa-1)R}{M_\kappa}\Big)>0$$ and 
%$\Phi^\eps$ satisfies 
estimate \eqref{eq:est} holds.
% In particular, for all $(t,\tau,\eps) \in [0,T_{\kappa}] \times \T \times ]0,\eps_0]$, $\Phi^\eps(t,\tau) \in K$. 
Now,  since $\pa_\tau G$ (resp. $\pa_uG$) is a continuous and locally bounded function from $[0,T] \times \T \times H^\sigma$ to $H^{\sigma-2}$ (resp. to ${\cal L} (H^{\sigma-2}, H^{\sigma-2})$), then $\psi^\eps(t,\tau) := \partial_\tau \varphi^\eps(t,\tau)$ is the unique solution on $[0,T_\kappa]$ of the {\em linear} differential equation in $H^{\sigma-2}$
\begin{align*}
\partial_t \psi^\eps(t,\tau) &= (\pa_\tau G)(t,\tau+t/\eps,\varphi^\eps(t,\tau)) +\pa_u G(t,\tau+t/\eps,\varphi^\eps(t,\tau)) \psi^\eps(t,\tau), \\
\psi^\eps(0,\tau) &= \pa_\tau \Phi^\eps_0(\tau) \in H^{\sigma-2}.
\end{align*}
%given that both terms $(\pa_\tau G)(t,\tau+t/\eps,\varphi^\eps(t,\tau))$ and $\pa_u G(t,\tau+t/\eps,\varphi^\eps(t,\tau))$ are bounded on $I \times \T \times K$ (provided that $\pa_\tau \Phi^\eps_0 \in H^{\sigma-2}$). 
  Hence $\varphi^\eps$ has first derivatives
$\pa_t \varphi^\eps$ in $H^\sigma$ and $\pa_\tau \varphi^\eps$ in $H^{\sigma-2}$. Finally, since $\Phi^\eps(t,\tau) = \varphi^\eps(t,\tau-t/\eps)$, we have 
$$
\pa_\tau \Phi^\eps (t,\tau) = \pa_\tau \varphi^\eps(t,\tau-t/\eps) \mbox{ and } \pa_t \Phi^\eps (t,\tau) = \pa_t \varphi^\eps(t,\tau-t/\eps) - \frac{1}{\eps} \pa_\tau \varphi^\eps (t,\tau-t/\eps)
$$
so that $\Phi^\eps$ has also first derivatives in $H^{\sigma-2}$.  Finally, the proof of the subsequent  assertions in Proposition \ref{prop:exist} can be done in the same way, the last estimate being a consequence of the Gronwall lemma.
\end{Proof}
\begin{remark}
From previous formulae, it appears that $\pa_t \Phi^\eps$ exists in $H^{\sigma-2}$ but is not necessarily uniformly bounded in $\eps$. In order to get a solution $\Phi^\eps$ with uniformly bounded first derivatives, we have to consider an appropriate $\eps$-dependent initial condition $\Phi^\eps_0$. In the next two subsections, we shall consider a formal expansion of $\Phi^\eps$ in $\eps$ so as to determine how this initial condition should be prescribed.
\end{remark}
%More generally, the $j$-th derivatives of $v^\eps$ w.r.t to $t$ and $\tau$ exists and $\pa_\alpha \pa_\beta v^\eps(t,\tau,\cdot)$ with $\alpha+\beta \leq j$ belongs to $X^{s-2(\alpha+\beta)}$. We thus get a unique solution $U^\eps(t,\tau,\cdot)$ in $\cap_{j=0}^k C^j([0,T_\kappa] \times \T, X^{s-2j})$.

\subsection{A formal Chapman-Enskog expansion}
\label{chapman}
In this subsection, we analyze {\em formally} the behavior of \eqref{eqF} in the limit $\eps\to 0$ under the 
assumption that its solution $U^\eps$ has uniformly bounded (in $\eps$)  derivatives up to order $3$. Following \cite{clm}, 
we thus consider the linear operator $L$, 
defined for all periodic (regular) function $\tau\in \T\mapsto h(\tau)$ by
$$
Lh=\pa_\tau h.
$$
This operator is skew-adjoint with respect to the $L^2(\T)$ scalar product and its kernel is the set of constant functions. The $L^2$-projector on this kernel is the averaging operator
$$\Pi h:=\frac{1}{P} \int_0^P h(\tau) d\tau$$
which obviously satisfies $\Pi L \equiv 0$.  
On the set of functions with vanishing average, $L$ is invertible with inverse defined by 
$$
(L^{-1}h)(\tau)=(\I-\Pi)\int_0^\tau h(\theta)d\theta.
$$
In order to alleviate notations, we further introduce $A:=L^{-1} (\I-\Pi)$ which operates on the set of periodic functions onto the set of zero-average periodic functions. \\

The Chapman-Enskog expansion (see for instance \cite{degond}) 
consists in writing the solution 
%$U^\eps(t,\tau,x)$ 
$U^\eps(t,\tau)$ 
in the form 
\be
\label{decomp}
U^\eps(t,\tau)=\underbar{U}^\eps(t)+ h^\eps(t,\tau),%\eps h^1(t,\tau,\underbar{U}^\eps(t)) + \eps^2 h^2 (t,\tau,\underbar{U}^\eps(t)) + {\cal O}(\eps^3) 
\ee
where 
$$
\underbar{U}^\eps(t)=\Pi \left(U^\eps(t,\tau)\right), \quad \Pi h^\eps= 0, 
$$
and then, under some regularity assumptions on $U^\eps$ with respect to $t$ and $\eps$, one seeks  the correction $h^\eps$ as an expansion 
in powers of $\eps$:
\be
\label{expansion}
h^\eps(t,\tau)= \eps  h_1(t,\tau,\underbar{U}^\eps(t)) + \eps^2 h_2 (t,\tau,\underbar{U}^\eps(t)) + \dots . 
\ee
Inserting the decomposition \eqref{decomp} into \eqref{eqF} leads to 
\be
\label{eqinsert}
\partial_t  \ub U^\eps  +  \pa_t h^\eps +\frac{1}{\eps} L h^\eps  = \cF(t, \tau,\ub U^\eps+h^\eps).
\ee
Projecting on the kernel of $L$ and taking into account that $\Pi h^\eps = 0$, we obtain \be
\label{eqmoy}
\partial_t \ub U^\eps =  \Pi \left(\cF(t,\tau,\ub U^\eps+ h^\eps)\right),
\ee
and then subtracting from \eqref{eqinsert}
\be
\label{eqh}
\partial_t h^\eps   + \frac{1}{\eps} Lh^\eps  = (\I-\Pi)\left(\cF(t,\tau,\ub U^\eps+ h^\eps)\right).
\ee
Since  $h^\eps$ belongs to the range of $L$, we get
\be
\label{h-esti1}
h^\eps=\eps A \left(\cF(t,\tau,\ub U^\eps+ h^\eps)\right) -\eps L^{-1}( \partial_t h^\eps).
\ee
Therefore, provided $h^\eps$ and its first time derivative are uniformly bounded w.r.t. $\eps$, we first deduce from this last equation that
 $h^\eps =  O(\eps)$.  Now if we additionally assume that the second and third time derivatives are uniformly bounded  w.r.t. $\eps$, then, by a simple induction on \eqref{h-esti1},   we get
$$h^\eps(t,\tau) = \eps h_1(t,\tau,\ub U^\eps) + \eps^2 h_2(t,\tau,\ub U^\eps)+ {\cal O}(\eps^3),$$
with $h_1$ and $h_2$ defined by
\begin{align} 
h_1(t,\tau,U) = &A \cF(t,\tau,U), \label{h1}\\
 h_2(t,\tau,U) =&  A \pa_u \cF(t,\tau, U) A \cF(t,\tau, U) -A^2 \Big(\pa_u \cF(t,\tau, U)\Pi \cF(t,\tau,U) + \pa_t \cF(t,\tau,U) \Big). \label{h2}
\end{align}
 Inserting these corrections into equation \eqref{eqmoy} yields the first and second order averaged models
\begin{align*}
\partial_t \ub U^\eps &=  \Pi \cF(t,\tau, \ub U^\eps) + \mathcal O(\eps), \\
\partial_t \ub U^\eps &=\Pi \cF(t,\tau, \ub U^\eps)+\eps \Pi \left(\pa_u \cF(t,\tau,\ub U^\eps)A\cF(t,\tau,\ub U^\eps)\right)+\mathcal O(\eps^2).
\end{align*}
Anticipating on next sections, let us now briefly address the crucial issue of the initial condition for \eqref{eqF}. According to the above calculations, one expects to get a smooth solution of \eqref{eqF} if  the initial condition $U^\eps_0$ follows the same expansion as above, i.e.
\begin{align}
U^\eps_0(\tau) &= \ub U^\eps_0 + \eps A \cF_0(\tau,\ub U^\eps_0)  \label{eq:initcond} \\
&\quad + \eps^2 \Big(A \pa_u \cF_0(\tau, \ub U^\eps_0) A \cF_0(\tau, \ub U^\eps_0) -A^2 \pa_u \cF_0(\tau, \ub U^\eps_0)\Pi \cF_0(\tau, \ub U^\eps_0) - A^2 (\pa_t \cF)_0 (\tau,\ub U^\eps_0))\Big)  \nonumber 
\end{align}
where we have denoted by a subindex $0$ the evaluation of functions at $t=0$ and where $\ub U^\eps_0 := \ub U^\eps(0)$ is  chosen so as to be compatible with $U^\eps_0(0)=u_0$ which is the initial condition for the original problem \eqref{eqftildegene}. Starting from $\ub U^\eps_0=u_0+{\cal O}(\eps)$ and inserting successively higher-order terms in the previous equation, we can obtain the expression of $\ub U^\eps_0$ and then of $U^\eps_0(\tau)$. For instance, we have 
\begin{equation*}
\ub U^\eps_0 = u_0 + \eps \Pi \int_0^\tau (\I-\Pi) \cF_0(\theta,u_0) d\theta +{\cal O}(\eps^2), 
\end{equation*}
so that 
$$
U^\eps_0(\tau) = u_0 + \eps \int_0^\tau (\I-\Pi) \cF_0(\theta,u_0) d\theta +{\cal O}(\eps^2)
$$
which provides an initial condition for our first order numerical scheme (see Subsection \ref{subfirst}). The explicit computation of second order terms is postponed to Subsection \ref{subsecond}. 
%\subsection{Discussion on the initial data and main result}
\subsection{Estimates of time derivatives}
\label{discuinit}
In this subsection, we indeed prove that the initial condition \eqref{eq:initcond} ensures  that  time derivatives of $U^\eps$ up to  order  $3$ are uniformly bounded in $\eps$. In the sequel, the following functional space will be useful:
\begin{equation}
\label{defX}
%X^\sigma=\bigcap_{0\leq 2\ell<\sigma-d/2}{\cal C}^\ell(H^{\sigma-2\ell}).
X^\sigma=\bigcap_{0\leq 2\ell<\sigma-d/2} C^\ell(\T; H^{\sigma-2\ell}).
\end{equation}

\begin{proposition}
\label{propcondinit} Suppose that $\cF$ satisfies Assumption \pref{mainassump} and let  $s > d/2+8$ and  $\kappa >1$. Consider the following initial condition 
\be
\label{condinit0}
\forall \tau \in \T,\qquad  U^\eps(0, \tau) = U_0^\eps(\tau)=\ub U^\eps_0 + (\eps h_1+\eps^2 h_2)(0,\tau, \ub U^\eps_0)+\eps^3r^\eps(\tau), 
\ee
where $\ub U^\eps_0 \in H^{s+2}$ is assumed to be uniformly bounded in $\eps\in ]0,\eps_0]$, where $h_1$ and $h_2$ are given by \eqref{h1} and \eqref{h2}, and where the remainder term $r^\eps$ is assumed to be bounded in $X^s$ uniformly in $\eps$. Then the following holds:
\begin{enumerate}[(i)]
\item $U^\eps_0$ is uniformly bounded in $L^\infty_\tau(H^s)$ and there exists $T_\kappa>0$ such that, for all $\eps\in ]0,\eps_0]$, equation \eqref{eqF}, subject to the initial condition \eqref{condinit0}, has a unique solution $U^\eps(t,\tau)\in C^0([0,T_\kappa]\times \T;H^s)$, which satisfies the uniform bound
\be
\label{estiunif}
\forall t\in [0,T_\kappa],\qquad \sup_{\eps\in ]0,\eps_0]}\|U^\eps(t)\|_{L^\infty_\tau(H^s)}\leq \kappa \sup_{\eps\in ]0,\eps_0]} \|U^\eps_0\|_{L^\infty_\tau (H^s)}.
\ee
\item Moreover, for any $T_\kappa$ for which \eqref{estiunif} holds, the solution $U^\eps$ satisfies the following estimates
\begin{equation}
\label{unifestideriv}
\forall t\in [0,T_\kappa],\quad \sup_{\eps\in ]0,\eps_0]}\|\partial_t^\alpha \pa_\tau^\beta U^\eps(t)\|_{L^{\infty}_\tau(H^{s-2(\alpha+\beta)})} \leq C,  \;\; \alpha=0,1,2,3, \;\; \beta=0,1, 
\end{equation}
for some constant $C>0$. 
\end{enumerate}
\end{proposition}

\begin{Proof} We prove this proposition in several steps.

\ms
\ni
{\bf Existence of $U^\eps$ and uniform bound.} Let us first estimate the initial condition defined by \eqref{condinit0}. From \eqref{h1}, \eqref{h2}, one gets
\be
\label{detU0eps}
U_0^\eps=\ub U^\eps_0+\eps A \cF_0+\eps^2A \pa_u \cF_0 A \cF_0 -\eps^2 A^2 \pa_u \cF_0\Pi \cF_0 - \eps^2A^2(\pa_t \cF)_0+\eps^3r^\eps,
\ee
where, for conciseness, we have further omitted the dependence\footnote{In the sequel, we explicitly mention the dependence $\cF_0(U^\eps_0)$ while $\cF_0$ stands for $\cF(0,\tau,\ub U_0^\eps)$ and $(\pa_t \cF)_0$ stands for $\pa_t \cF(0,\tau,\ub U_0^\eps)$.}  of $\cF_0$ in $\tau$ and $\ub U^\eps_0$. 
%We notice that, by Assumption \ref{mainassump}, we have
%$(\pa^\alpha_t\pa^\beta_\tau\pa^\gamma_u \cF)_0\in X^{s+2-2(\alpha+\beta)}$ (the spaces $X^\sigma$ are defined by \eqref{defX}) for $2(\alpha+\beta)<s+d/2$, with norms uniformly bounded w.r.t. $\eps$. Hence, observing that $A$ and $\pa_\tau A= (\I-\Pi)$ are bounded operators on $C^0(\T;H^s)$, 
We notice that, by Assumption \ref{mainassump}, we have
$$
\cF_0, \; \partial_u \cF_0 A \cF_0, \; \partial_u \cF_0 \Pi \cF_0 \in X^{s+2}, \;\;\; \mbox{ and } \;\;\; (\partial_t \cF)_0  \in X^s, 
$$
with norms uniformly bounded w.r.t. $\eps$. 
Hence, observing that $A$ and $\Pi$ are bounded operators on $C^0(\T;H^\sigma)$, for all $\sigma$, 
one deduces that $U_0^\eps$ belongs to $X^s$ and is uniformly bounded w.r.t. $\eps$. Hence, according to Proposition \ref{prop:exist}, $U^\eps$ exists on an interval $[0,T_\kappa]$ independent of $\eps$, and satisfies 
$$
\forall 0 \leq t \leq T_\kappa, \qquad \sup_{\eps \in ]0,\eps_0]} \|U^\eps(t,\cdot)\|_{L^\infty_\tau(H^s)} \leq \kappa \sup_{\eps \in ]0,\eps_0]}  \|U^\eps_0(\cdot)\|_{L^\infty_\tau(H^s)}. 
$$
Furthermore, derivatives of $U^\eps$ w.r.t. $t$ and $\tau$ exist and are functions with values in $H^{s-2}$. 

\ms
\ni
{\bf Estimate of the first derivative in $t$.} The first derivative $V^\eps=\pa_tU^\eps$ satisfies the equation  
\be
\label{eqv}
\partial_t V^\eps  + \frac{1}{\varepsilon}\pa_\tau V^\eps = \pa_u \cF(t,\tau,U^\eps) V^\eps
+ \pa_t \cF(t,\tau,U^\eps) 
\ee
with initial condition 
$$
V^\eps_0  = \cF_0(U^\eps_0) - \frac{1}{\varepsilon}L U^\eps_0.
$$
{}From \eqref{detU0eps} and $LA=(\I-\Pi)$, $LA^2=A$, we obtain
$$
V^\eps_0 = \cF_0(U^\eps_0) - \cF_0+\Pi \cF_0 - \eps (\I-\Pi) \pa_u \cF_0 A \cF_0 + \eps A \pa_u \cF_0 \Pi \cF_0 + \eps A (\pa_t \cF)_0-\eps^2Lr^\eps.
$$
Taylor-Lagrange expansions with integral remainder at order one and two give\footnote{The notation ${\cal O}_{X^\sigma}$ is used here for terms uniformly bounded in $\T$ with the appropriate $X^\sigma$-norm.}

\begin{align*}
&\cF_0(U_0^\eps) - \cF_0= \int_0^1 \pa_u \cF_0(\ub U^\eps_0 + \mu (U_0^\eps-\ub U^\eps_0)) d\mu \; \Big(U_0^\eps-\ub U^\eps_0\Big) = {\cal O}_{X^{s}}(\eps)
 \\
&= \pa_u \cF_0(\ub U^\eps_0) \; (U^\eps_0-\ub U^\eps_0) + \int_0^1 (1-\mu) \pa^2_u \cF_0(\ub U_0^\eps + \mu (U_0^\eps-\ub U^\eps_0)) d\mu \; \Big(U_0^\eps-\ub U^\eps_0, U_0^\eps-\ub U^\eps_0\Big) \\
% &= \pa_u \cF_0 \;\; (U^\eps_0-\ub U^\eps_0) + \int_0^1 (1-\mu) \pa^2_u \cF_0(\ub U^\eps + \mu (U_0^\eps-\ub U^\eps_0)) d\mu \; \Big(U_0^\eps-\ub U^\eps_0, U_0^\eps-\ub U^\eps_0\Big) \\
 &= \eps \pa_u \cF_0 A \cF_0+ {\cal O}_{X^{s}}(\eps^2),
\end{align*}
where we used that $r^\eps$ is uniformly bounded in $X^s$. Therefore\footnote{Notice that $LX^\sigma$ is continuously embedded in $X^{\sigma-2}$.}
\begin{align}
V^\eps_0 &= \Pi \cF_0 + {\cal O}_{X^{s-2}}(\eps) \\
&= \Pi \cF_0 + \eps \pa_u \cF_0 A \cF_0   - \eps (\I-\Pi) \pa_u \cF_0 A \cF_0 + \eps A \pa_u \cF_0 \Pi \cF_0 +  \eps A (\pa_t \cF)_0 + {\cal O}_{X^{s-2}}(\eps^2) \nonumber \\
&= \Pi \cF_0 + \eps \Pi \pa_u \cF_0 A \cF_0  + \eps A \pa_u \cF_0 \Pi \cF_0 + \eps A (\pa_t \cF)_0 + {\cal O}_{X^{s-2}}(\eps^2).  \label{eq:V0}
\end{align}
In particular, $V^\eps_0 \in C^0(\T;H^{s-2})\cap C^1(\T;H^{s-4})$ and is uniformly bounded in $L^{\infty}_{\tau}(H^{s-2})$ w.r.t. $\eps$. According to the second part of Proposition \ref{prop:exist} (for $\sigma=s-2$) with
$$
G(t,\tau,V) = \pa_u \cF(t,\tau,U^\eps(t,\tau)) V + \pa_t \cF(t,\tau,U^\eps(t,\tau))
$$
which is a map from $\R_+ \times \T \times H^{s-2}$ into $H^{s-2}$, we thus have an estimate of the form
$$
\forall t \in [0,T_\kappa], \quad \sup_{\eps \in ]0,\eps_0]} \|V^\eps(t,\cdot)\|_{L_\tau^\infty(H^{s-2})} \leq \Big( \sup_{\eps \in ]0,\eps_0]} \|V^\eps_0(\tau)\|_{L_\tau^\infty(H^{s-2})} + D_G T_\kappa\Big) e^{T_\kappa C_G}.
$$

\ms
\ni
{\bf Estimate of the second  derivative in $t$.} 
We proceed in an analogous way for  $W^\eps=\partial_t^2 U^\eps=\pa_t V^\eps \in H^{s-4}$ by considering 
\begin{align}
\label{eqw2}
\partial_t W^\eps + \frac{1}{\varepsilon}\partial_\tau W^\eps
=& \pa^2_u \cF(t,\tau,U^\eps) (V^\eps,V^\eps) + 2 \pa_t \pa_u \cF(t,\tau,U^\eps) V^\eps \nonumber  \\
&+ \pa_t^2 \cF(t,\tau,U^\eps) + \pa_u \cF(t,\tau,U^\eps) W^\eps.  
\end{align}
The initial condition for $W^\eps$ can be obtained from \eqref{eqv} at $t=0$ and \eqref{eq:V0}
\begin{align}
W^\eps_0  &= - \frac{1}{\varepsilon} L V^\eps_0 + \pa_u \cF_0(U^\eps_0) V^\eps_0
+ (\pa_t \cF)_0(U^\eps_0) \nonumber\\
&= - \frac{1}{\varepsilon} L \Big(\eps A\pa_u \cF_0 \Pi \cF_0 + \eps A (\pa_t \cF)_0\Big) + \pa_u \cF_0(U^\eps_0) V^\eps_0
+ (\pa_t \cF)_0(U^\eps_0) + {\cal O}_{X^{s-4}}(\eps) \nonumber\\
&= -(\I-\Pi) \pa_u \cF_0 \Pi \cF_0 -(\I-\Pi)( \pa_t \cF)_0+ \pa_u \cF_0 \Pi \cF_0 +  (\pa_t \cF)_0 +  {\cal O}_{X^{s-4}}(\eps) \nonumber\\
&= \Pi \pa_u \cF_0 \Pi \cF_0 +  \Pi (\pa_t \cF)_0 +  {\cal O}_{X^{s-4}}(\eps),  \label{eqW0}
\end{align}
which is uniformly bounded in the $H^{s-4}$-norm w.r.t. both $\eps$ and $\tau$. By Proposition \ref{prop:exist} applied to $W^\eps(t,\tau)$ with $\sigma=s-4$, one gets that $W^\eps$ is uniformly bounded in $L^\infty_\tau(H^{s-4})$.

\ms
\ni
{\bf Estimate of the third derivative in $t$.} Finally, we derive the equation for $Y^\eps(t,\tau) = \pa_t W^\eps(t,\tau)$, which reads
\begin{align}
\label{eqy}
\partial_t Y^\eps + \frac{1}{\varepsilon}\partial_\tau Y^\eps
&= \pa^3_u \cF(t,\tau,U^\eps) (V^\eps,V^\eps,V^\eps) + 3 \pa_t \pa^2_u \cF(t,\tau,U^\eps) (V^\eps,V^\eps) \nonumber \\
& +3 \pa^2_u \cF(t,\tau,U^\eps) (V^\eps,W^\eps) + 3 \pa_t^2 \pa_u \cF(t,\tau,U^\eps) V^\eps \nonumber  \\ 
& +3 \pa_t \pa_u \cF(t,\tau,U^\eps) W^\eps + \pa_t^3 \cF(t,\tau,U^\eps) + \pa_u \cF(t,\tau,U^\eps) Y^\eps. 
\end{align}
We then extract $Y^\eps_0(\tau)=\pa_t W^\eps(0,\tau)$ from \eqref{eqw2}:
\begin{equation*}
Y^\eps_0 = - \frac{1}{\varepsilon} L W^\eps_0
+ \pa^2_u \cF_0(U^\eps_0) (V^\eps_0,V^\eps_0) + 2 (\pa_t \pa_u \cF)_0(U^\eps_0) V^\eps_0 \nonumber  + (\pa_t^2 \cF)_0(U^\eps_0) + \pa_u \cF_0(U^\eps_0) W^\eps_0.  
\end{equation*}
The only source of concern could come from the term in factor of $\frac{1}{\eps}$. However, it is clear from the expression \eqref{eqW0} of $W^\eps_0$ that $L W^\eps_0 = {\cal O}_{X^{s-6}}(\eps)$, so that
\begin{equation}
\label{eq:Y0}
Y^\eps_0 = {\cal O}_{X^{s-6}}(1), 
\end{equation}
and according to Proposition   \ref{prop:exist}, $Y^\eps(t,\tau)$ is thus uniformly bounded in $\tau$ and $\eps$ in the $H^{s-6}$-norm, since the {\em source}-term $\pa_t^3 \cF(t,\tau,U^\eps)$ is uniformly bounded  in $H^{s-6}$.

\ms
\ni
{\bf Estimate of derivatives in $\tau$.} Using equation  \eqref{eqF} on $U^\eps$, \eqref{eqv} on $V^\eps$ 
and \eqref{eqw2} on $W^\eps$, we have 
\begin{align*}
\pa_\tau U^\eps &= \eps (-\partial_t U^\eps + \cF(t,\tau,U^\eps)=  \eps (- V^\eps + \cF(t,\tau,U^\eps), \\
\pa_\tau V^\eps &= \eps (-\pa_t V^\eps +  \pa_u \cF(t,\tau,U^\eps) V^\eps ) =  \eps (-W^\eps +  \pa_u \cF(t,\tau,U^\eps) V^\eps ), \\% = {\cal O}_{H^{s-4}}(\eps),
\partial_\tau W^\eps &=  \varepsilon ( -\partial_t W^\eps + \pa^2_u \cF(t,\tau,U^\eps) (V^\eps,V^\eps) \\
&+ 2 \pa_t \pa_u \cF(t,\tau,U^\eps) V^\eps + \pa_t^2 \cF(t,\tau,U^\eps) + \pa_u \cF(t,\tau,U^\eps) W^\eps).  
\end{align*}
Since the terms in the parentheses of the right-hand sides are uniformly bounded in $(t,\tau,\eps)$ respectively in $H^{s-2}$, 
$H^{s-4}$ and $H^{s-6}$ using Assumption  \ref{mainassump}, we deduce that 
$\partial_\tau U^\eps= {\cal O}_{H^{s-2}}(\eps)$, $\partial_\tau V^\eps= {\cal O}_{H^{s-4}}(\eps)$, 
and $\partial_\tau W^\eps= {\cal O}_{H^{s-6}}(\eps)$. 
%A similar argument using equation \eqref{eqw2} allows to assert that $\pa_\tau W^\eps = \pa_t^2 \pa_\tau U^\eps =  {\cal O}_{H^{s-6}}(\eps)$ is uniformly bounded as well. 
The estimate of $\pa_\tau Y^\eps = \pa_t^3 \pa_\tau U^\eps$ from \eqref{eqy} however requires some additional argument since $Z^\eps=\pa_t Y^\eps$ is not necessarily uniformly bounded. We thus consider the equation satisfied by $\tilde Z^\eps = \pa_\tau Y^\eps$ obtained by differentiation w.r.t. to $\tau$ of equation \eqref{eqy}. This equation is of the form 
$$
\pa_t \tilde Z^\eps + \frac{1}{\eps} \pa_\tau \tilde Z^\eps = S(t,\tau,U^\eps,V^\eps,W^\eps,Y^\eps,\pa_\tau U^\eps,\pa_\tau V^\eps,\pa_\tau W^\eps)+\pa_u \cF(t,\tau,U^\eps) \tilde Z^\eps, 
$$
and its initial condition can be obtained by differentiating $Y_0^\eps$ w.r.t. $\tau$, leading to $\tilde Z^\eps_0 = {\cal O}_{X^{s-8}}(1)$ and hence, once again by Proposition  \ref{prop:exist}, $\tilde Z^\eps$ is uniformly bounded in $L^\infty_\tau(H^{s-8})$, since the source-term $S$ lies in $C^0([0,T_\kappa]\times \T;H^{s-8})$. 
\end{Proof}
\begin{remark}\label{rempropcondinit}
If in Proposition \pref{propcondinit} we modify the hypotheses as follows: $s > d/2+4$ and we assume the following initial condition
\be
\label{condinit0bis}
\forall \tau \in \T,\qquad  U^\eps(0, \tau) = U_0^\eps(\tau)=\ub U^\eps_0 + \eps h_1(0,\tau, \ub U^\eps_0)+\eps^2r^\eps(\tau), 
\ee
where $\ub U^\eps_0 \in H^s$ is uniformly bounded in $\eps$, where $h_1$ is given by \eqref{h1}, and where the remainder term $r^\eps$ is assumed to be bounded in $X^s$ uniformly in $\eps$, then the conclusions of this Proposition are modified as follows: Item (i) is unchanged, and (ii) is replaced by:
\begin{enumerate}
\item[(ii)] For any $T_\kappa$ for which \eqref{estiunif} holds, the solution $U^\eps$ satisfies the following estimates
\begin{equation}
\label{unifestiderivbis}
\forall t\in [0,T_\kappa],\qquad \sup_{\eps\in ]0,\eps_0]}\|\partial_t^\alpha U^\eps(t)\|_{L^{\infty}_\tau(H^{s-2\alpha})} \leq C,  \qquad \alpha=0,1,2, 
\end{equation}
for some constant $C>0$. 
\end{enumerate}
\end{remark}

\section{First order scheme} \label{sect:fo}

This section is devoted to the construction and the numerical analysis of a first order numerical scheme.  We only  focus on the time discretization, and the other variable $\tau$ is kept at the continuous level.  We thus define a  uniform grid $t_n= n\Delta t$ of a time interval $[0,T_\kappa]$, $n=0,1,..., N,$ with $N\Delta t =T_\kappa$ (recall that $[0,T_\kappa]$ is the interval of time on which the solution $U^\eps$ of \eqref{eqF} is well-defined irrespectively  of $\eps$)   and consider the following numerical scheme for this equation. Denoting $U^\eps_{n} \approx U^\eps(t_{n},\tau)$ and $U^\eps_{n+1} \approx U^\eps(t_{n+1},\tau)$, it advances the solution from time $t_n$ to time $t_{n+1}$ through the inductive equation
\begin{equation}
\label{eq:sch1}
U^\eps_{n+1}(\tau) = U^\eps_n(\tau) + \Delta t \cF(t_n,\tau,U^\eps_n(\tau)) - \frac{\Delta t}{\eps} \pa_\tau U^\eps_{n+1}(\tau).
\end{equation}
This scheme does not define $U_{n+1}^\eps$ in a unique way, unless we impose --\,assuming that $U_n^\eps$ is periodic with period $P$\,-- that $U_{n+1}^\eps$ is periodic with the same period. Under this requirement, pre-multiplying \eqref{eq:sch1} by $e^{\mu \tau}$ with $\mu = \frac{\eps}{\Delta t}$, we have 
$$
\pa_\tau(e^{\mu \tau} U^\eps_{n+1})= \mu e^{\mu \tau} (U^\eps_{n+1} + \frac{1}{\mu} \pa_\tau U^\eps_{n+1})= \mu e^{\mu \tau}(U^\eps_{n} + \Delta t \cF(t_n,\tau,U^\eps_{n}))
$$
so that, upon integrating from $\tau$ to $\tau+P$, we obtain 
\begin{equation*}
(e^{\mu P}-1) e^{\mu \tau} U^\eps_{n+1}(\tau)= \mu \int_\tau^{\tau+P}  e^{\mu \theta} \Big(U^\eps_n(\theta) + \Delta t \cF(t_n,\theta,U^\eps_n(\theta))\Big)  d\theta,
\end{equation*}
or in more concise manner
\begin{equation}
\label{eq:sch1nv}
U^\eps_{n+1}(\tau)= \frac{\mu}{\exp(\mu P)-1}\int_\tau^{\tau+P}  e^{\mu (\theta-\tau)} \Big(U^\eps_n(\theta) + \Delta t \cF(t_n,\theta,U^\eps_n(\theta))\Big)  d\theta. 
\end{equation}
Note that it is straightforward to check that $U_{n+1}^\eps$, as given by formula \eqref{eq:sch1nv}, is periodic of period $P$. This last equation defines precisely the scheme whose convergence we now wish to investigate. From previous computations, we observe that the operator $Q_\mu=\I+\frac{1}{\mu} \pa_\tau$ play a central role. We thus briefly study its properties in the following proposition.
\begin{proposition} \label{Qmu}
Let $\sigma\in \R$. The operator $Q_\mu$, defined from the set $C^1(\T;H^\sigma)$ onto $C^0(\T;H^\sigma)$ by 
$$
\forall \tau \in \T, \quad (Q_\mu g)(\tau) = g(\tau) + \frac{1}{\mu} (\pa_\tau g)(\tau)
$$
is invertible and its inverse, which is defined on $C^0(\T;H^\sigma)$ with values in $C^1(\T;H^\sigma)$, can be explicitly written as 
$$
\forall \tau \in \T, \quad (Q_\mu^{-1} g)(\tau) = \frac{\mu}{\exp(\mu P)-1}\int_{\tau}^{\tau+P} e^{\mu(\theta-\tau)} g(\theta) d\theta.
$$
Moreover, it satisfies the following estimate:
$$
\forall g \in C^0(\T;H^\sigma),  \qquad \|Q_\mu^{-1} g\|_{L^\infty_\tau(H^\sigma)} \leq  \| g\|_{L^\infty_\tau(H^\sigma)}
$$
\end{proposition}
\begin{Proof}
The inversion formula has been proven above. As for the estimate, it stems from the identity
$$
 \frac{\mu}{\exp(\mu P)-1} \int_{\tau}^{\tau+P} e^{\mu(\theta-\tau)}  d\theta = 1. 
$$
\end{Proof}
\subsection{Local truncation error}
\begin{proposition} \label{prop:lerr} Let $s > d/2+4$, $\kappa >1$ and $\ub U_0^\eps$ uniformly bounded in $H^s$ with respect to $\eps$. Consider $U^\eps$ the solution of equation \eqref{eqF} with initial condition \eqref{condinit0bis}, and fix $\Delta t = \frac{T_\kappa}{N}$ for some $N \in \N^*$. 
The consistency error  $\elle^\eps_{n+1}$ for the numerical scheme (\ref{eq:sch1}) or (\ref{eq:sch1nv}), 
defined for $n=0,\ldots,N-1$ by 
\begin{equation}
\label{eq:cons}
\elle^\eps_{n+1} := U^\eps(t_{n+1}) - \tilde U^\eps_{n+1}
\end{equation}
where 
$$
(Q_\mu \tilde U^\eps_{n+1})(\cdot) = U^\eps(t_n,\cdot) + \Delta t \cF(t_n,\cdot,U^\eps(t_n,\cdot))
$$
satisfies the following estimate 
\begin{equation} \label{eq:lest1}
\sup_{\eps \in ]0,\eps_0]} \|\elle^\eps_{n+1} \|_{L^\infty_{\tau}(H^{s-4})} \leq  C \Delta t^2
\end{equation}
for some positive error constant $C$. 
\end{proposition}
\begin{Proof}
In the sequel, we shall omit the variable $\tau$ unless explicitly needed,  as it plays essentially no role in subsequent computations. From the equation satisfied by $U^\eps(t,\tau)$ (see \eqref{eqF}), we have
$$
Q_\mu U^\eps(t_{n+1}) = U^\eps(t_{n+1}) + \Delta t \Big( \cF(t_{n+1},U^\eps(t_{n+1})) - \pa_t U^\eps(t_{n+1})\Big), 
$$
so that 
\begin{equation} 
\label{eq:gnu}
Q_\mu (\elle^\eps_{n+1}) =g_n, 
\end{equation}
with 
$$
g_n := U^\eps(t_{n+1}) - U^\eps(t_n) + \Delta t \Big( \cF(t_{n+1},U^\eps(t_{n+1})) - \cF(t_n,U^\eps(t_n)) - \pa_t U^\eps(t_{n+1})\Big). 
$$
By Taylor expansion (with integral remainder) at $t=t_{n+1}$  we get, on the one hand,
$$
U^\eps(t_{n}) - U^\eps(t_{n+1}) = (-\Delta t)  \partial_t U^\eps(t_{n+1}) + \int_{t_{n+1}}^{t_{n}} (t_n-t) \partial^2_{t} U^\eps(t)dt,
$$
and on the other hand, 
\begin{equation*}
\cF(t_{n+1},U^\eps(t_{n+1})) = \cF(t_n, U^\eps(t_{n})) + \int_{t_n}^{t_{n+1}} \Big(\pa_t \cF(t,U^\eps(t)) + \pa_u \cF(t, U^\eps(t)) \partial_t U^\eps(t)  \Big) dt, 
\end{equation*}
so that 
$$
g_n = \Delta t \int_{t_n}^{t_{n+1}}  \Big(\pa_t \cF(t,U^\eps(t)) + \pa_u \cF(t,U^\eps(t)) \partial_t U^\eps(t)+  \frac{(t_n-t)}{\Delta t} \partial^2_{t} U^\eps(t)   \Big)dt.
$$
According to Remark \ref{rempropcondinit} and the Assumption \ref{mainassump} on $\cF$,  the quantities  $\sup_{t \in [0,T_\kappa]} \|\partial_t^\alpha U^\eps\|_{L^\infty_{\tau} (H^{s-2\alpha})}, \alpha=0,1,2$ are uniformly bounded in  $\eps$, so that 
\be
\label{bornebn}
\sup_{\eps \in ]0,\eps_0]} \|g_{n}\|_{L^\infty_\tau (H^{s-4})} \leq C\Delta t ^2,
\ee
and by Proposition \ref{Qmu}, we get the desired estimate. 
\end{Proof}
\subsection{Global  error}
%
%
% Main theorem
%
\begin{theorem}
\label{theo-o1}
Let $s > d/2+4$, $\kappa >1$ and $\ub U_0^\eps$ uniformly bounded in $H^s$ with respect to $\eps$. Let $U^\eps$ be the unique solution of \eqref{eqF} on $[0, T_\kappa]$ subject to the initial condition \eqref{condinit0bis}, and let $(U^\eps_n)_{0 \leq n \leq N}$ be defined for all $\tau\in\T$ by \eqref{eq:sch1nv}
with $U^\eps_0$ again defined by \eqref{condinit0bis}.  
%Under the additional  assumption
%\begin{eqnarray} \label{eq:boundderiv}
%\sup_{I \times \T \times K_{s-4}} \|(\cF')(t,\tau,u)\|_{{\cal L}_{1,s-4}} \leq B
%\end{eqnarray}
Then there exists $\Delta t_0 >0$ and $C>0$ such that the following estimate holds 
\begin{equation} 
\label{eq:est1}
\forall \, \Delta t < \Delta t_0, \qquad 
\sup_{\eps\in ]0,\eps_0]} \|U^\eps(t_n)-U^\eps_n\|_{L^\infty_{\tau}(H^{s-4})} \leq C\Delta t, 
\end{equation}
 for all $n=0, \dots, N$, where $T_\kappa=N\Delta t$ is the final time. 
\end{theorem}
\begin{Proof}
The global error $e^\eps_{n+1} := U^\eps(t_{n+1})-U^\eps_{n+1}$ can be decomposed into two parts
$$
e_{n+1}^\eps =  \underbrace{\left(U^\eps(t_{n+1}) -  \tilde U_{n+1}^\eps\right)}_{\mbox{local  error $\elle^\eps_{n+1}$}}+\underbrace{\left( \tilde U_{n+1}^\eps- U_{n+1}^\eps\right)}_{\mbox{transported error $\te^\eps_{n+1}$}}
$$
where 
\begin{equation*}
Q_{\mu} \tilde U^\eps_{n}= U^\eps(t_n)+ \Delta t \cF(t_n,U^\eps(t_n)).
\end{equation*}
We have 
$$
Q_{\mu} e_{n+1}^\eps = Q_{\mu} \elle^\eps_{n+1} + Q_{\mu} \te^\eps_{n+1}
$$
with (see equation \eqref{eq:gnu})
$$
Q_{\mu} \elle^\eps_{n+1} = g_n
$$
and 
$$
Q_{\mu} \te^\eps_{n+1} = e^\eps_n + \Delta t \Big( \cF(t_n,U^\eps(t_n))- \cF(t_n,U^\eps_n) \Big).  
$$
Let $R>0$ be fixed. As long as $\sup_{\eps\in ]0,\eps_0]} \|e^\eps_n\|_{L^\infty_\tau (H^{s-4})} < R$ (recall that $e^\eps_0=0$), we have for $\lambda \in [0,1]$
 $$\|\lambda U^\eps_n+(1-\lambda) U^\eps(t_n)\|_{L^\infty_\tau (H^{s-4})}\leq \|U^\eps (t_n)\|_{L^\infty_\tau (H^{s-4})}+\lambda \|e^\eps_n\|_{L^\infty_\tau (H^{s-4})} \leq C_0+ R$$
independently of $\eps$, where we used the uniform bound on $U^\eps$ given in Remark \ref{rempropcondinit}. Then we can use the local bound of $\pa_u \cF$ provided by Assumption \ref{mainassump} in order to obtain 
\begin{align*}
 \|\cF(t_n,U^\eps_n)- \cF(t_n, U^\eps(t_n))\|_{L^\infty_\tau (H^{s-4})} &\leq \int_0^1  \Big\|\pa_u \cF(t_n,\lambda U^\eps_n+(1-\lambda) U^\eps(t_n)) e^\eps_n  \Big\|_{L^\infty_\tau (H^{s-4})} d\lambda \\
&\leq  C_1 \|e^\eps_n\|_{L^\infty_\tau (H^{s-4})}.
\end{align*}
Finally, using estimate \eqref{eq:lest1} of Proposition \ref{prop:lerr} and Proposition \ref{Qmu}, we have 
$$
\|e^\eps_{n+1}\|_{L^\infty_\tau(H^{s-4})} \leq (1+ C_1 \Delta t) \|e^\eps_{n}\|_{L^\infty_\tau(H^{s-4})}  +   C_2 \Delta t^2 
$$
and a discrete Gronwall lemma provides us with the bound 
$$
\|e^\eps_{n+1}\|_{L^\infty_{\tau}(H^{s-4})} \leq \frac{C_2 \Delta t}{C_1} (\exp(C_1 T_\kappa)-1).
$$
By induction, we can a posteriori verify that $\|e^\eps_n\|_{H^{s-4}} < R$ for all $n=0,\cdots, N$, provided $\Delta t < \Delta t_0$
with $\Delta t_0 := R \, C_1 \, e^{-C_1 T_\kappa} \, /C_2.$ This completes the proof. 
\end{Proof}
\subsection{Choice of the initial data for the first order scheme}
\label{subfirst}
So far, we have addressed the question of numerical approximation of the augmented problem \eqref{eqF}, subject to an initial condition $U_0^\eps(\tau)$ satisfying \eqref{condinit0bis}. We now come back to our original problem \eqref{eqftildegene}, and recall that if $U^\eps_0(0)=u_0$, then $u^\eps(t)$ can be recovered as the diagonal $u^\eps(t)=U^\eps(t,t/\eps)$, so the global estimate \eqref{eq:est1} yields
\begin{equation*} 
\forall \, \Delta t < \Delta t_0, \qquad 
\sup_{\eps\in ]0,\eps_0]} \|u^\eps(t_n)-U^\eps_n(t_n/\eps)\|_{H^{s-4}} \leq C\Delta t. 
\end{equation*}
In practice, if the only known initial data is $u_0$, we proceed as follows to construct a suitable associated initial data $U_0^\eps$. We set
\begin{equation}
\label{firstID}
U_0^\eps(\tau):=u_0+\eps h_1(0,\tau,u_0)-\eps h_1(0,0,u_0)=u_0 + \eps \int_0^\tau (\I-\Pi) \cF_0(\theta,u_0) d\theta,
\end{equation}
and
$$
\ub U_0^\eps:=\Pi U_0^\eps=u_0 + \eps \Pi \int_0^\tau (\I-\Pi) \cF_0(\theta,u_0) d\theta.$$
In the numerical  Section \ref{sect:num}, this choice of initial data will be referred to as {\em first order initial data}.
Then, \eqref{condinit0bis} is satisfied if the remainder term is defined by
$$r^\eps(\tau):=\frac{1}{\eps^2}\left(U_0^\eps(\tau)-\ub U^\eps_0 - \eps h_1(0,\tau, \ub U^\eps_0)\right)=\frac{1}{\eps}\left(h_1(0,\tau,u_0)-h_1(0,\tau,\ub U_0^\eps)\right).$$
From Assumption \ref{mainassump}, it is easy to prove that, as soon as $u_0\in H^s$, the function $\ub U^\eps_0$ and $r^\eps$ are uniformly bounded respectively in $H^s$ and $X^s$ (this space is defined by \eqref{defX}), which is enough to apply Theorem \ref{theo-o1}.

%%%%%%%%%%%%
%
% Second-oder scheme
%
%%%%%%%%%%%%
\section{A second order scheme} \label{sect:so}
We now present a two-stage second order scheme. The first stage is  composed of half-a-step of the first-order scheme presented in previous section:
\be
\label{s1}
U^\eps_{n+1/2} = U^\eps_n +\frac{ \Delta t}{2} \cF(t_n,U^\eps_n) - \frac{\Delta t}{2\varepsilon}\partial_\tau U^\eps_{n+1/2}, 
\ee
while the second stage computes the updated approximation $U^\eps_{n+1}$ as follows: 
\be
\label{s2}
U^\eps_{n+1}  = U^\eps_n + \Delta t \cF(t_{n+1/2}, U^\eps_{n+1/2}) - \frac{\Delta t}{2\varepsilon}\partial_\tau (U^\eps_n+U^\eps_{n+1}).  
\ee
As in the previous section, a more concise version of the scheme reads
\begin{align}
\label{eq:sch2nv}
Q_{2 \mu} U^\eps_{n+1/2} &= U^\eps_n +\frac{ \Delta t}{2} \cF(t_n,U^\eps_n), \\
Q_{2 \mu} U^\eps_{n+1}  &= (2 \I- Q_{2 \mu}) U^\eps_{n}  + \Delta t \cF(t_{n+1/2}, U^\eps_{n+1/2}). \label{eq:sch2nv2}
\end{align}
%Note that since $\pa_\tau$ and $Q_\mu$ commute, the scheme intrinsically defines an approximation of the derivative w.r.t. $\tau$ through the formulae
%\begin{eqnarray}
%\label{eq:sch2nvtau}
%Q_{2 \mu} V^\eps_{n+1/2} &=& V^\eps_n +\frac{ \Delta t}{2} \Big(\pa_\tau \cF(t_n,U^\eps_n) + \pa_u \cF(t_n,U^\eps_n) V^\eps_n\Big), \\
%Q_{2 \mu} V^\eps_{n+1}  &=& (2 I- Q_{2 \mu}) V^\eps_{n}  + \Delta t \Big(\pa_\tau \cF(t_{n+1/2}, U^\eps_{n+1/2})+\pa_u \cF(t_{n+1/2}, U^\eps_{n+1/2}) V^\eps_{n+1/2} \Big) \nonumber
%\end{eqnarray}
%where $V^\eps_n = \pa_\tau U^\eps_n$ can be regarded as an approximation of the solution of equation
%\begin{eqnarray} \label{eq:tev}
%\partial_t V^\eps  + \frac{1}{\varepsilon}\pa_\tau V^\eps = 
%\pa_\tau \cF(t,\tau,U^\eps) + \pa_u \cF(t,\tau,U^\eps) V^\eps, \qquad V^\eps(0) = \pa_\tau U^\eps_0.
%\end{eqnarray}
%Finally, let us remark that the scheme intrinsically defines an approximation of the derivative w.r.t. $\tau$ through the formula
%\begin{eqnarray} \label{eq:sch2tau}
%\frac{1}{2 \mu} \pa_\tau U^\eps_{n+1}(\tau) &=& \varphi^P_{2\mu} \int_\tau^{\tau+P}  e^{2 \mu (s-\tau)} \Big(U^\eps_n(\tau) - U^\eps_n(s) ) ds\\
%&+&\varphi^P_{2\mu} \int_\tau^{\tau+P}  e^{2 \mu (s-\tau)} \Delta t  \Big(\cF(t_{n+1/2},\tau,U^\eps_{n+1/2}(\tau))-\cF(t_{n+1/2},s,U^\eps_{n+1/2}(s)) \Big) ds  \nonumber \\
%&-& \varphi^P_{2\mu} \int_\tau^{\tau+P}  e^{2 \mu (s-\tau)} \frac{1}{2 \mu}  \Big(\pa_\tau U^\eps_{n}(\tau)-\pa_\tau U^\eps_{n}(s)\Big) ds \nonumber
%\end{eqnarray}

%%%%%%%%%%%%%%%%%%%%%%%%%%%%%%%%%
% A preliminary lemmma
%%%%%%%%%%%%%%%%%%%%%%%%%%%%%%%%%
Prior to proving our main result, we state an elementary result which is an essential ingredient of subsequent proofs.
\begin{lemma} \label{lem:prelim2}
Consider a function $g \in L^2(\T ;H^\sigma)$ with $\sigma\in \R$. The following estimates hold true for the operator $Q_{\mu}$ defined in Proposition \pref{Qmu}:
\begin{equation} \label{eq:pidi}
\forall \beta \in [-1,1],  \qquad \|((1+\beta) Q_{\mu}^{-1} - \beta \I) g \|_{L^2_\tau(H^\sigma)} \leq  \|g\|_{L^2_\tau(H^\sigma)}, 
\end{equation}
with equality for $|\beta|=1$, and
\begin{equation} \label{eq:pidi2}
\|\pa_\tau Q_{\mu}^{-1}g \|_{L^2_\tau(H^\sigma)} \leq  2\mu\|g\|_{L^2_\tau(H^\sigma)}. 
\end{equation}
\end{lemma}
\begin{Proof}
Let $g \in L^2(\T ;H^\sigma)$ and $((1+\beta) Q_{\mu}^{-1} - \beta \I) g=f$. This last equality is clearly equivalent to
$$g-f = \frac{1}{\mu} \pa_{\tau} (\beta g +f).$$
Taking the inner product against $\beta g +f$ in the real Hilbert space $L^2_\tau(H^\sigma)$ and using the skew-symmetry of $\pa_\tau$, 
we get
$$ \|f\|_{L^2_\tau(H^\sigma)}^2- \beta \|g\|_{L^2_\tau(H^\sigma)}^2= (1-\beta) \langle f, g\rangle_{L^2_\tau(H^\sigma)} \leq (1-\beta)\|f\|_{L^2_\tau(H^\sigma)} \|g\|_{L^2_\tau(H^\sigma)} .$$
This proves the case of equality in \eqref{eq:pidi} and also implies  
$$ (\|f\|_{L^2_\tau(H^\sigma)} + \beta \|g\|_{L^2_\tau(H^\sigma)})(\|f\|_{L^2_\tau(H^\sigma)}- \|g\|_{L^2_\tau(H^\sigma)})\leq 0,$$
which ends the proof of inequality \eqref{eq:pidi}. To prove \eqref{eq:pidi2}, we simply remark that
$$\pa_\tau Q_{\mu}^{-1} g=\mu(g-Q_{\mu}^{-1}g),$$
and use \eqref{eq:pidi} with $\beta=0$.
\end{Proof}
\subsection{Local truncation error}
\begin{proposition}\label{prop:soest} Let $s > d/2+8$, $\kappa >1$ and $\ub U_0^\eps$ uniformly bounded in $H^{s+2}$ with respect to $\eps$. Consider $U^\eps$ the solution of equation \eqref{eqF} with initial condition \eqref{condinit0}  given by Proposition \pref{propcondinit}, and fix $\Delta t = \frac{T_\kappa}{N}$ for some $N \in \N^*$. 
The local truncation error $\elle_{n+1}$ for the scheme (\ref{s1})-(\ref{s2}) 
defined for $n=0,\ldots,N-1$ by $\elle^\eps_{n+1} := U^\eps(t_{n+1}) - \tilde U^\eps_{n+1}$ where
\begin{align*}
Q_{2 \mu} \tilde U^\eps_{n+1/2} &= U^\eps(t_n) +\frac{ \Delta t}{2} \cF(t_n, U^\eps(t_n)),\\
Q_{2 \mu} \tilde U^\eps_{n+1} &= (2 \I - Q_{2 \mu}) U^\eps(t_n) + \Delta t \cF(t_{n+1/2}, \tilde U^\eps_{n+1/2}), 
\end{align*}
satisfies the following estimates
\begin{equation} \label{eq:lest}
\sup_{\eps \in ]0,\eps_0]} \|\elle^\eps_{n+1} \|_{L^\infty_\tau(H^{s-4-2k})} \leq  C (\Delta t)^{2+k},
\end{equation}
and
\begin{equation} \label{eq:LEst}
\sup_{\eps \in ]0,\eps_0]} \|\pa_\tau \elle^\eps_{n+1} \|_{L^\infty_\tau(H^{s-6-2k})} \leq  C (\Delta t)^{2+k},
\end{equation}
for $k\in \{0,1\}$ and for some positive error constant $C$.
\end{proposition}
\begin{Proof} 
From the equation satisfied by $U^\eps(t,\tau)$ at $t_n$ and $t_{n+1}$ we have
\begin{align*}
(2\I-Q_{2\mu}) U^\eps(t_{n}) &= U^\eps(t_{n}) - \frac{\Delta t}{2}  \Big( \cF(t_{n},U^\eps(t_{n})) - \pa_t U^\eps(t_{n})\Big), \\ 
Q_{2\mu} U^\eps(t_{n+1}) &= U^\eps(t_{n+1}) + \frac{\Delta t}{2}  \Big( \cF(t_{n+1},U^\eps(t_{n+1})) - \pa_t U^\eps(t_{n+1})\Big)
\end{align*}
so that 
%NC
\begin{equation*} %\label{eq:gnu}
(Q_{2 \mu} \elle^\eps_{n+1}) =g_n
\end{equation*}
with 
%$$
%g_n := U^\eps(t_{n+1}) -  U^\eps(t_n) + \frac{\Delta t}{2} \Big( \cF(t_{n+1},U^\eps(t_{n+1})) +  \cF(t_{n+1/2},U^\eps(t_{n+1/2})) - \pa_t U^\eps(t_{n+1})- \pa_t U^\eps(t_{n+1})\Big)
%$$
%NC
\begin{align}
g_n &:= U^\eps(t_{n+1}) -  U^\eps(t_n) + \frac{\Delta t}{2} \Big( \cF(t_{n+1},U^\eps(t_{n+1})) + \cF(t_{n},U^\eps(t_{n})) -2 \cF(t_{n+1/2},\tilde{U}^\eps_{n+1/2})\Big) \nonumber\\
&\qquad - \frac{\Delta t}{2}  \Big( \pa_t U^\eps(t_{n+1}) + \pa_t U^\eps(t_{n})\Big). \label{gn}
\end{align}
Now, by a symmetric Taylor expansion, it is straightforward to show that 
$$
U^\eps(t_{n+1}) - U^\eps(t_n) = \frac{\Delta t}{2} \Big(\pa_t U^\eps(t_{n+1}) +\pa_t U^\eps(t_n) \Big)
+  R_n,  
$$
with 
$$
R_n:= \frac{1}{2}  \int_{t_n}^{t_{n+1}} (t_n-t)(t_{n+1}-t) \pa_t^3 U^\eps(t) dt.
$$
This remainder term can be estimated in two different ways.
First, using the uniform bound of $\|\pa_t^3 U^\eps\|_{L^\infty_\tau(H^{s-6})}$ in \eqref{unifestideriv}, we get directly
$$
\|R_n\|_{L^\infty_{\tau}(H^{s-6})} \leq C \Delta t^3.
$$
Second, integrating by parts $R_n= \frac{1}{2}  \int_{t_n}^{t_{n+1}} (t_n+t_{n+1}-2t) \pa_t^2 U^\eps(t) dt,$
and using the uniform bound of $\|\pa_t^2 U^\eps\|_{L^\infty_\tau(H^{s-4})}$ in  \eqref{unifestideriv}, we obtain
$$\|R_n\|_{L^\infty_\tau(H^{s-4})}\leq C\Delta t^2.$$
Here and in the sequel, $C$ denotes a generic positive constant independent of $\eps\in ]0,\eps_0]$. 

Next, from the proof of Proposition \ref{prop:lerr}, we have the estimate
$$
\|\tilde U^\eps_{n+1/2}-U^\eps(t_{n+1/2})\|_{L^\infty_\tau(H^{s-4})} \leq C \Delta t^2
$$
from which we get
$$
\|\cF(t_{n+1/2}, U^\eps(t_{n+1/2})) - \cF(t_{n+1/2},  \tilde{U}^\eps_{n+1/2}) \|_{L^\infty_\tau(H^{s-4})} \leq C \, \Delta t^2,
$$
where we have used the local boundedness of $\pa_u \cF$ and the uniform boundedness of $U^\eps$.
Using once more a Taylor expansion (of the function $t \mapsto \cF(t,U^\eps(t))$ around $t=t_{n+1/2}$), it stems from the local boundedness of the  first and second derivatives of $\cF$,  and the uniform boundedness of the first and second time derivatives of $U^\eps$,
$$
\|\cF(t_{n+1},U^\eps(t_{n+1})) + \cF(t_{n},U^\eps(t_{n})) - 2 \cF(t_{n+1/2}, U^\eps(t_{n+1/2}))\|_{L^\infty_\tau(H^{s-4})} \leq C\Delta t^2.
$$
This eventually leads to 
$$
\sup_{\eps\in ]0,\eps_0]}\|g_n\|_{L^\infty_\tau(H^{s-4-2k})} \leq C \Delta t^3+\sup_{\eps\in ]0,\eps_0]}\|R_n\|_{L^\infty_\tau(H^{s-4-2k})}\leq C(\Delta t)^{2+k}
$$
for $k\in\{0,1\}$.
The estimate \eqref{eq:lest} then follows from the properties of the inversion formula for $Q_{2 \mu}$. 

The estimate \eqref{eq:LEst} can be obtained in a similar way. Denoting 
$$\cG(t,\tau,U^\eps,V^\eps) = \pa_\tau \cF(t,\tau,U^\eps)+ \pa_u \cF(t,\tau,U^\eps) V^\eps,$$
we see that the derivative w.r.t. $\tau$ of the scheme \eqref{eq:sch2nv}, \eqref{eq:sch2nv2} yields the following scheme on the unknown $V^\eps=\pa_\tau U^\eps$:
\begin{align*}
Q_{2 \mu} V^\eps_{n+1/2} &= V^\eps_n +\frac{ \Delta t}{2} \cG(t_n,U^\eps_n,V^\eps_n), \\
Q_{2 \mu} V^\eps_{n+1}  &= (2 \I- Q_{2 \mu}) V^\eps_{n}  + \Delta t \cG(t_{n+1/2}, U^\eps_{n+1/2},V^\eps_{n+1/2}).
\end{align*}
The only  difference is now that the r.h.s. $\cG$ is a now a map from $\R_+ \times \T \times H^{s-2}$ into $H^{s-2}$ which has, according to Proposition \ref{propcondinit}, uniformly  bounded derivatives $\pa_t^\alpha \cG$ for the $L^\infty_{t,\tau}(H^{s-2(\alpha+1)})$ norm for $\alpha=0,1,2,3$.
\end{Proof}

\subsection{Global error}
%%%%%%%%%%%%%%%%%%%%%%%%%%%%%%%%%
% Convergence second-order scheme 
%%%%%%%%%%%%%%%%%%%%%%%%%%%%%%%%%
\begin{theorem}
\label{theo-o2}
Let $s > d/2+8$, $\kappa >1$ and $\ub U_0^\eps$ uniformly bounded in $H^{s+2}$ with respect to $\eps$. Let $U^\eps$ be the unique solution of \eqref{eqF} on $[0, T_\kappa]$ subject to the initial condition \eqref{condinit0} given by Proposition \pref{propcondinit}, and let $(U^\eps_n)_{0 \leq n \leq N}$ be defined for all $\tau\in\T$ by \eqref{eq:sch2nv}
with $U^\eps_0$ again defined by \eqref{condinit0}.  Then there exists $\Delta t_0 >0$ and $C>0$ such that the following estimate holds 
\begin{equation} 
\label{eq:est2}
%\sup_{\eps\in ]0, \eps_0]} \|U^\eps(t_n)-U^\eps_n\|_{L^\infty_{\tau}(H^{s-8})} \leq C\Delta t^2, 
\forall \, \Delta t < \Delta t_0, \qquad \sup_{\eps\in ]0,\eps_0]} \|U^\eps(t_n)-U^\eps_n\|_{L^\infty_{\tau}(H^{s-8})} \leq C\Delta t^2, 
\end{equation}
 for all $n=0, \dots, N$, where $T_\kappa=N\Delta t$ is the final time. 
\end{theorem}

\begin{Proof} Let $R>0$ and let us assume for the time being, that $e^\eps_n$ satisfies the uniform estimate
\begin{equation}
\label{boots}
\sup_{\eps\in ]0,\eps_0]}\|e^\eps_n\|_{L^\infty_\tau(H^{s-6})}\leq R,
\end{equation}
for all $n\leq N$, so that all terms like $\pa_u \cF(t,\cdot,\lambda U^\eps(t_n,\cdot) + (1-\lambda) U^\eps_n(\cdot))$ are also uniformly bounded from Assumption \ref{mainassump} and Proposition \ref{propcondinit}. This hypothesis can eventually be justified by induction as in Theorem \ref{theo-o1}, using \eqref{first} that we will obtained in the midst of the proof. 

The global error $e_{n+1}^\eps=U^\eps(t_{n+1})-U_{n+1}^\eps$ can be decomposed into two parts as
$$
e_{n+1}^\eps =  \underbrace{\left(U^\eps(t_{n+1}) -  \tilde U_{n+1}^\eps \right)}_{\mbox{local  error $\elle^\eps_{n+1}$}}+\underbrace{\left( \tilde U_{n+1}^\eps - U_{n+1}^\eps \right),}_{\mbox{transported error $\te^\eps_{n+1}$}}
$$
where 
\begin{align*}
Q_{2 \mu} \tilde U^\eps_{n+1/2} &= U^\eps(t_n) + \frac{\Delta t}{2} \cF(t_n,U^\eps(t_n)), \\
Q_{2 \mu}  \tilde U^\eps_{n+1} &= (2 \I -  Q_{2 \mu}) U^\eps(t_n)  + \Delta t \cF(t_{n+1/2},\tilde U^\eps_{n+1/2}).
\end{align*}
We thus have 
$$
Q_{2 \mu} e_{n+1}^\eps = Q_{2 \mu} \elle^\eps_{n+1} + Q_{2 \mu} \te^\eps_{n+1}
$$
with 
$$
Q_{2 \mu} \elle^\eps_{n+1} = g_n\,,
$$
$g_n$ being still defined by \eqref{gn}, and 
\begin{equation}
\label{Qtransport}
Q_{2 \mu} \te^\eps_{n+1} = \Big(2 \I - Q_{2 \mu}\Big)  e^\eps_{n} + \Delta t \Big(\cF(t_{n+1/2},\tilde{U}^\eps_{n+1/2})-\cF(t_{n+1/2},U^\eps_{n+1/2})\Big). 
\end{equation}

\ms
\ni
{\bf Step 1: second order estimate of the global error in $L^2_\tau(H^{s-6})$.} In this first step, we prove the following estimate on the global error:
\begin{equation}
\label{first0}
\forall \, \Delta t < \Delta t_0, \qquad \sup_{\eps\in ]0,\eps_0]} \|e_n^\eps\|_{L^2_{\tau}(H^{s-6})} \leq C\Delta t^2. 
\end{equation}
The second term of the r.h.s. of \eqref{Qtransport} can be bounded as follows
\begin{align*}
&\Big\|\cF(t_{n+1/2},\tilde{U}^\eps_{n+1/2})-\cF(t_{n+1/2},U^\eps_{n+1/2})\Big\|_{L^2_\tau(H^{s-6})}\\
&\qquad = \Big\|\int_{0}^{1}  \pa_u \cF\left(t_{n+1/2},\lambda \tilde U^\eps_{n+1/2} + (1-\lambda)   U^\eps_{n+1/2} \right) \Big(\tilde U_{n+1/2}^\eps -U_{n+1/2}^\eps \Big) d \lambda\Big\|_{L^2_\tau(H^{s-6})} \\
&\qquad \leq C_1 \|\tilde e^\eps_{n+1/2}\|_{L^2_\tau(H^{s-6})} 
\end{align*} 
where $\tilde e^\eps_{n+1/2} = \tilde U_{n+1/2}^\eps -U_{n+1/2}^\eps$ can also be bounded by using the properties of $Q_{2 \mu}^{-1}$ (Lemma \ref{lem:prelim2} for $\beta=0$) 
\begin{equation}
\label{enp12}
\|\tilde e^\eps_{n+1/2}\|_{L^2_\tau(H^{s-6})} \leq (1+C_2\Delta t) \|e^\eps_{n}\|_{L^2_\tau(H^{s-6})}.
\end{equation}
Using inequality \eqref{eq:pidi} in Lemma \ref{lem:prelim2} for $\beta=0$ and $\beta=1$, together with the local error estimate \eqref{eq:lest} established above, we finally have 
$$
\|e^\eps_{n+1}\|_{L^2_\tau(H^{s-6})} \leq \Big(1+C_1 \Delta t (1+C_2 \Delta t)\Big) \|e^\eps_{n}\|_{L^2_\tau(H^{s-6})} +C_3 \Delta t^3
$$
and owing to a discrete version of Gronwall lemma, we end up with a $L^2$ estimate for $e^\eps_{n}$
$$
\|e^\eps_{n}\|_{L^2_\tau(H^{s-6})} \leq  C \Delta t^2.
$$

\ms
\ni
{\bf Step 2: first order estimate of the global error in $L^\infty_\tau(H^{s-6})$.} In this second step, we prove:
\begin{equation}
\label{first}
\forall \, \Delta t < \Delta t_0, \qquad \sup_{\eps\in ]0,\eps_0]} \|e_n^\eps\|_{L^\infty_{\tau}(H^{s-6})} \leq C\Delta t. 
\end{equation}
To this aim, we estimate $\pa_\tau e_n^\eps$ in $L^2_\tau(H^{s-6})$. Using \eqref{eq:pidi} for $\beta=1$ (remarking that $Q_{2\mu}^{-1}$ and $\pa_\tau$ commute) and \eqref{eq:pidi2}, we deduce from \eqref{Qtransport} that
$$\|\pa_\tau \te_{n+1}^\eps\|_{L^2_\tau(H^{s-6})}\leq \|\pa_\tau e_n^\eps\|_{L^2_\tau(H^{s-6})}+4\mu \Delta t C_1(1+C_2\Delta t) \|e^\eps_{n}\|_{L^2_\tau(H^{s-6})}.$$
Therefore, since $\mu=\eps/\Delta t$, we deduce from \eqref{first0}, from the local truncation error \eqref{eq:LEst} (with $k=0$) and, again, from Lemma \ref{lem:prelim2} that
$$\|\pa_\tau e_{n+1}^\eps\|_{L^2_\tau(H^{s-6})}\leq \|\pa_\tau e_n^\eps\|_{L^2_\tau(H^{s-6})}+C\Delta t^2,$$
which gives directly \eqref{first} after summation and after using the Sobolev embedding of  $H^1_\tau(H^{s-6})$ into $L^\infty_\tau(H^{s-6})$.

%%%%%
At this stage, we can already ensure that \eqref{boots} is satisfied, provided that $\Delta t<\Delta t_0$ with $\Delta t_0:=R/C$, the constant $C$ being given in \eqref{first}.  

\ms
\ni
{\bf Step 3: second order estimate of the global error in $L^\infty_\tau(H^{s-8})$.} The last step of the proof consists in proving the same estimate as \eqref{first0} for $\pa_\tau e^\eps_{n}$, but in $L^2_\tau(H^{s-8})$. This proof is similar as Step 1, up to replacing the vector field $\cF(t,U)$ by $\cG(t,U,V)=\pa_\tau \cF(t,U) + \pa_u \cF(t,U)V$ and the local truncation error \eqref{eq:lest} by \eqref{eq:LEst} (with $k=1$). Therefore, we point out the following estimate on $\cG$. 
If $U$, $\tilde U$ belong to a bounded set of $L^\infty_\tau(H^{s-6})$ and if $V$, $\tilde V$ belong to a bounded set of $L^2_\tau(H^{s-8})$, then, for all $t$,
\begin{align*}
&\Big\|\cG(t,\tilde{U},\tilde{V})-\cG(t,U,V)\Big\|_{L^2_\tau(H^{s-8})}\\
&\qquad \leq \Big\|\pa_\tau \cF(t,\tilde{U})-\pa_\tau \cF(t,U)\Big\|_{L^2_\tau(H^{s-8})}+\Big\|\pa_u\cF(t,\tilde{U})\tilde{V}-\pa_u \cF(t,U)V\Big\|_{L^2_\tau(H^{s-8})}\\
&\qquad \leq C \|\tilde U-U\|_{L^2_\tau(H^{s-8})} +\Big\| (\pa_u\cF(t,\tilde{U})-\pa_u \cF(t,U))\tilde{V}\Big\|_{L^2_\tau(H^{s-8})} +\Big\| \pa_u\cF(t,\tilde{U})(\tilde{V}-V)\Big\|_{L^2_\tau(H^{s-8})}\\
%&\qquad = C \|\tilde U-U\|_{L^\infty_\tau(H^{s-8})} +\Big\| \int_0^1 (\pa^2_u\cF(t,\lambda \tilde{U} + (1-\lambda)U) (\tilde{U}-U, \tilde{V}) d\lambda \Big\|_{L^2_\tau(H^{s-8})}\\
%&+\Big\| \pa_u\cF(t,\tilde{U})(\tilde{V}-V)\Big\|_{L^2_\tau(H^{s-8})}\\
&\qquad \leq C \|\tilde U-U\|_{L^2_\tau(H^{s-8})} + C \| \tilde{U}-U  \|_{L^\infty_\tau(H^{s-8})} \| \tilde{V} \|_{L^2_\tau(H^{s-8})}+ C\| \tilde{V}-V\|_{L^2_\tau(H^{s-8})}\\
&\qquad \leq C \|\tilde U-U\|_{L^2_\tau(H^{s-8})} +C \| \partial_\tau \tilde{U}- \partial_\tau U  \|_{L^2_\tau(H^{s-8})}+ C\| \tilde{V}-V\|_{L^2_\tau(H^{s-8})}, \\
%&\qquad \leq C \|\tilde U-U\|_{L^\infty_\tau(H^{s-6})} +C \|\tilde V-V\|_{L^2_\tau(H^{s-8})}.
\end{align*} 
where the Sobolev embedding of $H^1_\tau(H^{s-8})$ into $L^\infty_\tau(H^{s-8})$ has been used for the last  
inequality. With $\tilde{U} = \tilde{U}^\eps_{n+1/2}, U = U^\eps_{n+1/2}, \tilde{V} = \partial_\tau \tilde{U}^\eps_{n+1/2}, V = \partial_\tau U^\eps_{n+1/2}$ and using \eqref{enp12} and \eqref{first0}, we have 
$$
\Big\|\cG(t,\tilde{U},\tilde{V})-\cG(t,U,V)\Big\|_{L^2_\tau(H^{s-8})} \leq C\Delta t^2 + C\|\partial_\tau \tilde{e}_{n+1/2}^\eps\|_{L^2_\tau(H^{s-8})}. 
$$
Arguing like in Step 1, we deduce that $\|\partial_\tau e^\eps_{n+1}\|_{L^2_\tau(H^{s-8})}\leq (1+C\Delta t)\|\partial_\tau e^\eps_{n}\|_{L^2_\tau(H^{s-8})} +C\Delta t^3$ and we finally conclude that 
$$
\sup_{\eps\in]0,\eps_0]}\|{e}^\eps_{n}\|_{H^1_\tau(H^{s-8})} \leq C\Delta t^2. 
$$
The Sobolev embedding of  $H^1_\tau(H^{s-8})$ into $L^\infty_\tau(H^{s-8})$ allows again to obtain the desired error estimate
$$
\sup_{\eps\in]0,\eps_0]}\|{e}^\eps_{n}\|_{L^\infty_\tau(H^{s-8})} \leq C\Delta t^2.
$$
\end{Proof}
\subsection{Choice of the initial data for the second order scheme}
\label{subsecond}
As in Subsection \ref{subfirst}, let us explain our strategy to construct the initial condition $U_0^\eps(\tau)$ for the augmented problem \eqref{eqF}, as soon as the initial condition $u_0$ for \eqref{eqftildegene} is known.

We set\footnote{For simplicity, we omit here the $t$ dependency in $h_1$ and $h_2$ which are all evaluated at $t=0$.}
\begin{align}
U_0^\eps(\tau):=&u_0+\eps h_1(\tau,u_0)-\eps h_1(0,u_0)+\eps^2 h_2(\tau,u_0)-\eps^2 h_2(0,u_0)\nonumber\\
&-\eps^2\pa_u h_1(\tau,u_0)h_1(0,u_0)+\eps^2\pa_u h_1(0,u_0)h_1(0,u_0),\label{secondID}
\end{align}
where $h_1$ and $h_2$ are defined by \eqref{h1} and \eqref{h2}, and 
$$\ub U_0^\eps:=\Pi U_0^\eps=u_0-\eps h_1(0,u_0)-\eps^2 h_2(0,u_0)+\eps^2\pa_u h_1(0,u_0)h_1(0,u_0).$$
In the numerical  Section \ref{sect:num}, this choice of initial data will be referred to as {\em second order initial data}.
Assuming that $u_0\in H^{s+4}$, one deduces from Assumption \ref{mainassump} that $\ub U_0^\eps$ is uniformly bounded in $H^{s+2}$. Next, consider the remainder term defined according to \eqref{condinit0} by
$$
r^\eps(\tau):=\frac{1}{\eps^3}\left(U_0^\eps(\tau)-\ub U^\eps_0 - \eps h_1(\tau, \ub U^\eps_0)- \eps^2 h_2(\tau, \ub U^\eps_0)\right).
$$
By Taylor expansions, one gets
\begin{align*}
r^\eps(\tau)&=\frac{1}{\eps^2}\big(h_1(\tau,u_0)-\eps\pa_u h_1(\tau,u_0)h_1(0,u_0)-h_1(\tau,\ub U_0^\eps)+\eps h_2(\tau,u_0)-\eps h_2(\tau,\ub U_0^\eps)\big)\\
&=\frac{1}{\eps^2}\big(h_1(\tau,u_0)-\eps\pa_u h_1(\tau,u_0)h_1(0,u_0)-h_1(\tau,u_0-\eps h_1(0,u_0))+{\mathcal O}_{X^s}(\eps^2)\big)\\
&=\mathcal O_{X^s}(1).
\end{align*}
The above choice of initial data ensures the validity of the assumptions of Theorem \ref{theo-o2}.

\section{Applications and numerical results}\label{sect:num}

In this section, we apply our two-scale technique in two situations. We first present numerical experiments for the nonlinear Klein-Gordon equation in the nonrelativistic limit regime, then we consider a stiffer problem, the nonlinear Schr\"odinger equation in a highly oscillatory regime.

A special care will be given to the choice of the initial condition $U^\eps_0(\tau)$ for the augmented problem, that will be constructed from the initial condition $u_0$. The possible choices will be:
\begin{itemize}
\item[--] the {\em uncorrected initial data} $U^\eps_0=u_0$,
\item[--] the {\em first order initial data} $U^\eps_0$ given by \eqref{firstID},
\item[--] the {\em second order initial data} $U^\eps_0$ given by \eqref{secondID},
\item[--] a {\em third order initial data} $U^\eps_0$ obtained by pushing the Chapman-Enskog expansion of Subsection \ref{chapman} to the order three.
\end{itemize}

\subsection{The nonlinear Klein-Gordon equation in the nonrelativistic limit regime}

In this subsection, we consider the following nonlinear Klein-Gordon (NKG) equation:
\be
\label{kg}
\eps \pa_{tt}u-\Delta u+\frac{1}{\eps}u+f(u)=0,\qquad x\in \R^d,\quad t>0,
\ee
with initial conditions given as
\be
\label{kginit}
u(0,x)=\phi(x),\qquad \pa_tu(0,x)=\frac{1}{\eps}\gamma(x),\qquad x\in \R^d.
\ee
The unknown is the function $u(t,x):\,\R^{1+d}\to \C$ and the parameter $\eps>0$ is inversely proportional to the square of the speed of light. The nonlinearity is assumed to be a smooth function $f$ satisfying $f(0)=0$, $f(\R)\subset \R$  and  the gauge invariance $f(e^{is}u)=e^{is}f(u)$ for all $u\in \C$, $s\in \R$.  This implies in particular that $f$ satisfies $f(\bar z)= \overline {f(z)}$ for all $z\in \C$. The limit $\eps\to 0$ in equation \eqref{kg}, \eqref{kginit}, referred to as the nonrelativistic limit, has been studied in \cite{kg1,kg2,kg3}. The regime of small $\eps$ (but not zero) is highly oscillatory and has been recently explored numerically in \cite{bao-dong} by Gauschi-type exponential methods, allowing for time steps of order $\mathcal O(\eps)$. In \cite{faou-schratz}, a different approach is proposed, based on asymptotic expansions with respect to $\eps$. None of these methods is uniformly accurate in $\eps\in (0,1]$. We apply here our two-scale reformulation technique which naturally leads to uniformly accurate numerical schemes.

It is convenient to rewrite \eqref{kg} under the equivalent form of a first order system in time (see e.g. \cite{kg3,faou-schratz}). Setting
\be
\label{change}
v_+=u-i\eps(1-\eps\Delta)^{-1/2}\pa_t u,\qquad v_-=\overline{u}-i\eps(1-\eps\Delta)^{-1/2}\pa_t \overline{u},
\ee
and denoting
$$\widetilde f(v_+,v_-)=\left(f(\frac{1}{2}(v_++\overline{v_-}))\,,\,f(\frac{1}{2}(\overline{v_+}+v_-))\right),$$
we obtain that \eqref{kg}, \eqref{kginit} is equivalent to the following system on the unknown $v=(v_+,v_-)$:
\be
i\pa_t v=-\frac{1}{\eps}(1-\eps\Delta)^{1/2}v-(1-\eps\Delta)^{-1/2}\widetilde f(v),
\label{kg2}
\ee
\be
\label{kg2init}
v(0,\cdot)=(v_+(0,\cdot),v_-(0,\cdot))=\left(\phi-i(1-\eps\Delta)^{-1/2}\pa_t \gamma\,,\,\overline{\phi}-i(1-\eps\Delta)^{-1/2}\pa_t \overline{\gamma}\right).
\ee
Finally, introducing the filtered unknown
$$\widetilde u=e^{-i\frac{t}{\eps}\sqrt{1-\eps\Delta}}v,$$
we obtain the following equation:
\be
\label{kg3}
\pa_t\widetilde u=i(1-\eps\Delta)^{-1/2}e^{-i\frac{t}{\eps}\sqrt{1-\eps\Delta}}\widetilde f\left(e^{i\frac{t}{\eps}\sqrt{1-\eps\Delta}}\widetilde u\right).
\ee
Note that \eqref{kg3} is under the form \eqref{eqftildegene} if we set
\be
\label{cF1}
\cF(t,\tau,u,\eps)=i(1-\eps\Delta)^{-1/2}e^{-i\tau}e^{-itA_\eps}\widetilde f\left(e^{i\tau}e^{itA_\eps}u\right)
\ee
with
$$A_\eps=\frac{1}{\eps}\left(\sqrt{1-\eps\Delta}-1\right).$$
This self-adjoint operator is not singular as $\eps\to 0$: for all $\eps>0$, we have
$$0\leq A_\eps\leq -\frac{1}{2}\Delta$$
in the sense of operators. It can be checked that this vector field $\cF$ satisfies Assumption \ref{mainassump}.

\bs
For our numerical tests, we borrow an example in dimension $d=1$ from \cite{bao-dong}. We choose
$$f(u)=4|u|^2u,\qquad \phi(x)=\frac{2}{e^{x^2}+e^{-x^2}},\qquad \gamma(x)=0.$$
The final time of the simulations is $T_f=0.4$. The computational domain in $x$ is $[-8,8]$ and is large enough for  periodic boundary conditions to be taken, with no significative error compared to the solution in the whole space (this feature is a posteriori checked). For the numerical evaluations of the function $\cF$ in our scheme, a spectral method is used in the $x$ variable and the fast Fourier transform is used in the practical implementation. For the series of tests, the reference solution is computed as follows. For $\eps\geq 10^{-2}$, we use the Yoshida fourth order splitting method \cite{yoshida} with $\Delta x=16/256=0.0625$, $\Delta t=\eps\,T_f/2000$. For smaller values of $\eps$, we rather use our second order uniformly accurate scheme, with small grid steps: $\Delta x=16/256=0.0625$, $\Delta t=2\pi/512000\approx 1.2\times 10^{-5}$, $\Delta \tau= 2\pi/128\approx 0.05$.
We shall use the $H^s$ relative error of a given numerical scheme which we define as
\be
{\cal E}_{s}= \frac{\|u^{ref}(t_{final},\cdot)-u^{num}(t_{final},\cdot)\|_{H^s}}{\|u^{ref}(t_{final},\cdot)\|_{H^s}},
\ee
where $u^{num}(t_{final},\cdot)$ is the approximated solution  obtained by the considered numerical scheme, at the final time $t_{final}$ of the simulation.
In order to validate the reference solution and show the behavior of a non uniformly accurate scheme, we first compare the reference solution $u^{ref}(t_{final},x)$ to the numerical solution $u^{Strang}(t_{final},x)$ computed with the following Strang splitting algorithm for \eqref{kg2}:

\ss
-- Step 1 for $t \in [t_n, t_n+\frac{\Delta t}{2}]$: we solve 
$$
i\pa_t v_1=-\frac{1}{\eps}(1-\eps\Delta)^{1/2}v_1,\quad v_{1\mid t=t_n}=v^n
$$
which has an explicit solution in the Fourier space.

\ss
-- Step 2 for $t \in [t_n, t_n+\Delta t]$: we solve 
$$
i\pa_t v_2=-(1-\eps\Delta)^{-1/2}\widetilde f(v_2),\quad v_{2\mid t=t_n}=v_{1\mid t= t_n+\frac{\Delta t}{2}}
$$
which has also an explicit solution (remark indeed that the solution $v_2=(v_{2+},v_{2-})$ of this equation satisfies $v_{2+}+\overline {v_{2-}}=${\em constant}).

\ss
-- Step 3 for $t \in [t_n+\frac{\Delta t}{2},t_n+\Delta t]$: we solve 
$$
i\pa_t v_3=-\frac{1}{\eps}(1-\eps\Delta)^{1/2}v_3,,\quad v_{3\mid t=t_n+\frac{\Delta t}{2}}=v_{2\mid t= t_n+\Delta t}
$$
We set finally $v^{n+1}=v_{3\mid t=t_n+\Delta t}$.

\begin{figure}[!htbp]
  \centerline{
  \subfigure[Error with respect to $\Delta t$]{\includegraphics[width=.55\textwidth]{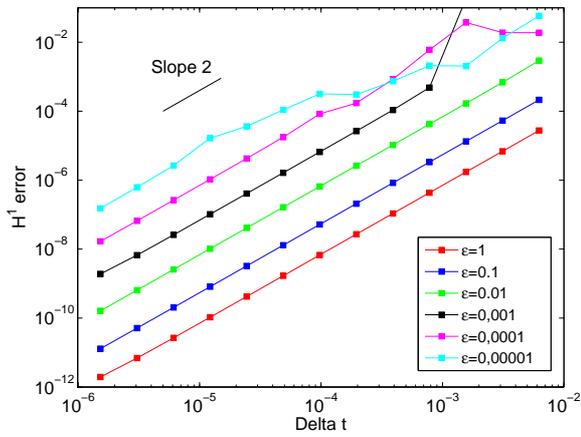}}\hspace*{-7mm}
  \subfigure[Error with respect to $\eps$]{\includegraphics[width=.717\textwidth]{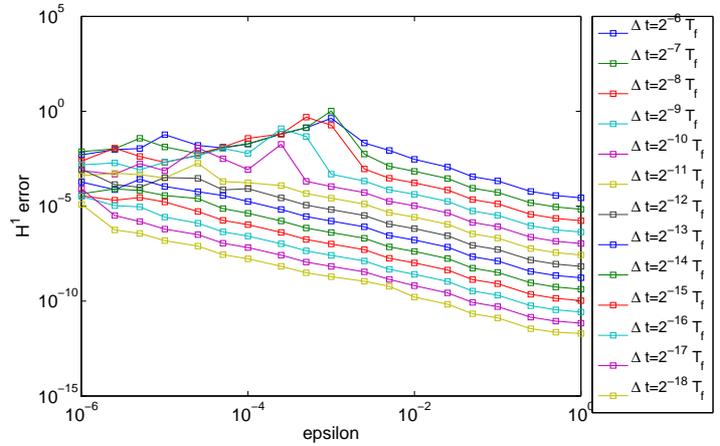}}}
    \caption{(NKG) $H^1$ relative error (in log-log scale)for the Strang splitting scheme.}
  \label{fig0}
\end{figure}

\begin{figure}[!htbp]
  \centerline{\includegraphics[width=.60\textwidth]{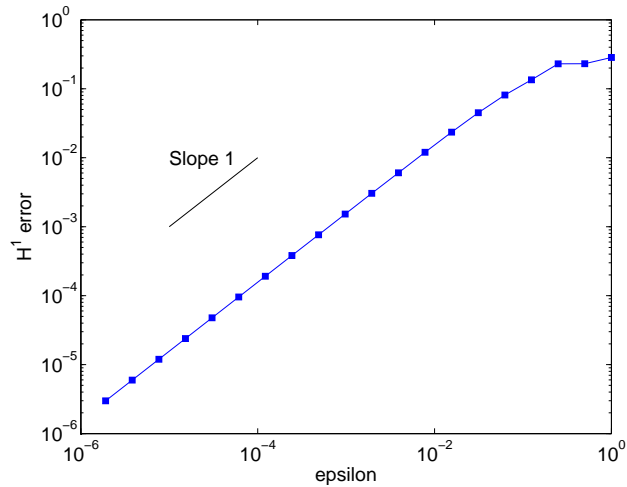}}
    \caption{(NKG) $H^1$ relative error (in log-log scale) between the reference solution and the limiting averaged model.}
  \label{fig9}
\end{figure}

On Figure \ref{fig0}, we represent the $H^1$ error between the reference solution and the numerical solution computed with this Strang splitting scheme, with a fixed number of grid points in $x$, $N_x=200$, for various values $\Delta t=2^{-K}\,T_f$ with $K\in\{6,\ldots,18\}$ and for various values of $\eps$ from $\eps=1$ to $\eps=10^{-6}$. It appears numerically that the $H^1$ relative error behaves asymptotically like $C\frac{\Delta t^2}{\eps}$, 
where $C$ is constant which does not depend on $\Delta t$ and $\eps$.
%$$\|u^{ref}(t_{final},\cdot)-u^{Strang}(t_{final},\cdot)\|_{L^2}\sim C\frac{\Delta t^2}{\eps}.$$
The Strang splitting scheme becomes inefficient for small values of $\eps$. A natural idea is to use instead an asymptotic model as $\eps\to 0$, which is not stiff with respect to the parameter $\eps$.  Let us illustrate the limitation of the use of the limit averaged model. As $\eps\to 0$, the solution $v$ of the nonlinear Klein-Gordon equation \eqref{kg2} behaves asymptotically as follows:
\be
\label{estimave}
\|v(t,x)-e^{it/\eps}w(t,x)\|\leq C\eps,
\ee
where $w$ solves the averaged equation
\be
i\pa_tw=-\frac{1}{2}\Delta w+\frac{1}{2\pi}\int_0^{2\pi}e^{-i\tau}\widetilde f\left(e^{i\tau}w\right)d\tau.
\label{kglimit}
\ee
On Figure \ref{fig9}, we check numerically the error estimate \eqref{estimave}: we plot, with respect to $\eps$, the $H^1$ error between the reference solution and the numerical solution of \eqref{kglimit} (where the integral is discretized with the rectangle quadrature method), computed with small time and space grid steps. Clearly, the averaged model can only be used as an approximation of the original problem for very small values of $\eps$.

\bs
Instead, our two-scale method naturally leads to uniformly accurate numerical schemes. Let us now illustrate this property by studying the behavior of our first and second order schemes with respect to the various numerical parameters. 

\begin{figure}[!htbp]
  \centerline{
  \subfigure[Error with respect to $N_x$]{\includegraphics[width=.55\textwidth]{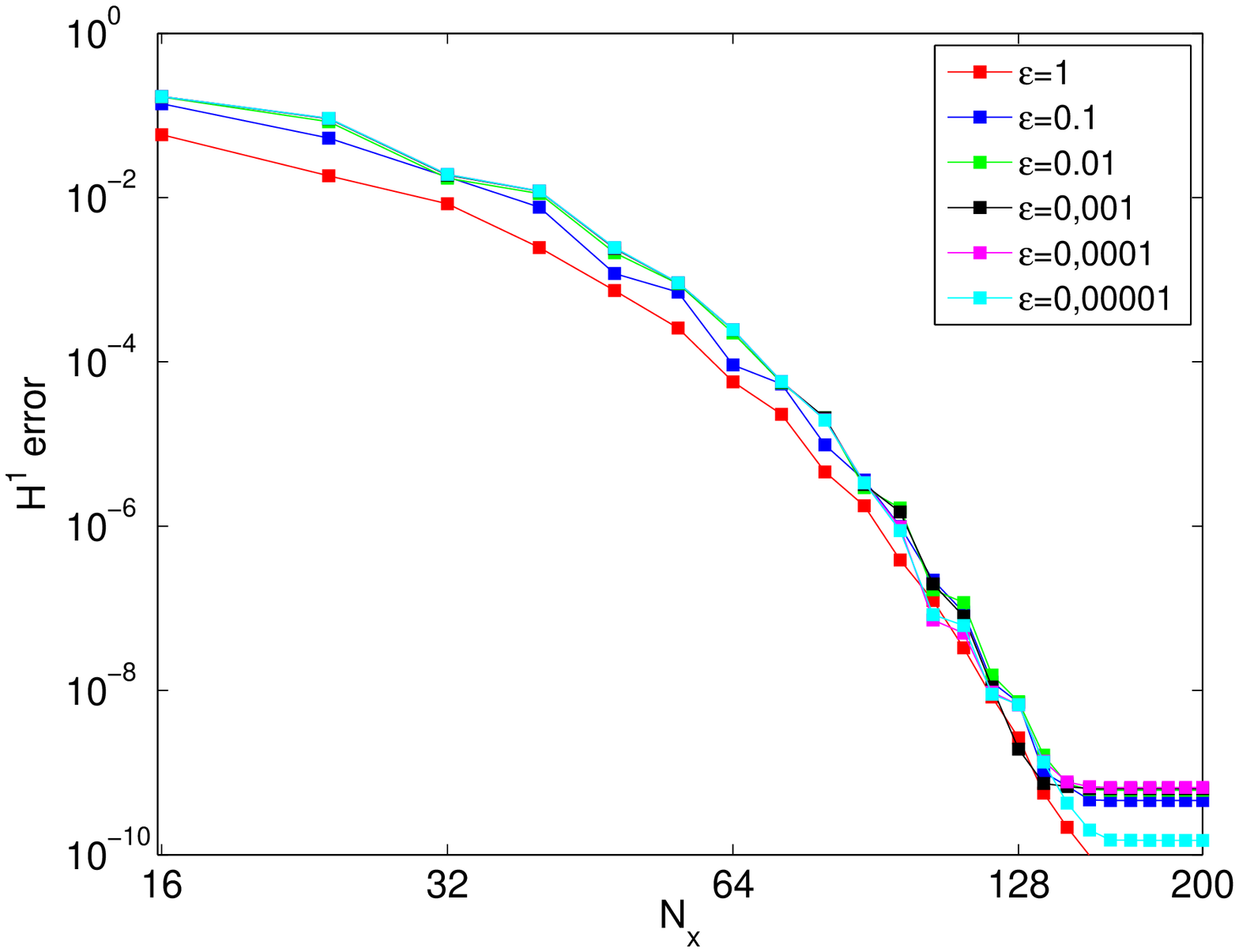}}\hspace*{-7mm}
  \subfigure[Error with respect to $N_\tau$]{\includegraphics[width=.55\textwidth]{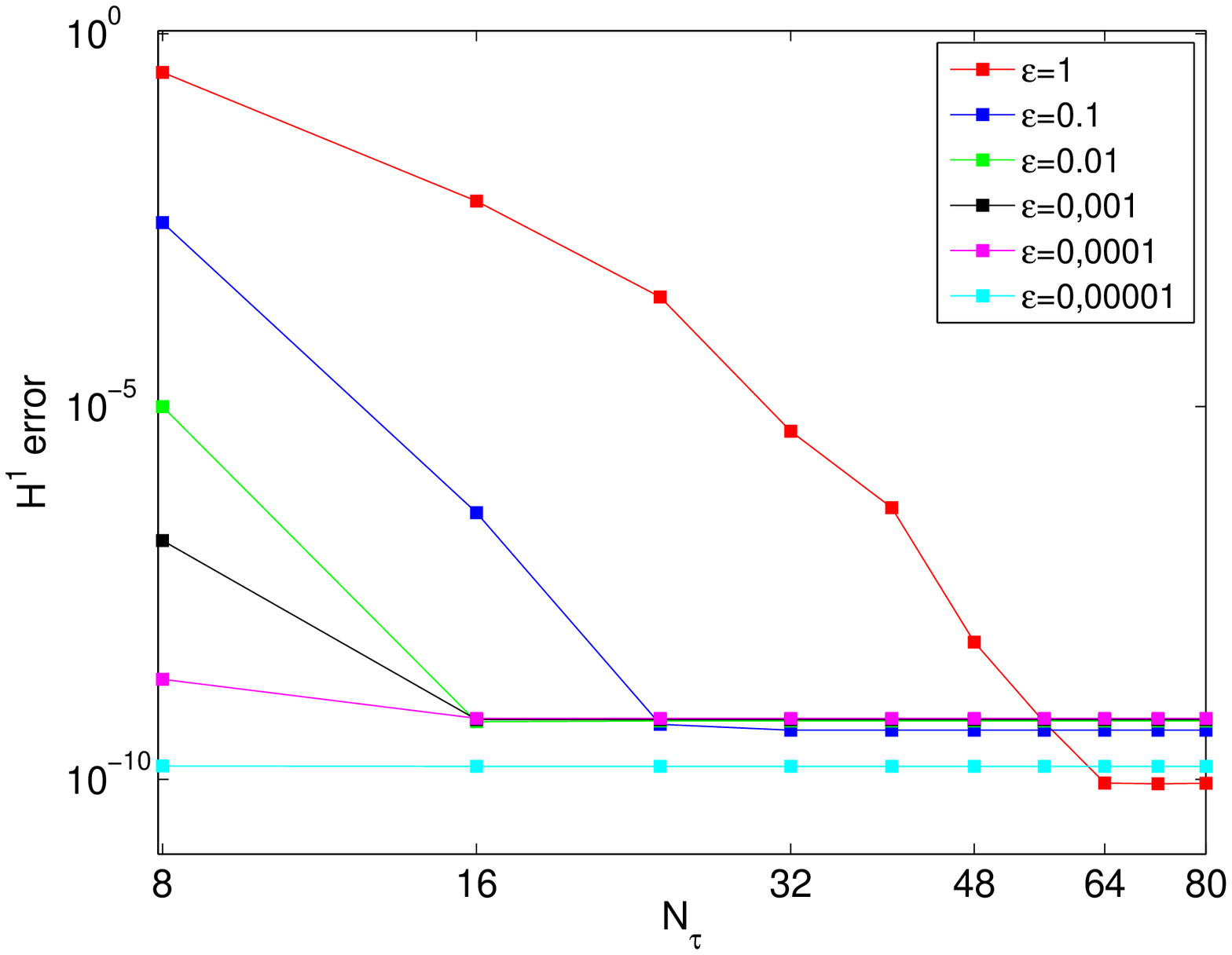}}}
    \caption{(NKG) $H^1$ relative error (in log-log scale) for the second order UA scheme in $\Delta t$ with the third order initial data.}
  \label{fig8}
\end{figure}

On Figure \ref{fig8}, we show that our scheme has a spectral accuracy with respect to the variables $x$ and $\tau$ (here, the time step is fixed $\Delta t=2\times 10^{-5}$). On the left part, we plot the $H^1$ error for our second order (in time) scheme with respect to the number $N_x$ of grid points in $x$. This error appears to be independent of $\eps$ and has a spectral behavior. On the right part of Figure \ref{fig8}, we plot the $H^1$ error for our scheme with respect to the number $N_\tau$ of gridpoints in the $\tau$ variable, illustrating also the spectral accuracy in this variable. Note that this error decreases rapidly when $\eps$ becomes small: for instance, for $\eps\leq 0.01$, $N_\tau=16$ would be sufficient. 

\begin{figure}[!htbp]
  \centerline{
  \subfigure[With the third order initial data]{\includegraphics[width=.55\textwidth]{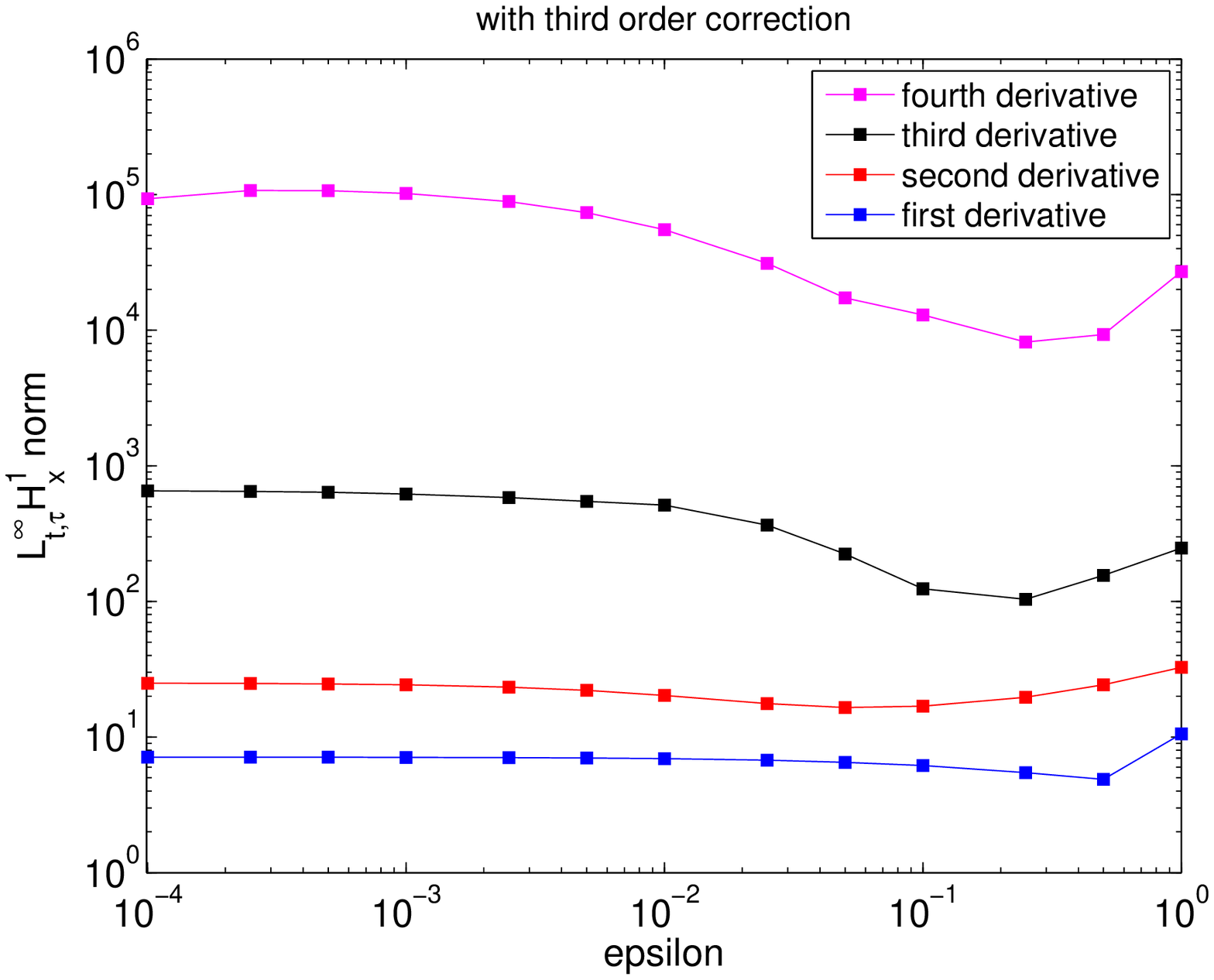}}\hspace*{-7mm}
  \subfigure[With the second order initial data]{\includegraphics[width=.55\textwidth]{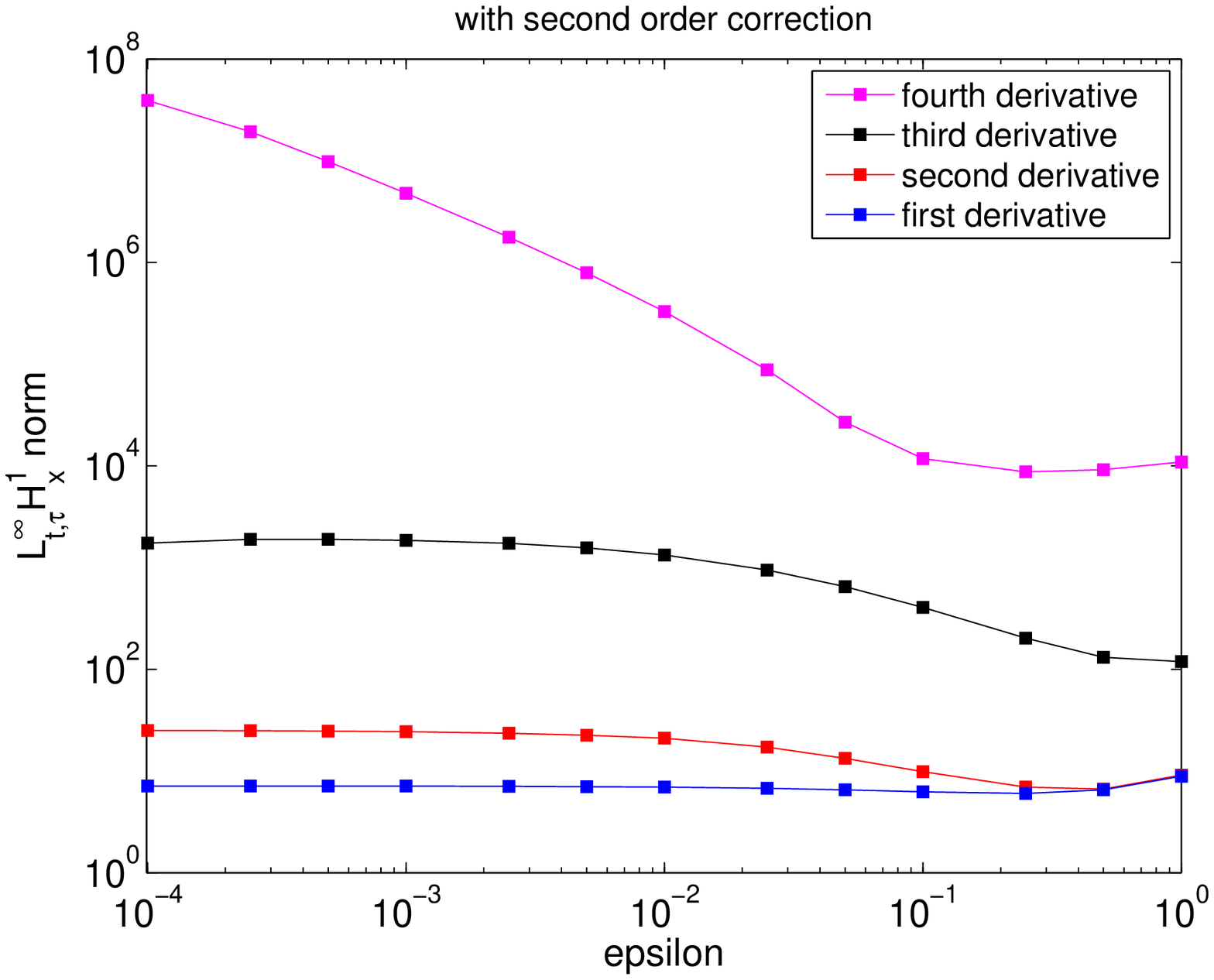}}}
    \caption{(NKG) $L^\infty_tL^\infty_{\tau}H^1_x$ norm  (in  log-log scale) of $\pa^k_tU^\eps(t,\tau,x)$ with respect to $\eps$, for $k\in\{1,2,3,4\}$.}
  \label{fig10}
\end{figure}

\begin{figure}[!htbp]
  \centerline{
  \subfigure[With the first order initial data]{\includegraphics[width=.55\textwidth]{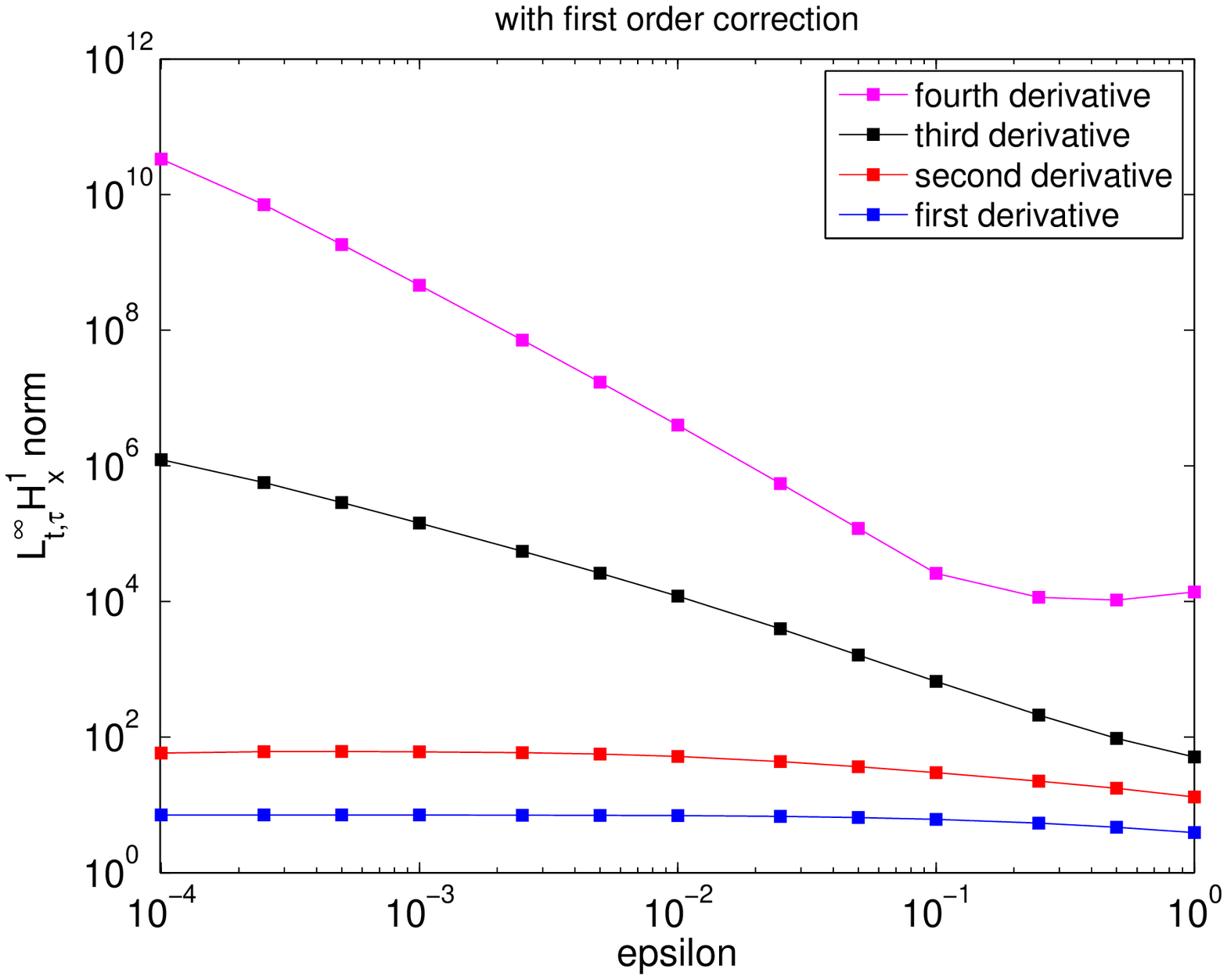}}\hspace*{-7mm}
  \subfigure[With the uncorrected initial data]{\includegraphics[width=.55\textwidth]{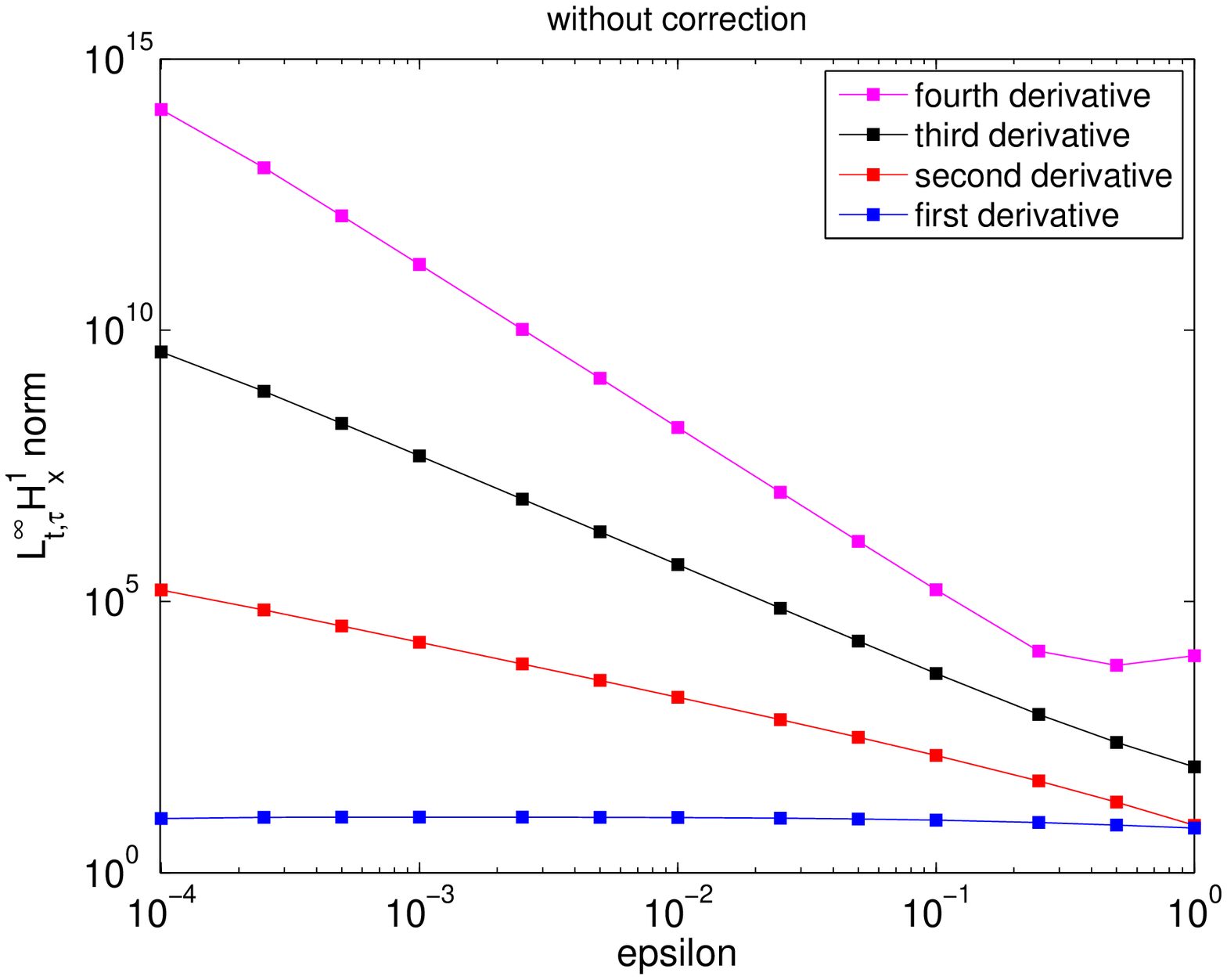}}}
    \caption{(NKG) $L^\infty_tL^\infty_\tau H^1_x$ norm of $\pa^k_tU^\eps(t,\tau,x)$ with respect to $\eps$, for $k\in\{1,2,3,4\}$  (in log-log scale).}
  \label{fig11}
\end{figure}

In the sequel, the space and $\tau$ grid steps are fixed: $N_x=200$ and $N_\tau=64$ are chosen. We now concentrate on the behavior with respect to the time step $\Delta t$. The above numerical analysis of our schemes shows that the optimal accuracy in $\Delta t$ can only be obtained if the initial data $U^\eps_0(\tau,x)$ for the augmented problem is chosen with enough correction terms in the asymptotic formula obtained by Chapman-Enskog expansion. On Figures \ref{fig10} and \ref{fig11}, we illustrate the importance of this choice by plotting the $L^\infty_tL^\infty_{\tau}H^1_x$ norms of the derivatives $\pa^k_tU^\eps(t,\tau=0,x)$, for $k\in\{1,2,3,4\}$, with respect to $\eps$ and with different choices of initial data. These curves indicate that, if the initial data is taken with $n$ correction terms, then we have the following behavior as $\eps\to 0$:
$$\pa^k_tU^\eps=\mathcal O(\eps^{n+1-k}).$$

\begin{figure}[!htbp]
  \centerline{
  \subfigure[Error with respect to $\Delta t$]{\includegraphics[width=.55\textwidth]{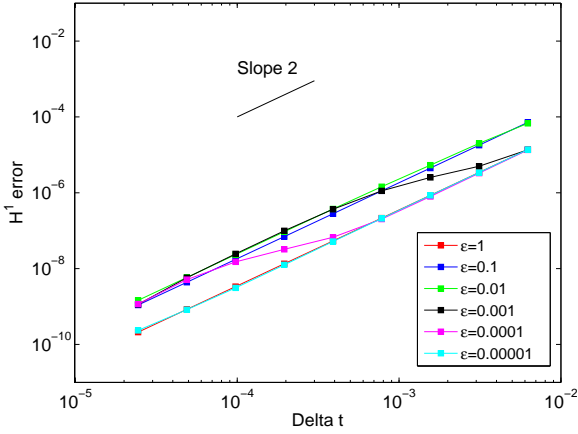}}\hspace*{-7mm}
  \subfigure[Error with respect to $\eps$]{\includegraphics[width=.717\textwidth]{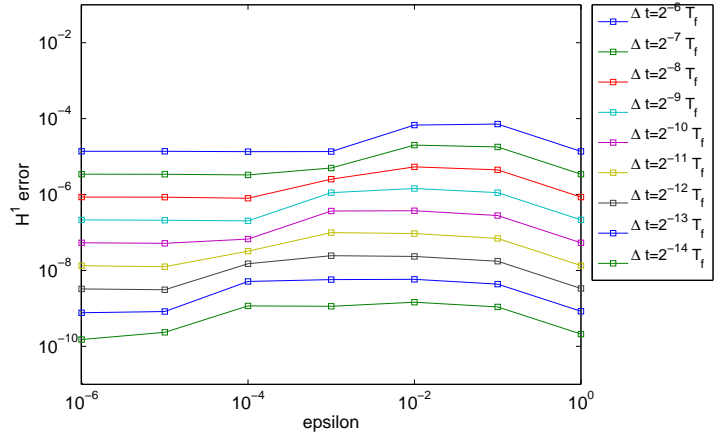}}}
    \caption{(NKG) $H^1$ relative error (in log-log scale) for the second order UA scheme with the third order initial data.}
  \label{fig1}
\end{figure}

\begin{figure}[!htbp]
  \centerline{
  \subfigure[Error with respect to $\Delta t$]{\includegraphics[width=.55\textwidth]{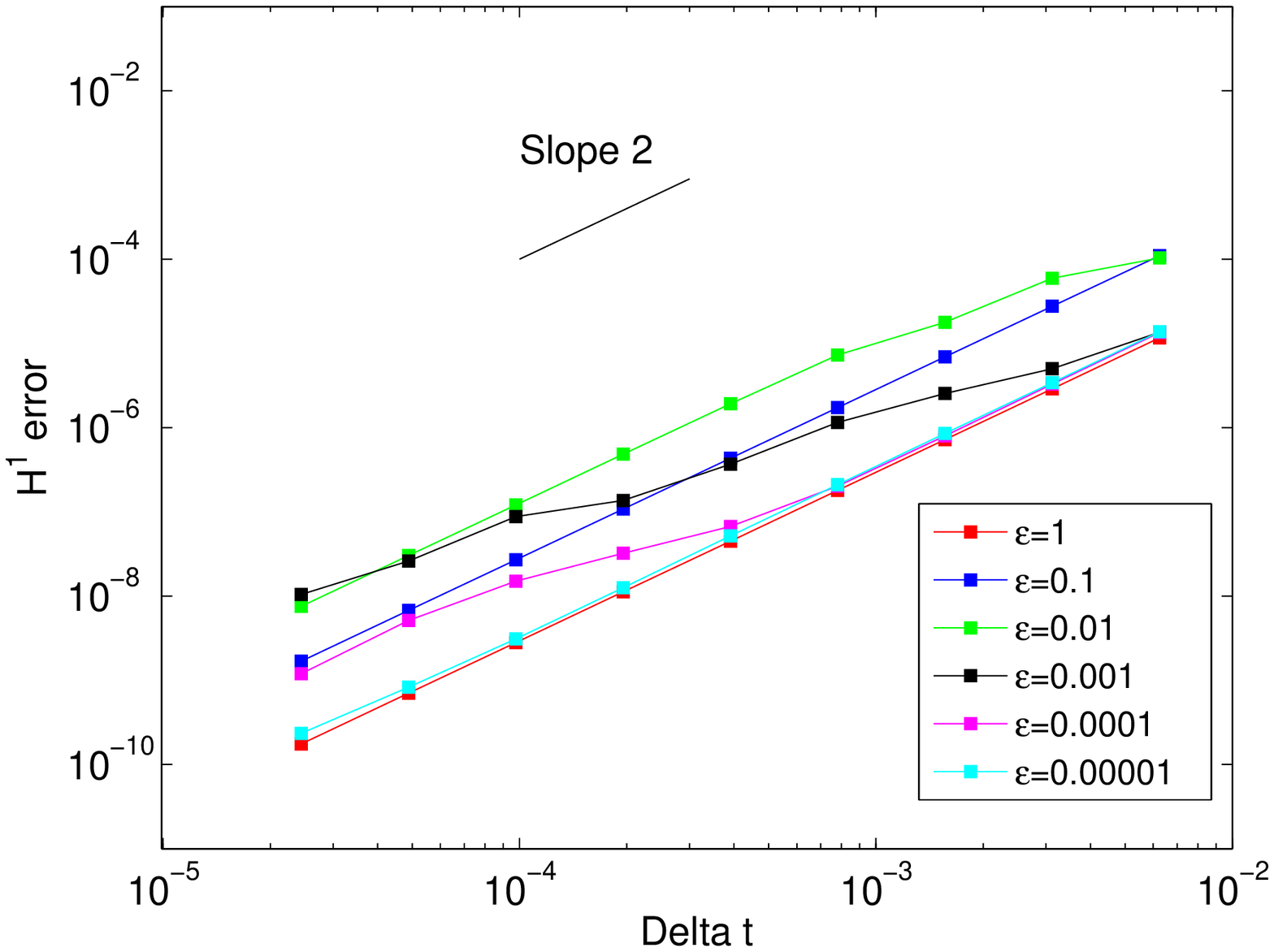}}\hspace*{-7mm}
  \subfigure[Error with respect to $\eps$]{\includegraphics[width=.717\textwidth]{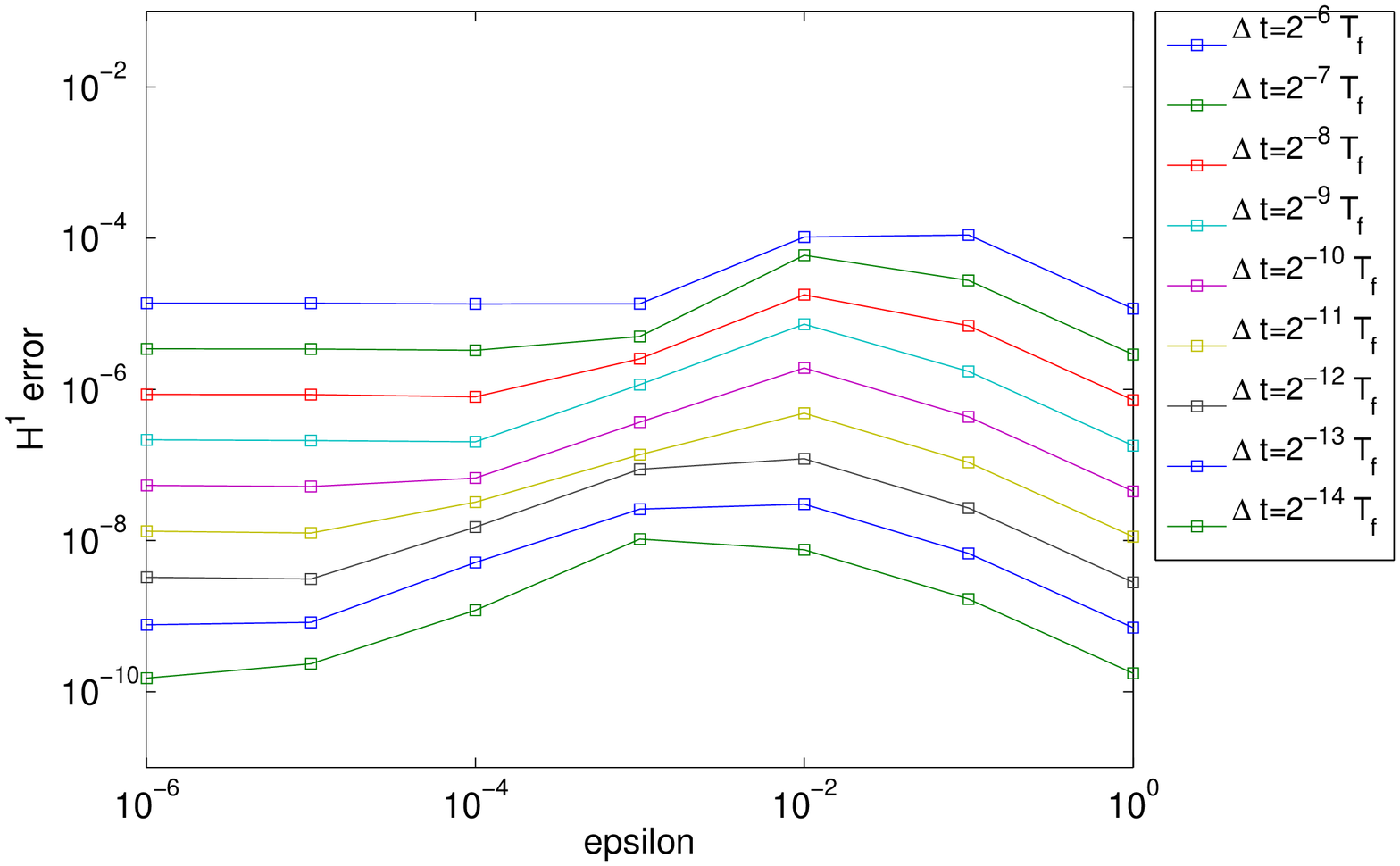}}}
    \caption{(NKG) $H^1$ relative error (in the log-log scale) for the second order UA scheme with the second order initial data.}
  \label{fig2}
\end{figure}

\begin{figure}[!htbp]
  \centerline{
  \subfigure[Error with respect to $\Delta t$]{\includegraphics[width=.55\textwidth]{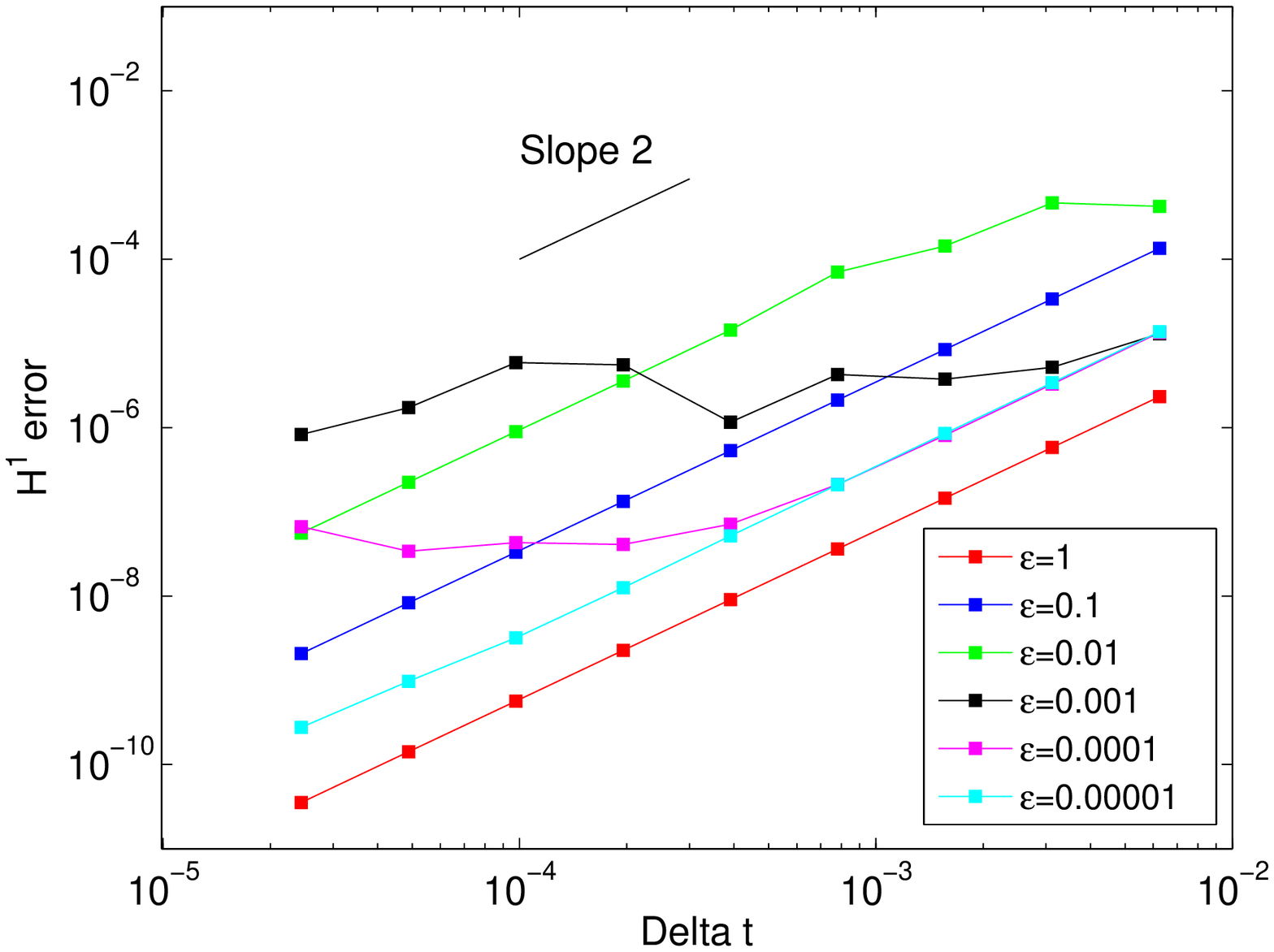}}\hspace*{-7mm}
  \subfigure[Error with respect to $\eps$]{\includegraphics[width=.717\textwidth]{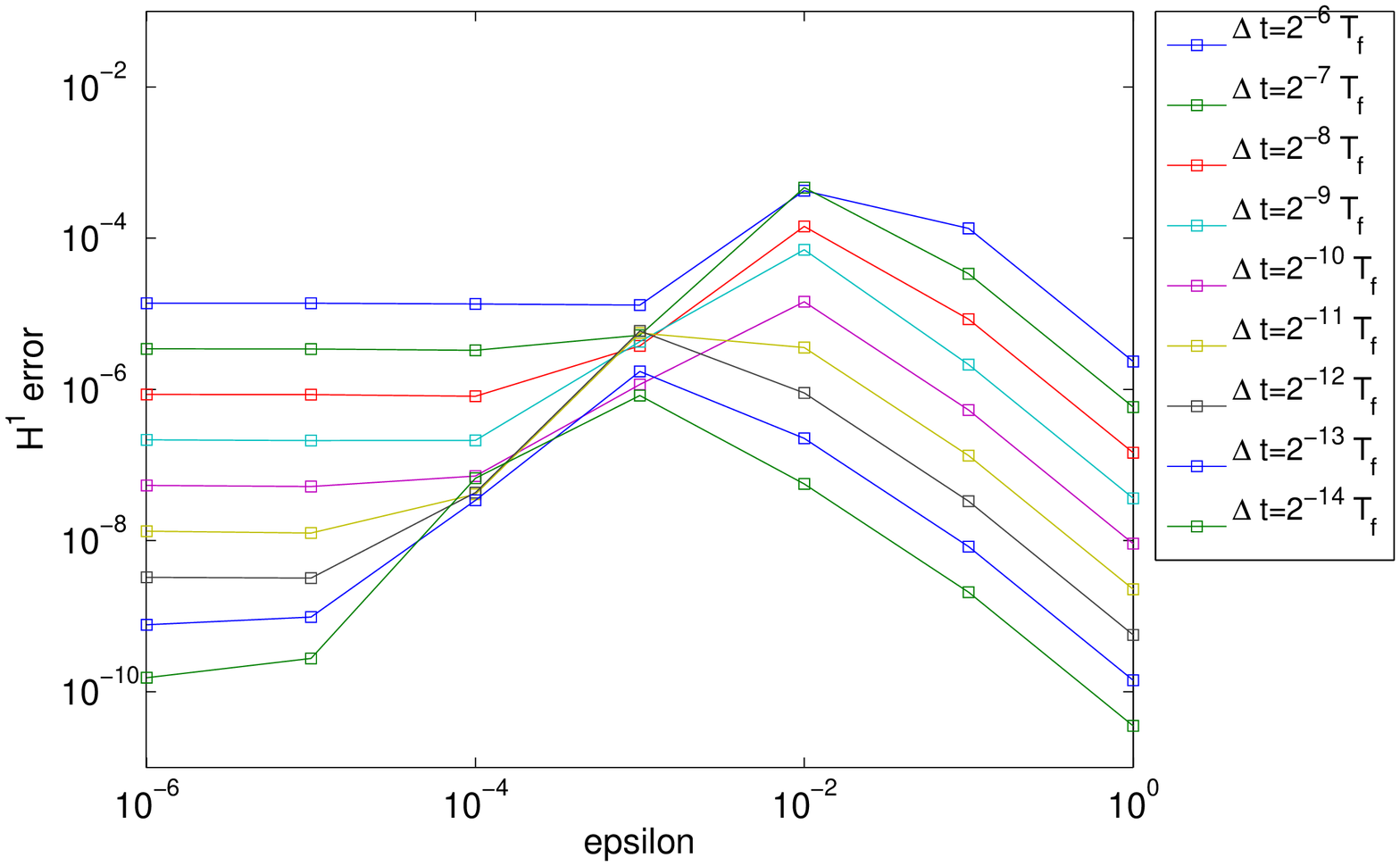}}}
    \caption{(NKG) $H^1$ relative error  (in log-log scale) for the second order UA scheme with the first order initial data.}
  \label{fig3}
\end{figure}

\begin{figure}[!htbp]
  \centerline{
  \subfigure[Error with respect to $\Delta t$]{\includegraphics[width=.55\textwidth]{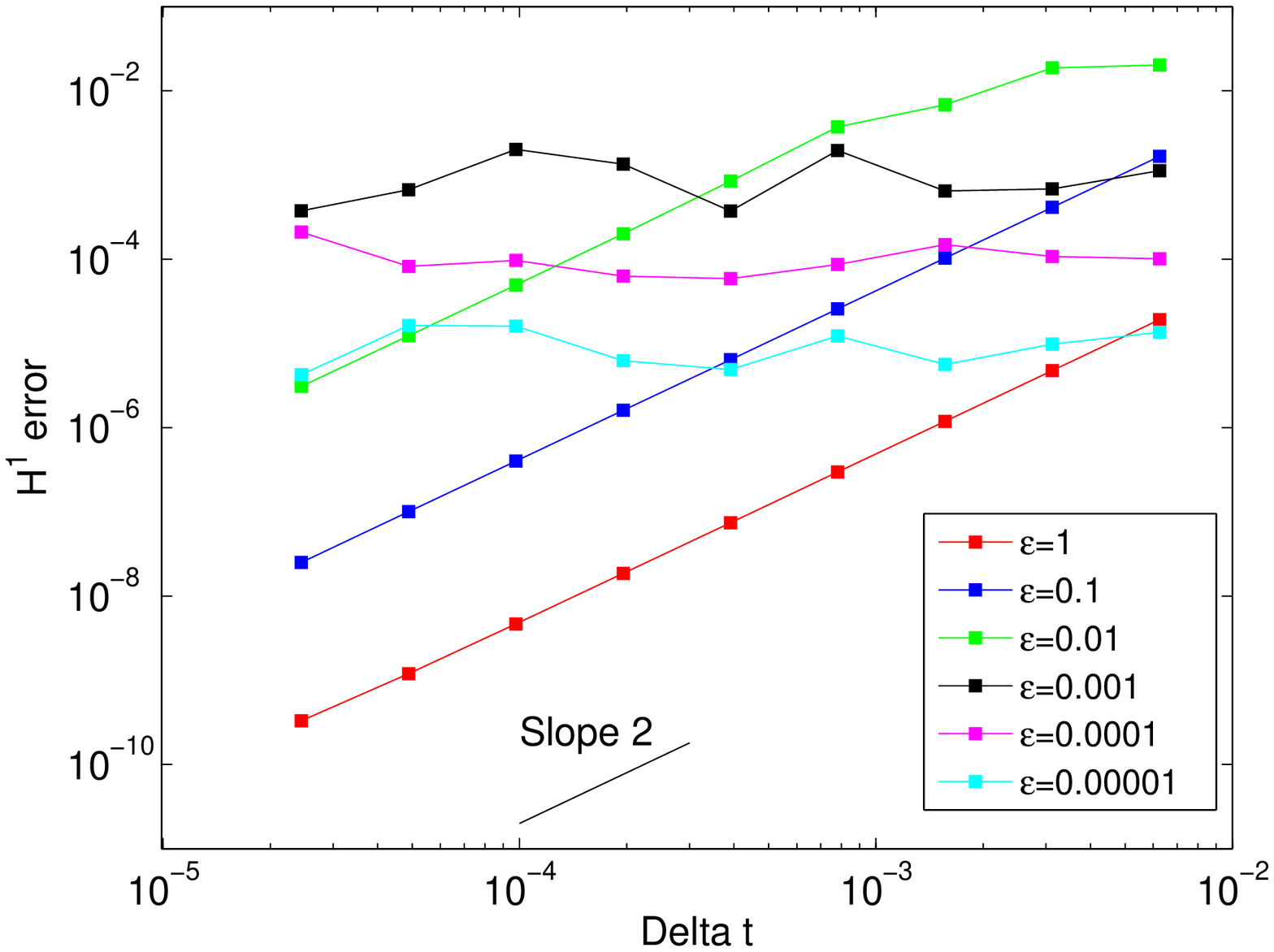}}\hspace*{-7mm}
  \subfigure[Error with respect to $\eps$]{\includegraphics[width=.717\textwidth]{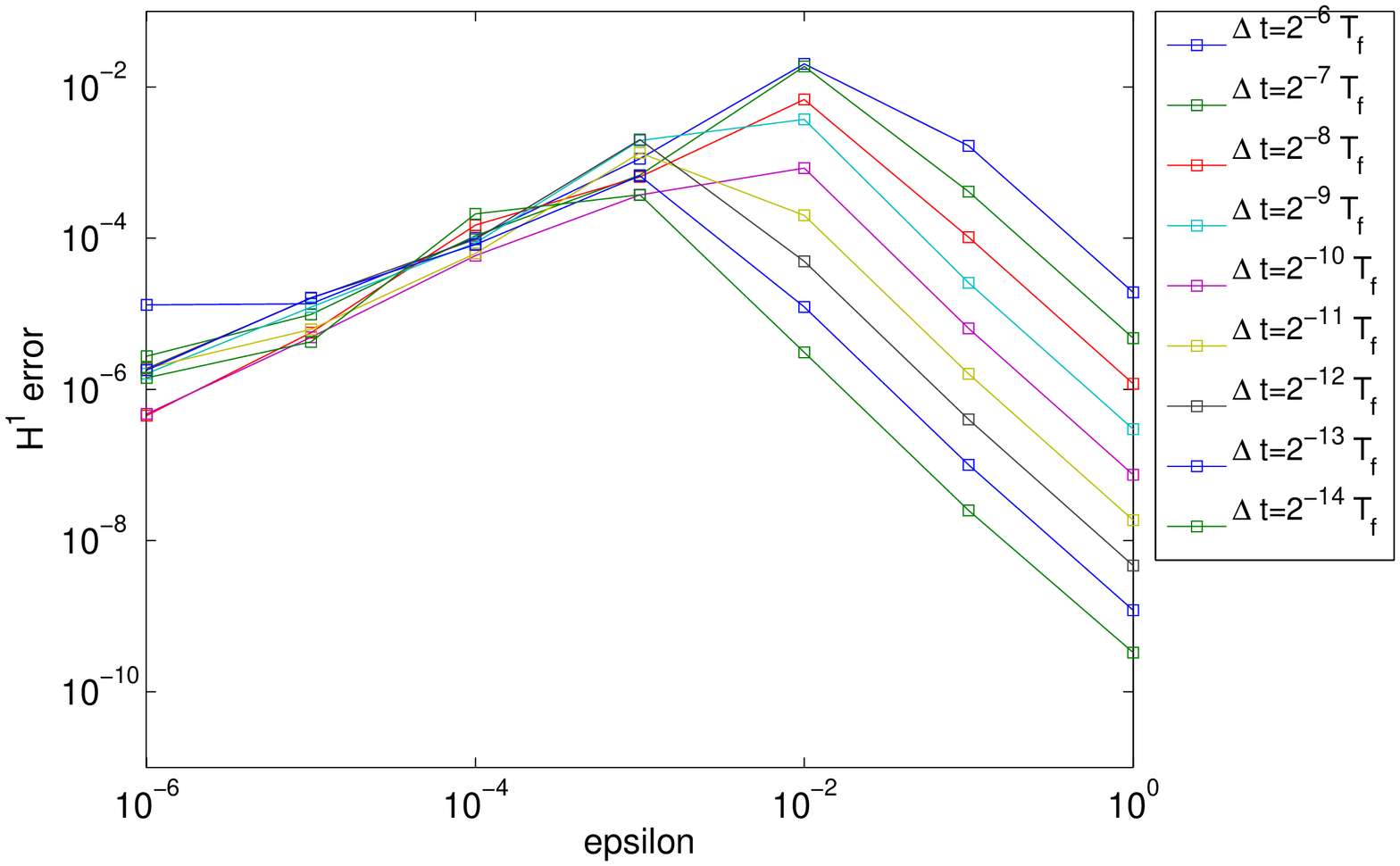}}}
    \caption{(NKG) $H^1$ relative error (in  log-log scale) for the second order UA scheme with the uncorrected intial data.}
  \label{fig4}
\end{figure}

On Figures \ref{fig1}, \ref{fig2}, \ref{fig3}, \ref{fig4}, we plot the $H^1$ error for our second order scheme with respect to $\Delta t$ and $\eps$, for four choices of initial data $U_0$. It appears clearly that, as expected, the uniform second order accuracy is obtained for the second or third order corrected initial data, with better results in the case  of the third order initial data (that we explain by the fact that the fourth derivative in time has an influence on the constants in the error estimate). If the initial data $U_0$ is not taken with enough correction terms, the second order accuracy is lost for intermediate regimes of $\eps$ (see Figures \ref{fig3} and \ref{fig4}).

\begin{figure}[!htbp]
  \centerline{
  \subfigure[Error with respect to $\Delta t$]{\includegraphics[width=.55\textwidth]{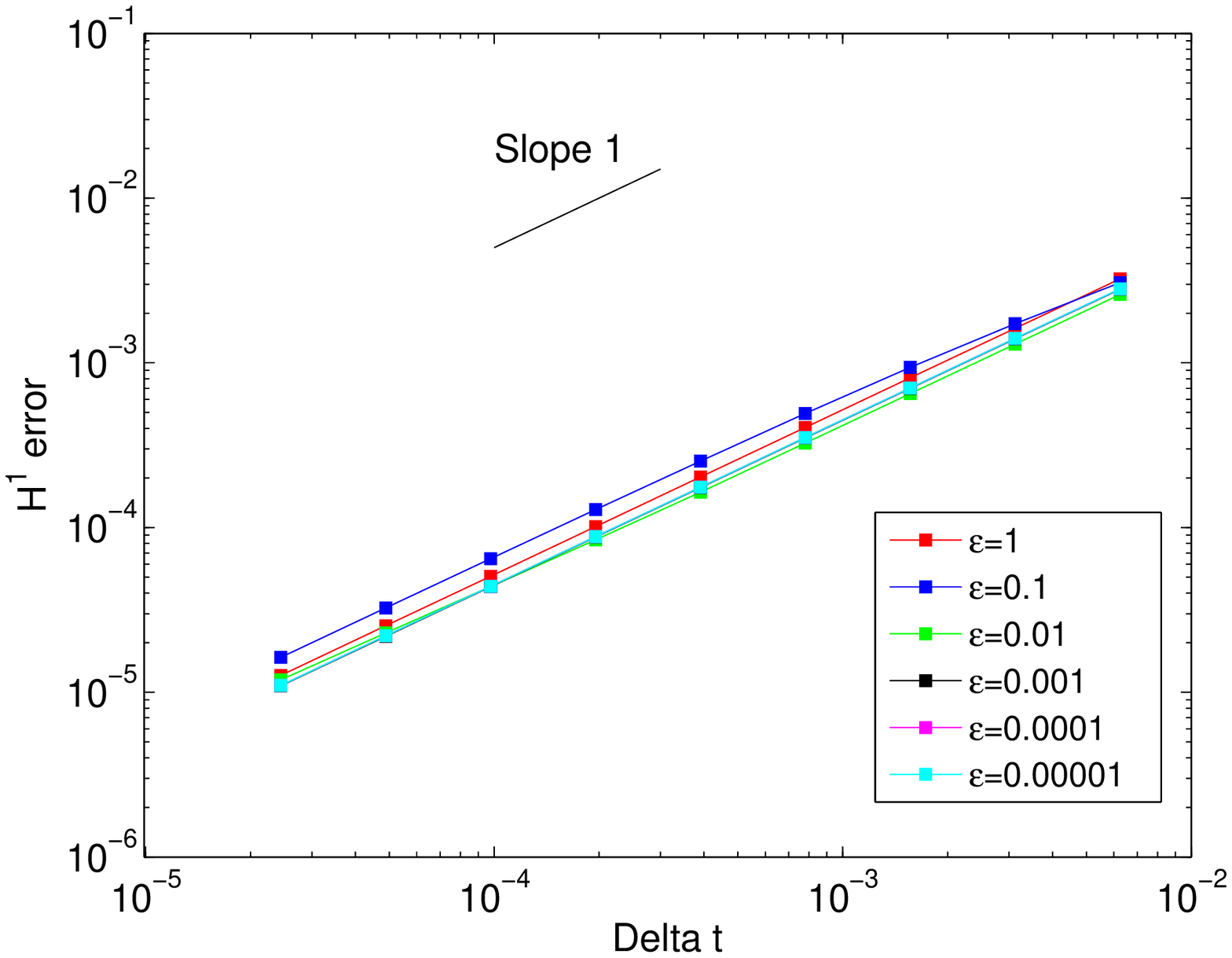}}\hspace*{-7mm}
  \subfigure[Error with respect to $\eps$]{\includegraphics[width=.717\textwidth]{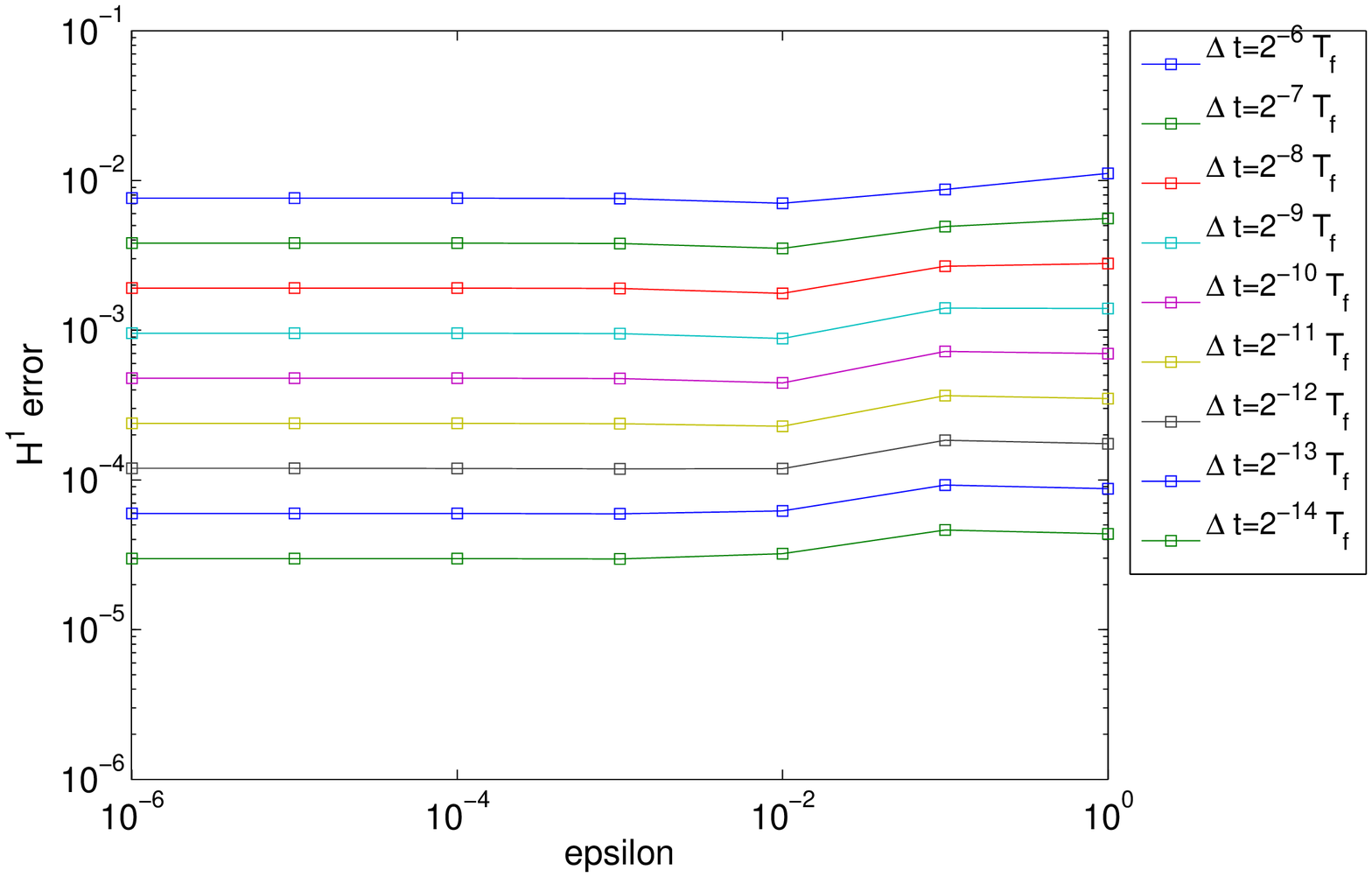}}}
    \caption{(NKG) $H^1$ relative error (in log-log scale) for the first order UA scheme with the second order initial data.}
  \label{fig5}
\end{figure}

\begin{figure}[!htbp]
  \centerline{
  \subfigure[Error with respect to $\Delta t$]{\includegraphics[width=.55\textwidth]{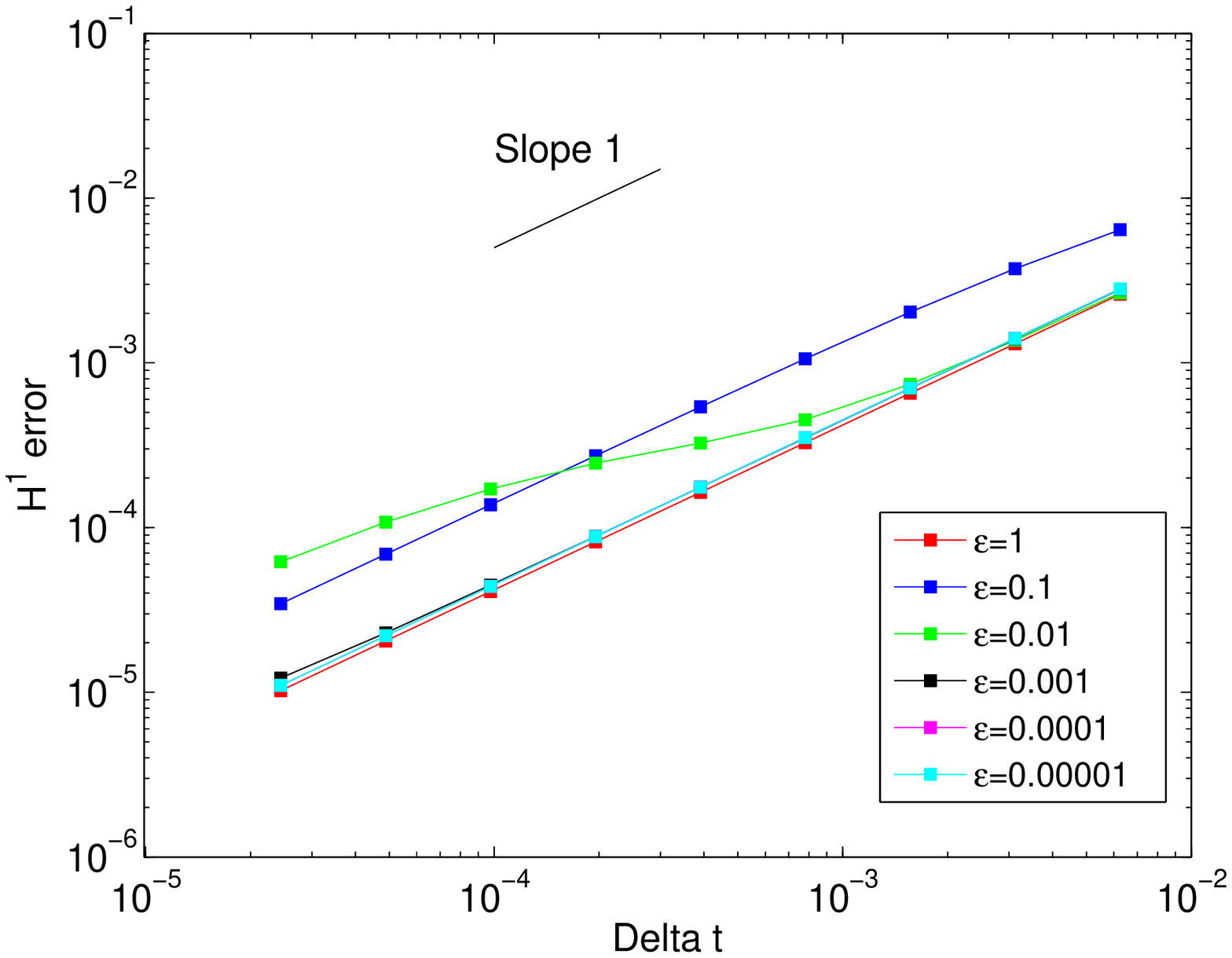}}\hspace*{-7mm}
  \subfigure[Error with respect to $\eps$]{\includegraphics[width=.717\textwidth]{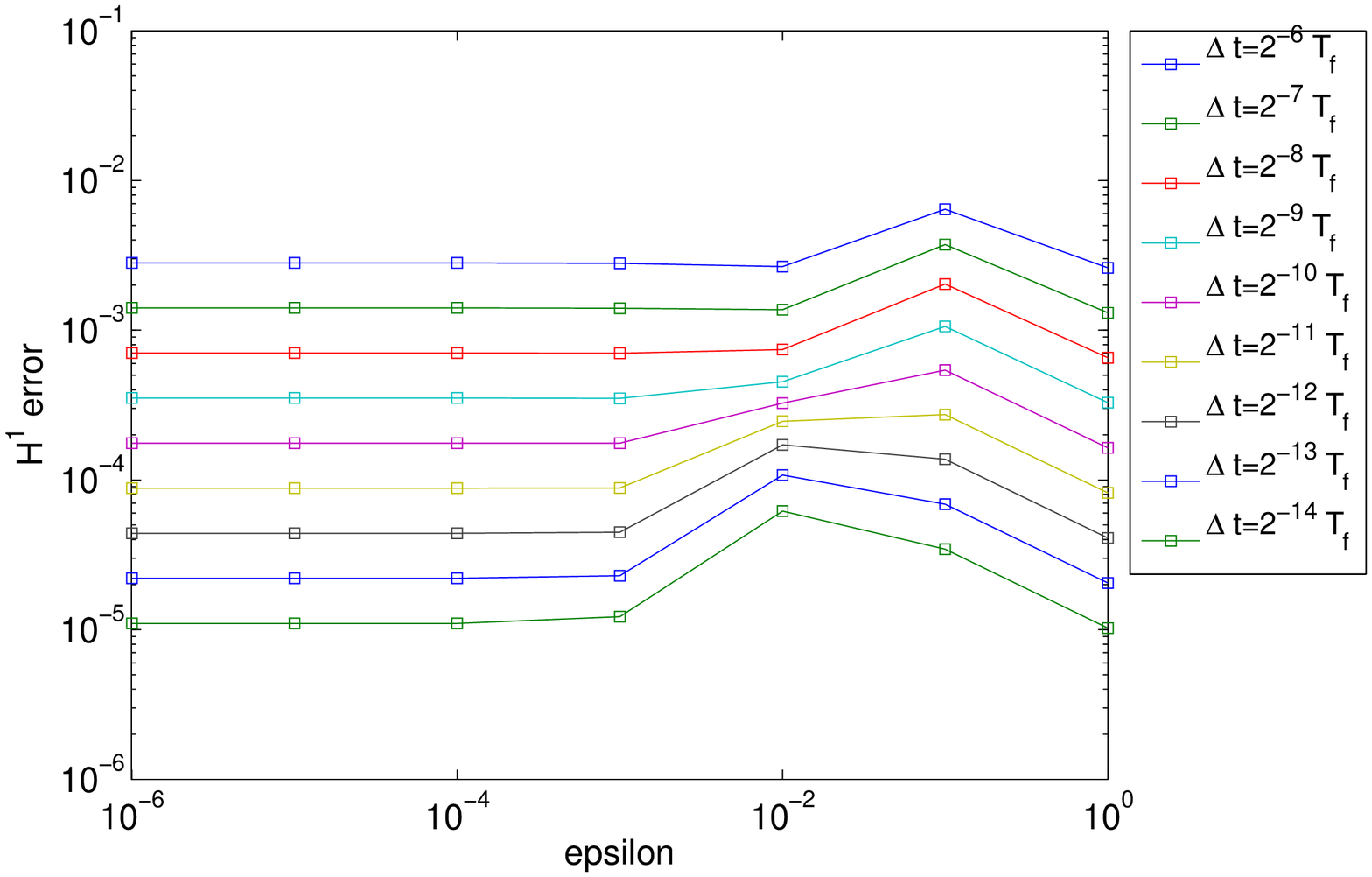}}}
    \caption{(NKG) $H^1$ relative error (in  log-log scale) for the first order UA scheme with the first order initial data.}
  \label{fig6}
\end{figure}

\begin{figure}[!htbp]
  \centerline{
  \subfigure[Error with respect to $\Delta t$]{\includegraphics[width=.55\textwidth]{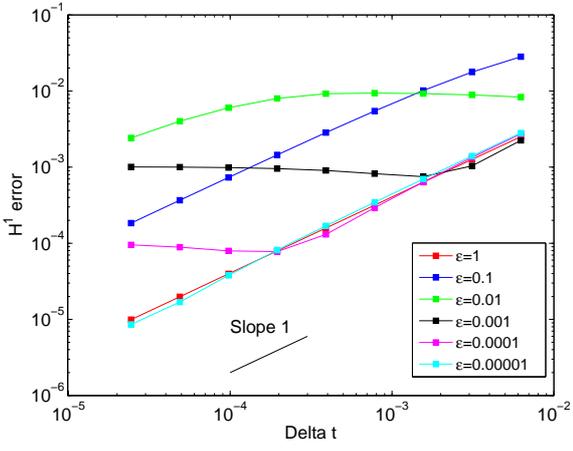}}\hspace*{-7mm}
  \subfigure[Error with respect to $\eps$]{\includegraphics[width=.717\textwidth]{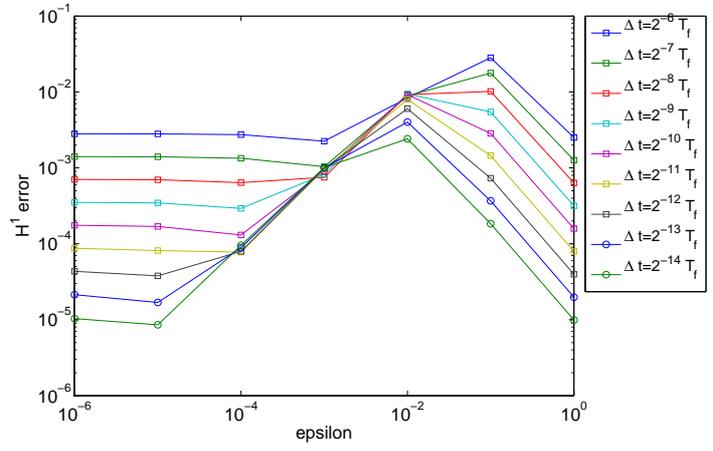}}}
    \caption{(NKG) $H^1$ relative  error (in log-log scale) for the first order UA scheme with the uncorrected initial data.}
  \label{fig7}
\end{figure}

On Figures \ref{fig5}, \ref{fig6} and \ref{fig7} we plot the $H^1$ error for our first order scheme with respect to $\Delta t$ and $\eps$. It appears  that the uniform first order accuracy is obtained for the first or second order corrected initial data, again with better results in the case of the second order initial data. If the initial data $U_0$ is taken with no correction term, the first order accuracy is lost for intermediate regimes of $\eps$ (see Figure \ref{fig7}).

\clearpage

\subsection{The nonlinear Schr\"odinger equation in a highly oscillatory regime}

In this subsection, we consider the cubic nonlinear Schr\"odinger (NLS) equation under the following form:
\be
\label{nls} 
i \partial_t u = -\frac{1}{\eps}\Delta u + \gamma(x) |u|^2 u,\qquad u(0,x)=u_0(x),
\ee
on the torus $x\in [0,a]^d$. 

For numerical simulations, our precise example is the one in dimension $d=1$ studied in \cite{grebert11} and in \cite{strobo}. We take $\gamma(x)=2  \cos(2x)$, the space domain in $x$ is $[0,2\pi]$, the initial data is
$$u_0(x)= \cos x+\sin x$$ and the final time of the simulation is $T_f=0.4$. 

As for the nonlinear Klein-Gordon case, let us first show  that \eqref{nls} fits with our general framework. The filtered wavefunction
$$\widetilde u=e^{-i\frac{t}{\eps}\Delta}u,$$
satisfies the equation
$$
i\pa_t \widetilde u=e^{-i\frac{t}{\eps}\Delta}\left(\gamma(x) \left|e^{i\frac{t}{\eps}\Delta}\widetilde u\right|^2e^{i\frac{t}{\eps}\Delta}\widetilde u\right),
$$
which is again under the form \eqref{eqftildegene} with
\be
\label{cF2}
\cF(t,\tau,u,\eps)=-ie^{-i\tau\Delta}\left(\gamma \left|e^{i\tau \Delta}u\right|^2e^{i\tau \Delta}u\right).
\ee
Here again,  it can be checked that this vector field $\cF$ satisfies Assumption \ref{mainassump}. The spectrum of the Laplace operator $-\Delta$ on the torus $x\in [0,a]^d$ is $$\left\{(2\pi/a)^2|k|^2=(2\pi/a)^2\,(k_1^2+\cdots+k_d^2) \, ; \, k \in \Z^d\right\}\subset (2\pi/a)^2\N$$
so that $\tau\mapsto e^{i\tau \Delta}$ is periodic, with period $P=\frac{a^2}{2\pi}$, and $\cF$ is also periodic w.r.t. $\tau$.

\bs
In order to validate our approach, we now proceed with similar numerical tests as in the case of the NKG equation. The reference solution is computed as follows. For $\eps\geq 10^{-2}$, we use the Yoshida fourth order splitting method \cite{yoshida} with $\Delta x=2\pi/128$, $\Delta t=\eps\,T_f/32768$. For smaller values of $\eps$, we rather use our second order scheme, with the following parameters: $\Delta x=2\pi/128\approx 0.05$, $\Delta t=2\pi/512000\approx 1.2\times 10^{-5}$, $\Delta \tau= 2\pi/4096\approx 1.5\times 10^{-3}$.

\begin{figure}[!htbp]
  \centerline{
  \subfigure[Error with respect to $\Delta t$]{\includegraphics[width=.55\textwidth]{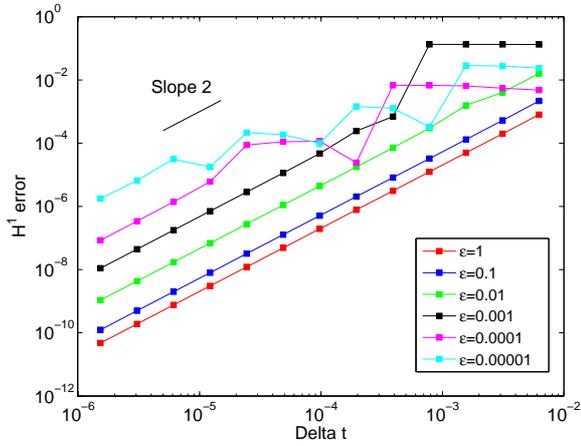}}\hspace*{-7mm}
  \subfigure[Error with respect to $\eps$]{\includegraphics[width=.717\textwidth]{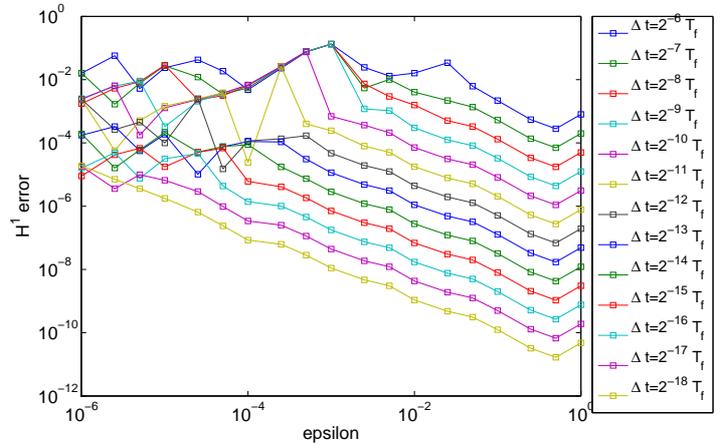}}}
    \caption{(NLS case) $H^1$ relative error for the Strang splitting scheme.}
  \label{fig1000}
\end{figure}

\bs
We first plot on Figure \ref{fig1000} the $H^1$ error between the numerical solution computed with the standard Strang splitting scheme for NLS (with a fixed, large enough, number of points in $x$, $N_x=128$) and the reference solution. It appears again that the error behaves asymptotically like $C\frac{\Delta t^2}{\eps}$, where $C$ does not depend on $\Delta t$ and $\eps$.
%$$\|u^{ref}(t_{final},\cdot)-u^{Strang}(t_{final},\cdot)\|_{L^2}\sim C\frac{\Delta t^2}{\eps}.$$

\begin{figure}[!htbp]
  \centerline{\includegraphics[width=.60\textwidth]{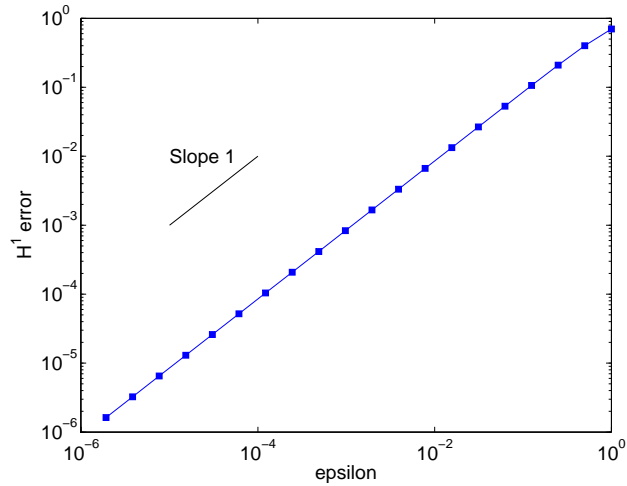}}
    \caption{(NLS case) $H^1$ relative error between the reference solution and the limiting averaged model.}
  \label{fig1009}
\end{figure}

As $\eps\to 0$, the solution of \eqref{nls} behaves as
\be
\label{estimaveNLS}
\|u(t,x)-e^{i\frac{t}{\eps}\Delta}w(t,x)\|\leq C\eps,
\ee
where $w$ solves the averaged equation
\be
i\pa_t \widetilde w=\frac{1}{P}\int_0^P e^{-i\tau \Delta}\left(\gamma(x) \left|e^{i\tau\Delta}w\right|^2e^{i\tau\Delta}w\right)d\tau.
\label{nlslimit}
\ee
On Figure \ref{fig1009}, we illustrate this asymptotic behavior by plotting the error between the solution of the limiting averaged model and the reference solution.

\begin{figure}[!htbp]
  \centerline{
  \subfigure[Error with respect to $N_x$]{\includegraphics[width=.55\textwidth]{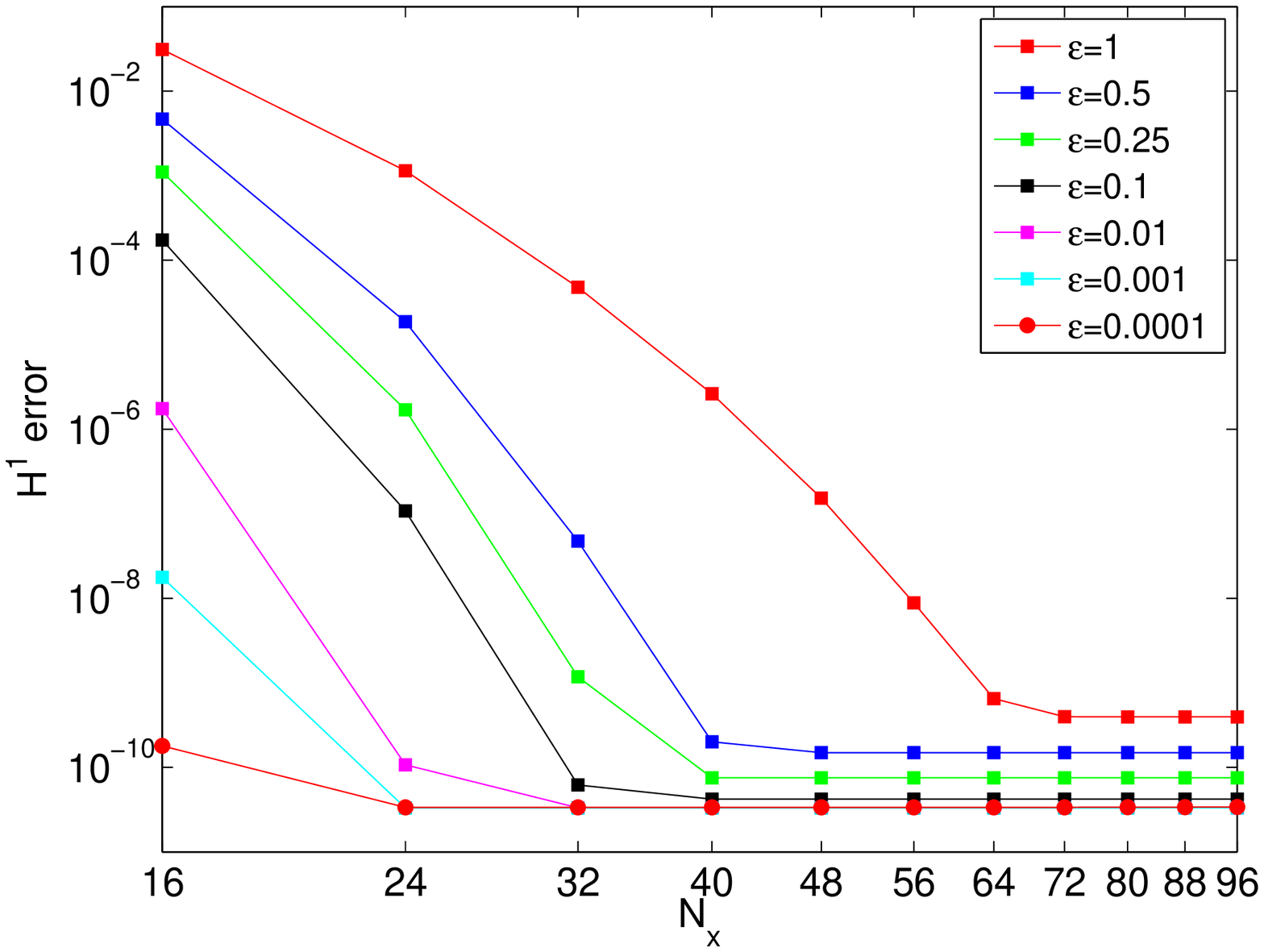}}\hspace*{-7mm}
  \subfigure[Error with respect to $N_\tau$]{\includegraphics[width=.55\textwidth]{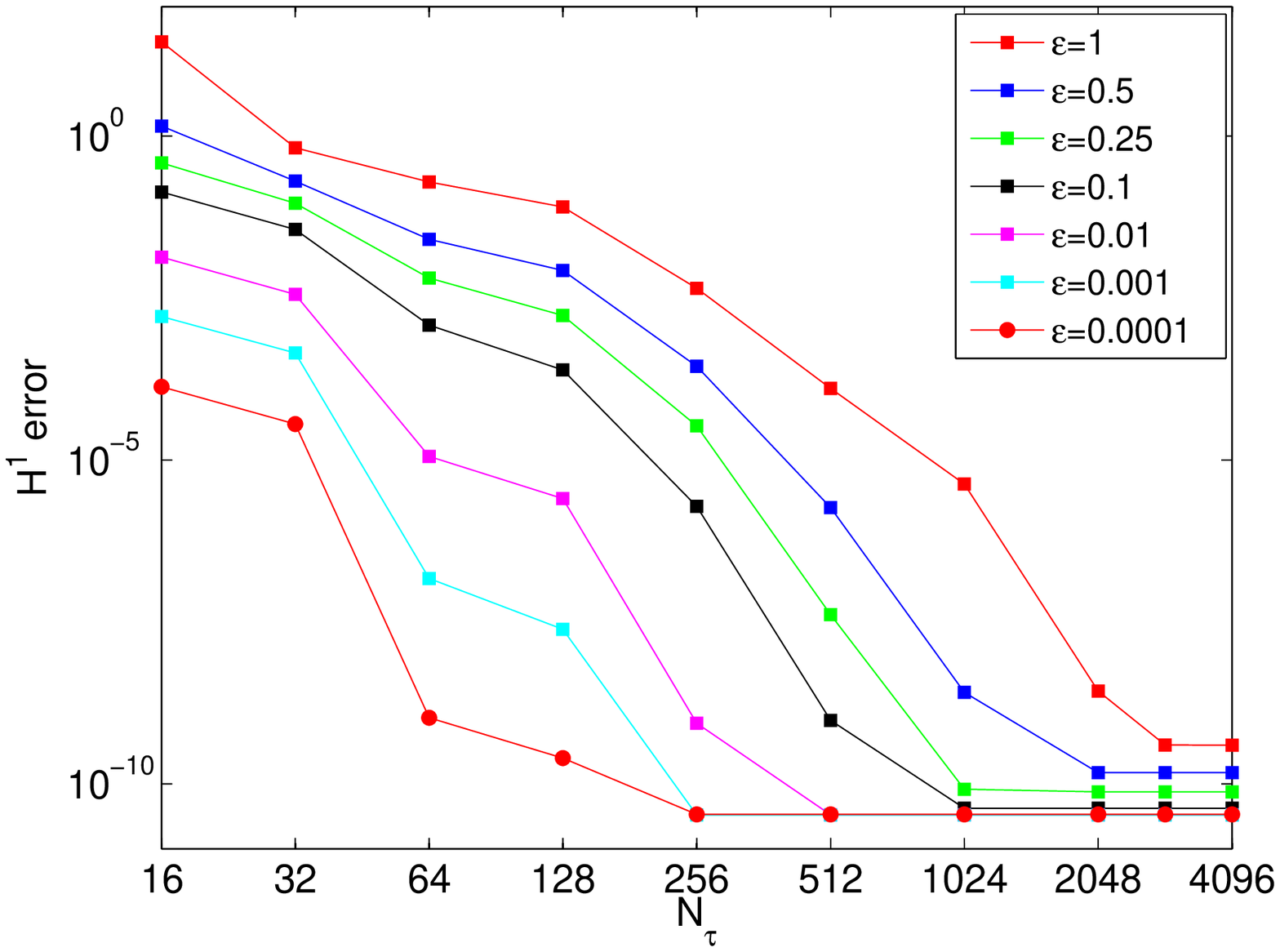}}}
    \caption{(NLS case) $H^1$ relative error for the second order UA scheme with the third order initial data.}
  \label{fig1008}
\end{figure}

\bs
Let us now characterize the behavior of our uniformly accurate numerical schemes with respect to the numerical parameters. We first plot on Figure \ref{fig1008} the $H^1$ error with respect to the number of grid points $N_x$ in the $x$ variable (left figure, for which we take $N_\tau=2048$ and $\Delta t=2\times 10^{-5}$) and with respect to the number of grid points $N_\tau$ in the $\tau$ variable (right figure, for which we take $N_x=64$ and $\Delta t=2\times 10^{-5}$). As in the NKG case, we observe that our scheme has a spectral accuracy in $x$ and in $\tau$. However, two main differences can be observed between the NKG and the NLS cases. First, the error in $N_x$ becomes smaller as $\eps$ decreases. Second, we have to take much smaller steps $\Delta \tau$ in the NLS case than in the NKG case. This is due to the operators $e^{i\tau \Delta}$ in the function $\cF$: the NLS problem is  stiffer than the NKG problem and involves high frequencies in the $\tau$ variable. In the sequel, we fix $N_x=64$ and $N_\tau= 2048$.

\begin{figure}[!htbp]
  \centerline{
  \subfigure[Error with respect to $\Delta t$]{\includegraphics[width=.55\textwidth]{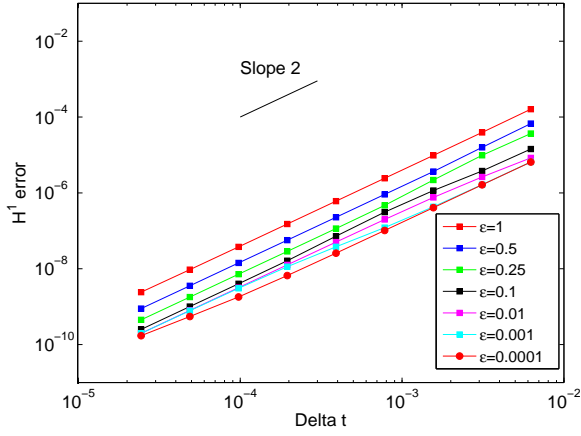}}\hspace*{-7mm}
  \subfigure[Error with respect to $\eps$]{\includegraphics[width=.717\textwidth]{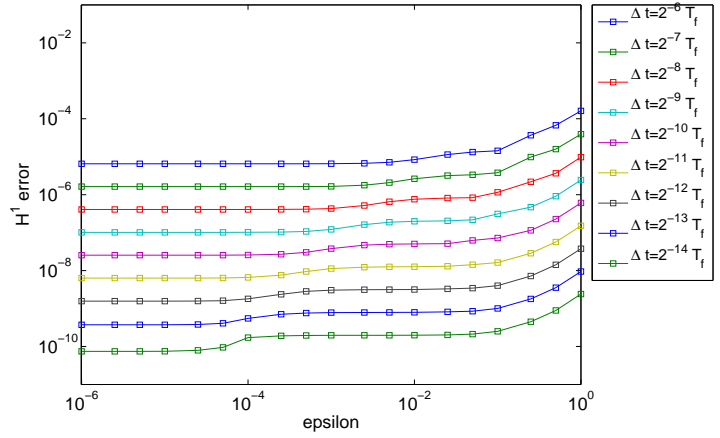}}}
    \caption{(NLS case) $H^1$ relative error for the second order UA scheme with the third order initial data.}
  \label{fig1001}
\end{figure}

\begin{figure}[!htbp]
  \centerline{
  \subfigure[Error with respect to $\Delta t$]{\includegraphics[width=.55\textwidth]{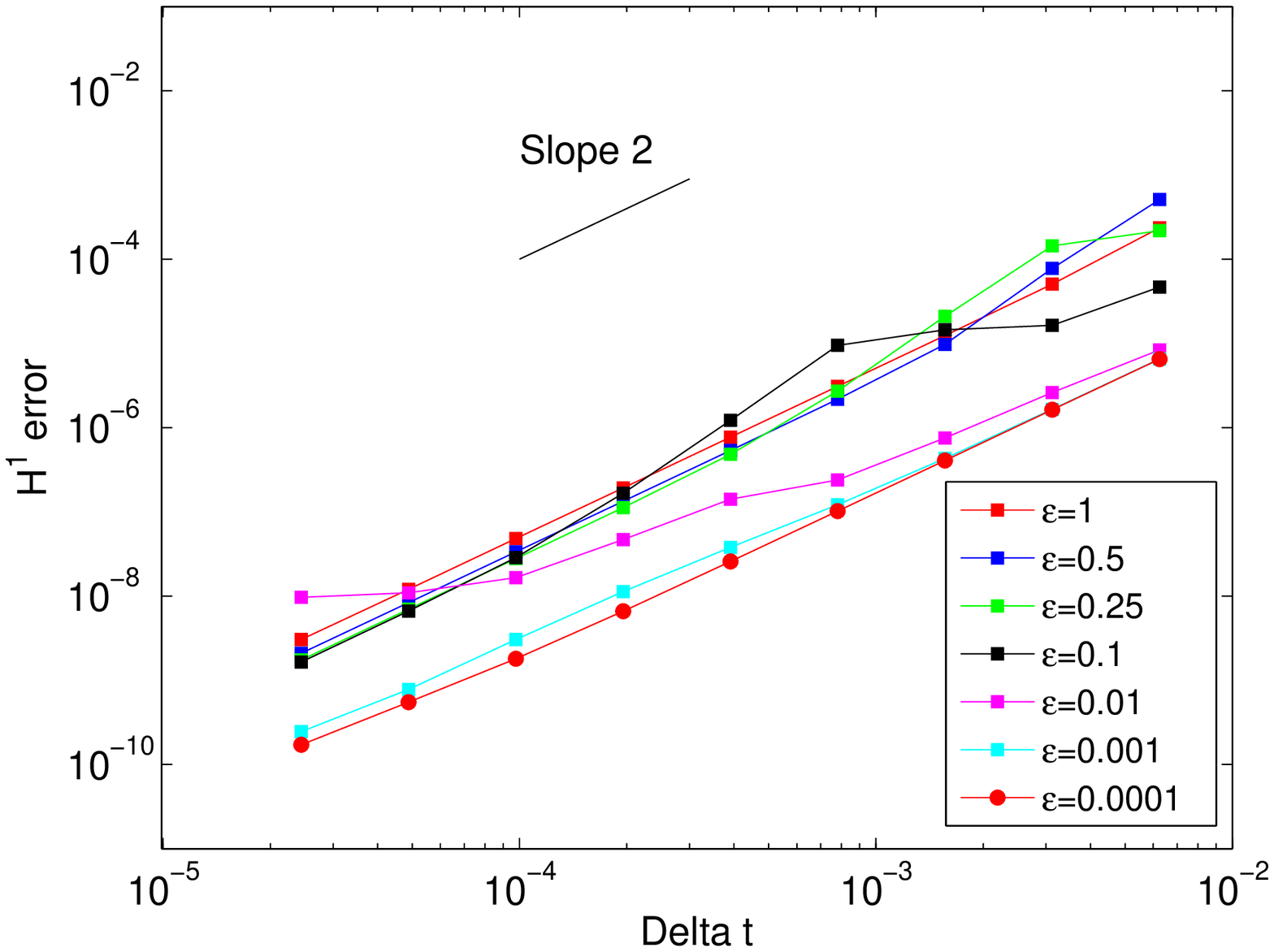}}\hspace*{-7mm}
  \subfigure[Error with respect to $\eps$]{\includegraphics[width=.717\textwidth]{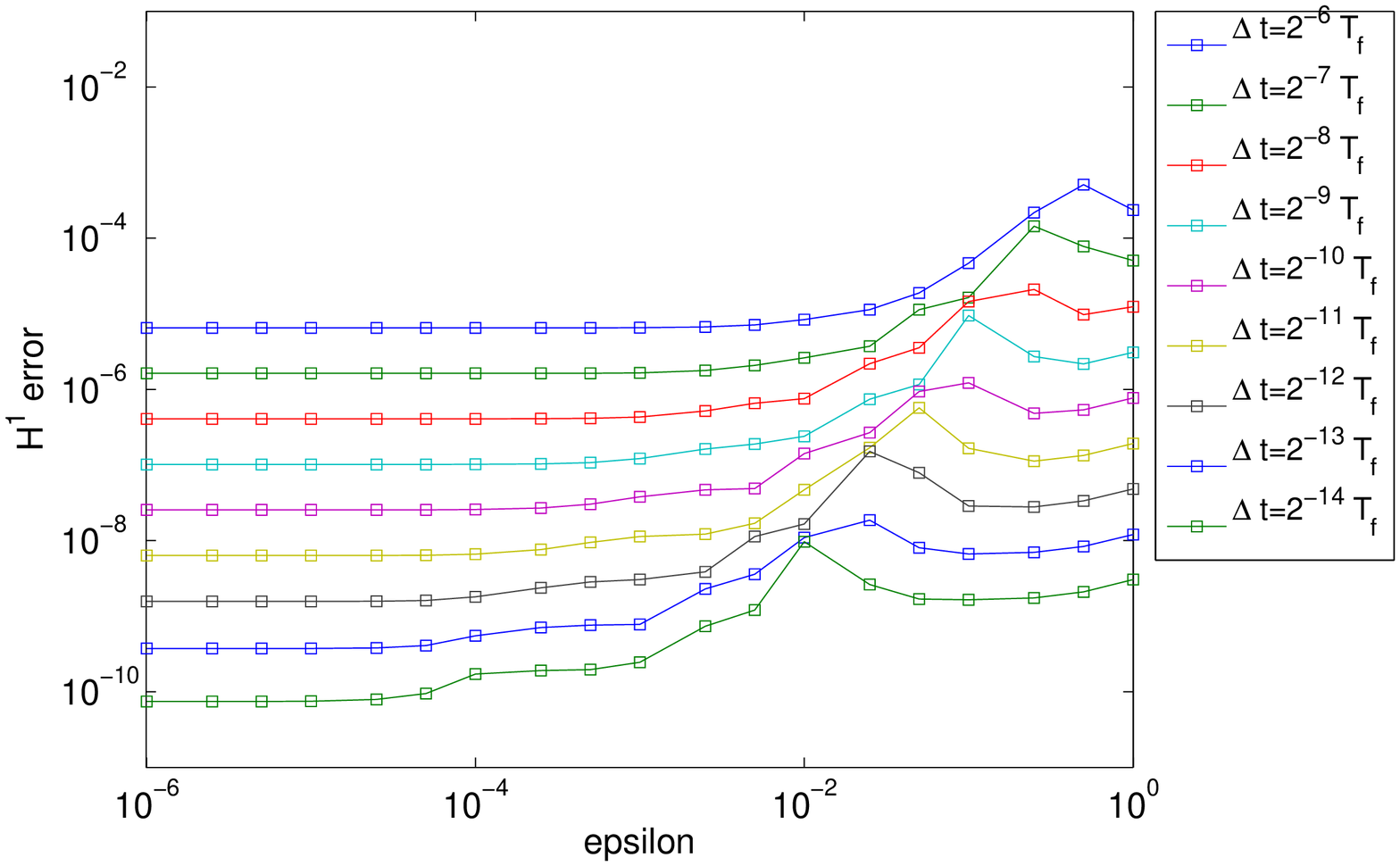}}}
    \caption{(NLS case) $H^1$ relative error for the second order UA scheme with the second order initial data.}
  \label{fig1002}
\end{figure}

\begin{figure}[!htbp]
  \centerline{
  \subfigure[Error with respect to $\Delta t$]{\includegraphics[width=.55\textwidth]{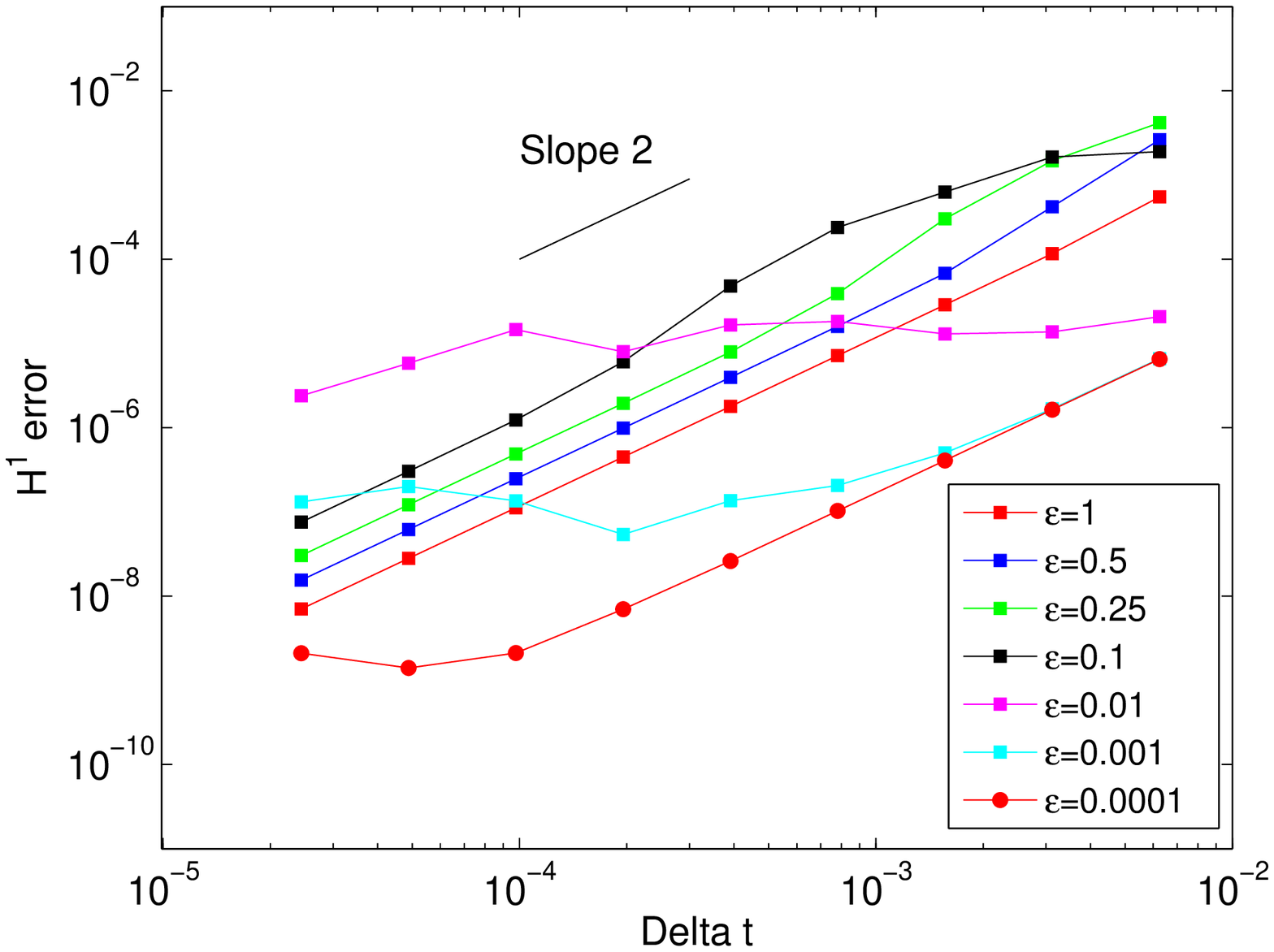}}\hspace*{-7mm}
  \subfigure[Error with respect to $\eps$]{\includegraphics[width=.717\textwidth]{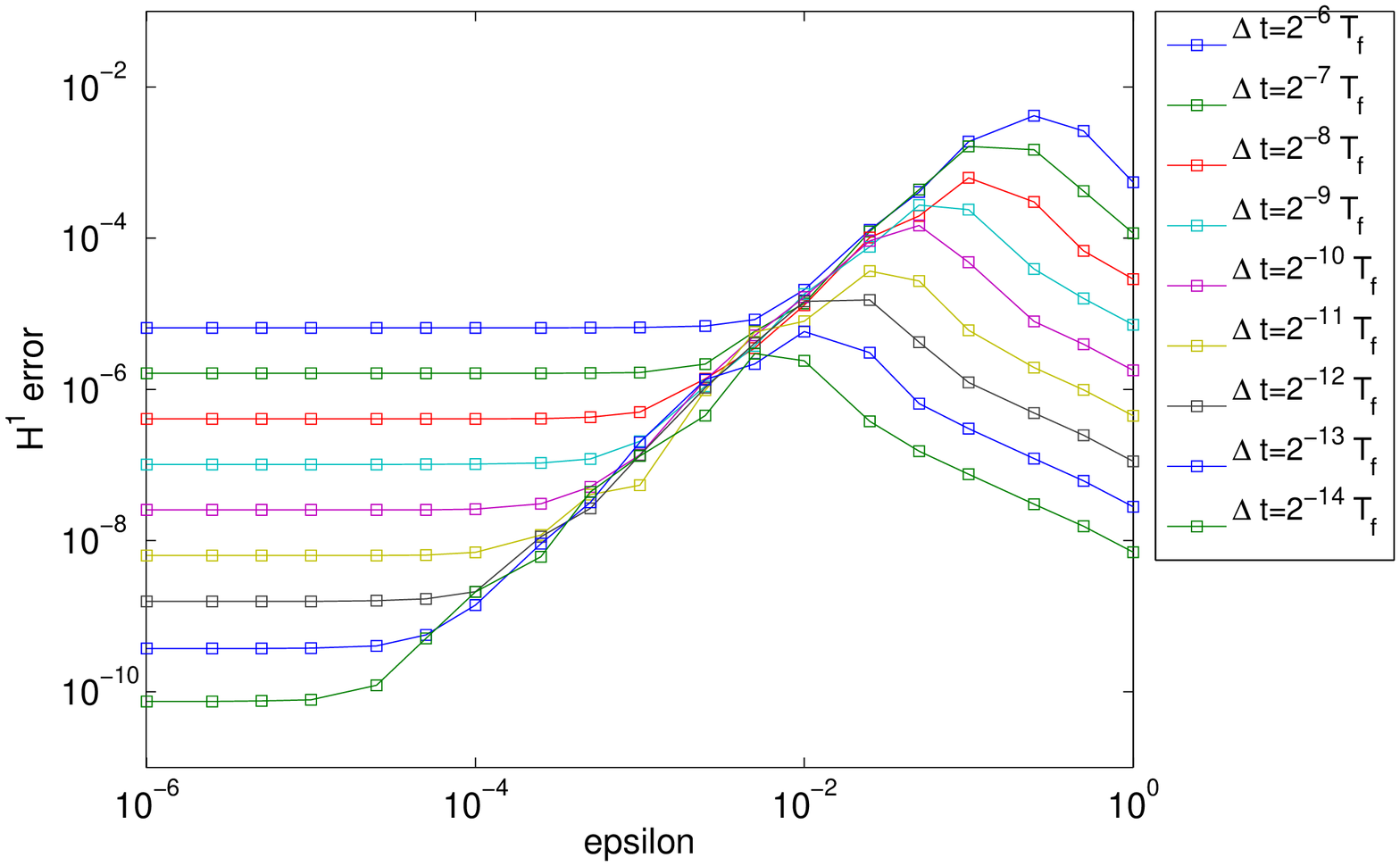}}}
    \caption{(NLS case) $H^1$ relative error for the second order UA scheme with the first order initial data.}
  \label{fig1003}
\end{figure}

\begin{figure}[!htbp]
  \centerline{
  \subfigure[Error with respect to $\Delta t$]{\includegraphics[width=.55\textwidth]{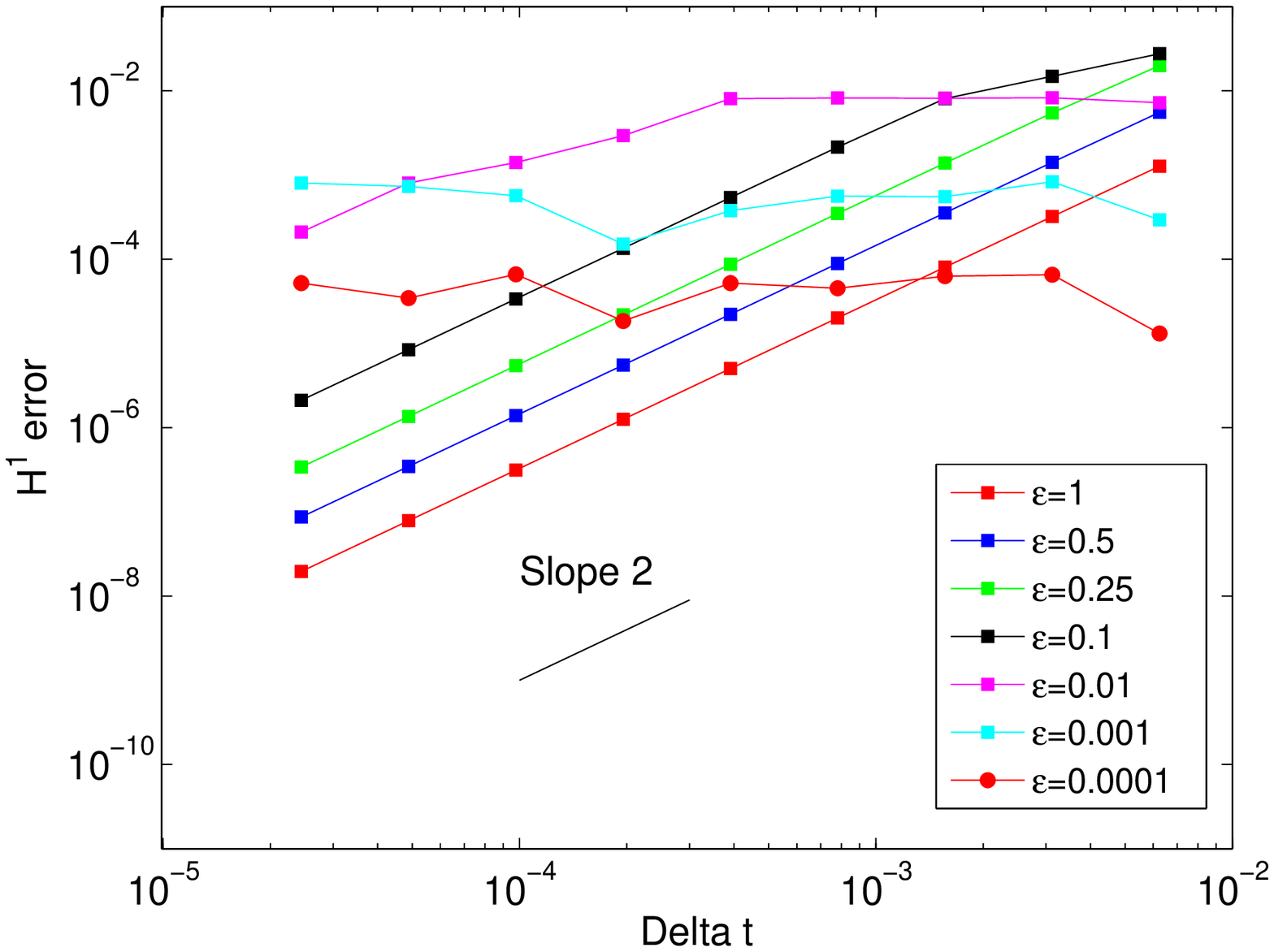}}\hspace*{-7mm}
  \subfigure[Error with respect to $\eps$]{\includegraphics[width=.717\textwidth]{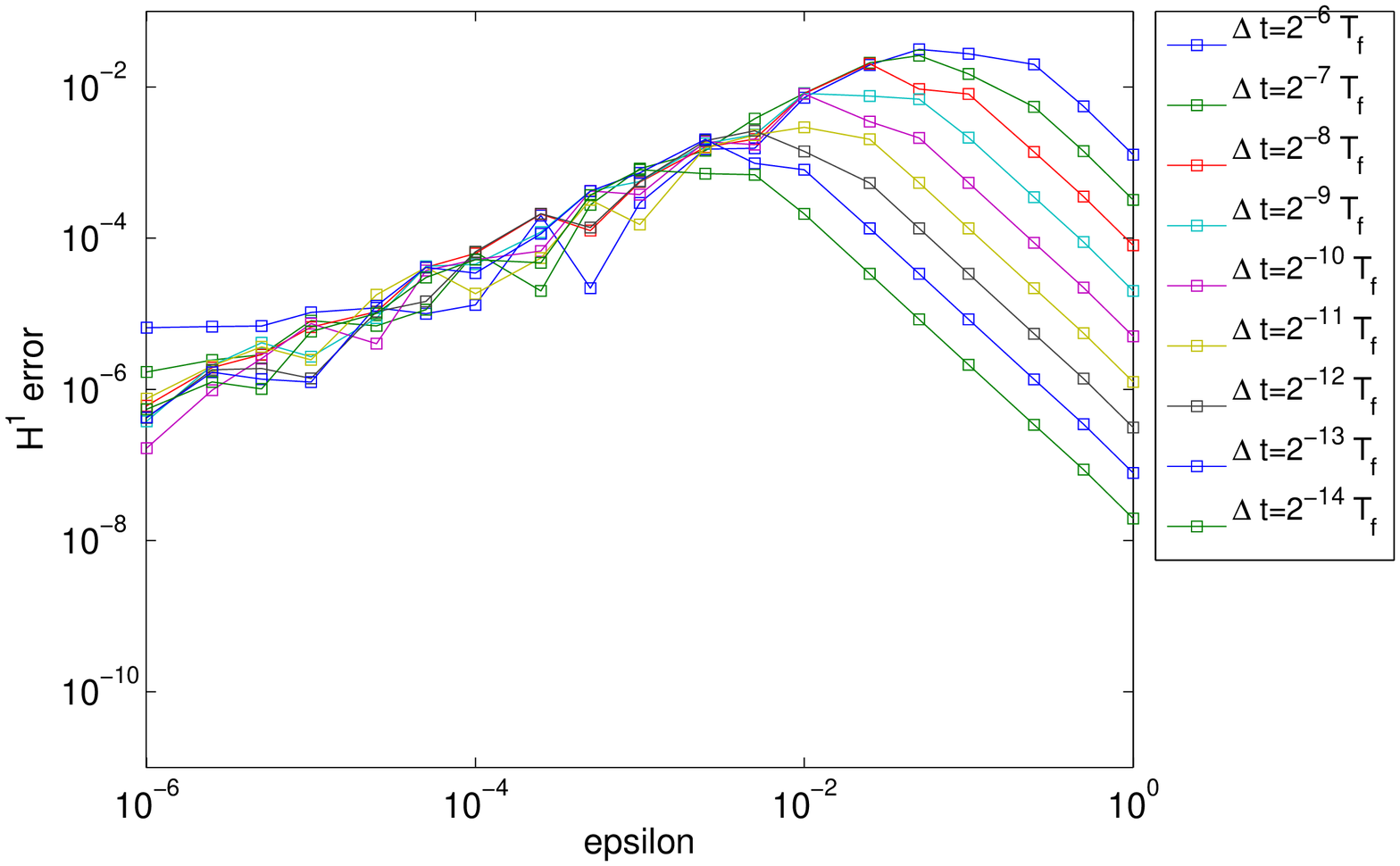}}}
    \caption{(NLS case) $H^1$ relative error for the second order UA scheme with the uncorrected intial data.}
  \label{fig1004}
\end{figure}

We now observe the behavior of our schemes with respect to the time step $\Delta t$. On Figures \ref{fig1001} and \ref{fig1002}, we plot the error between the reference solution and the numerical solution of our second order numerical scheme, for the third order and the second order initial data $U_0$. As in the NKG case, our numerical scheme displays a uniform second order error, with a slightly better result in the case of the third order initial data. If the initial data is not taken with enough terms, the uniform accuracy is lost for intermediate regimes, see Figures \ref{fig1003} (initial data with first order correction) and \ref{fig1004} (initial data with no correction).

\begin{figure}[!htbp]
  \centerline{
  \subfigure[Error with respect to $\Delta t$]{\includegraphics[width=.55\textwidth]{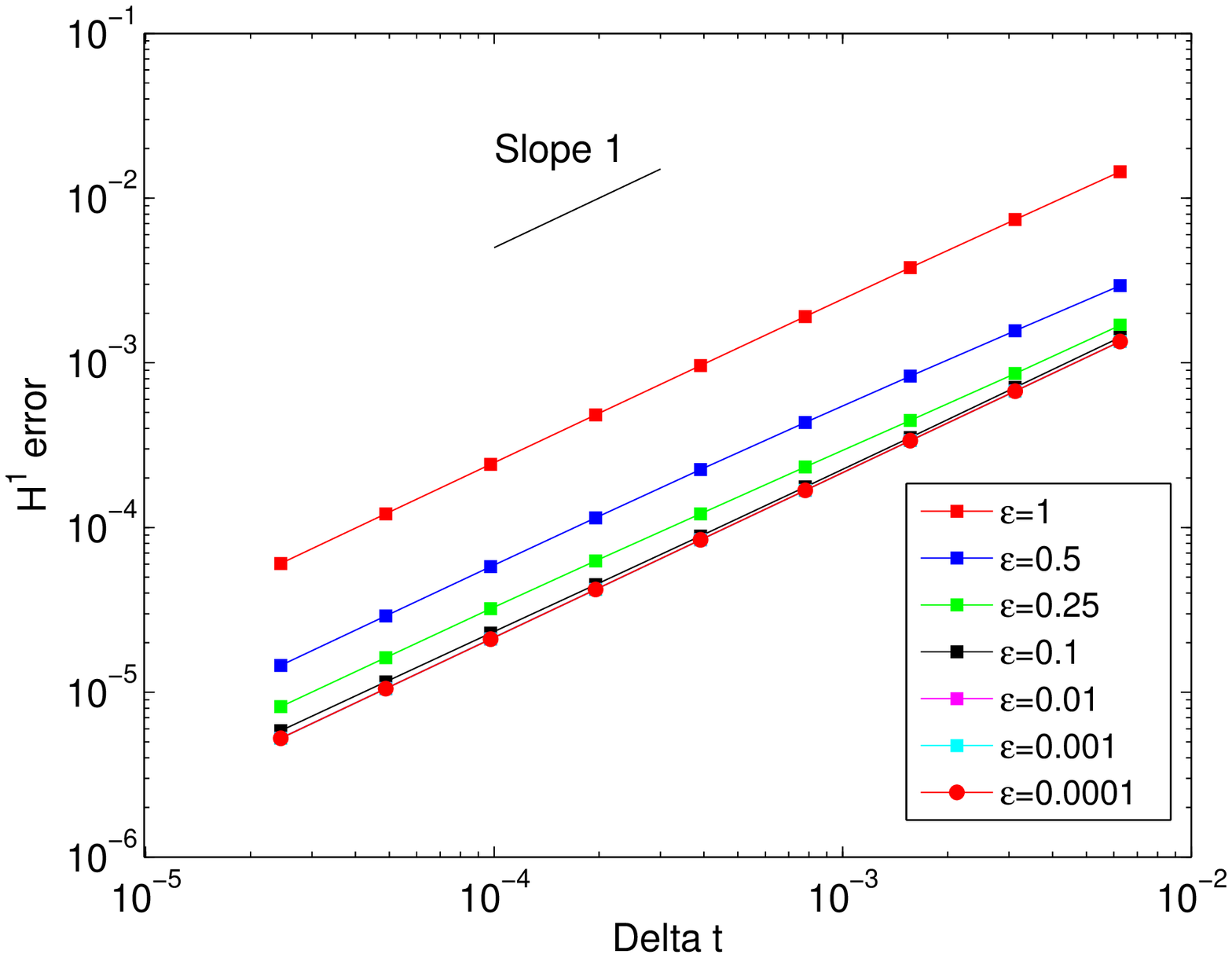}}\hspace*{-7mm}
  \subfigure[Error with respect to $\eps$]{\includegraphics[width=.717\textwidth]{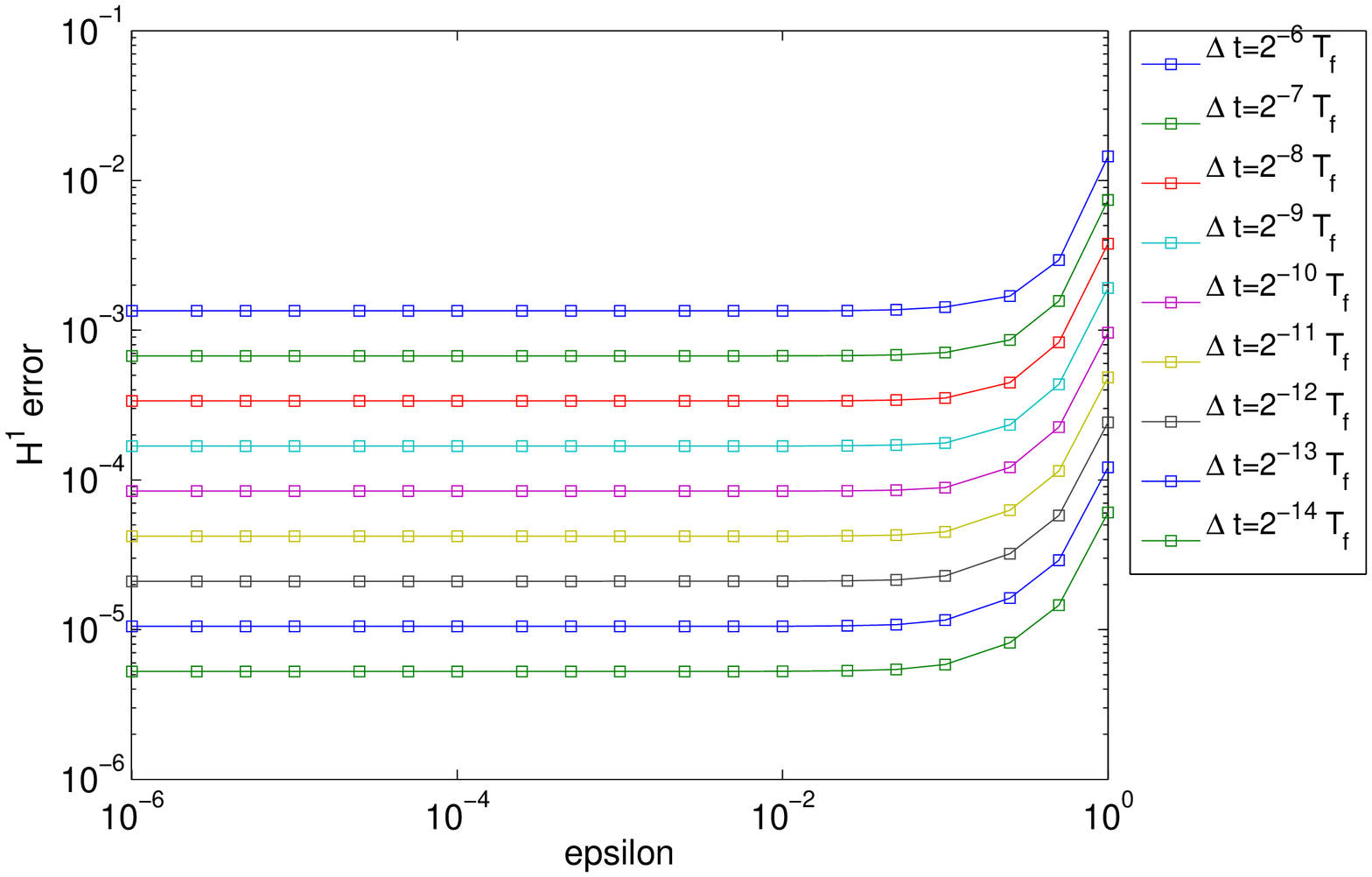}}}
    \caption{(NLS case) $H^1$ relative error for the first order UA scheme with the second order initial data.}
  \label{fig1005}
\end{figure}

\begin{figure}[!htbp]
  \centerline{
  \subfigure[Error with respect to $\Delta t$]{\includegraphics[width=.55\textwidth]{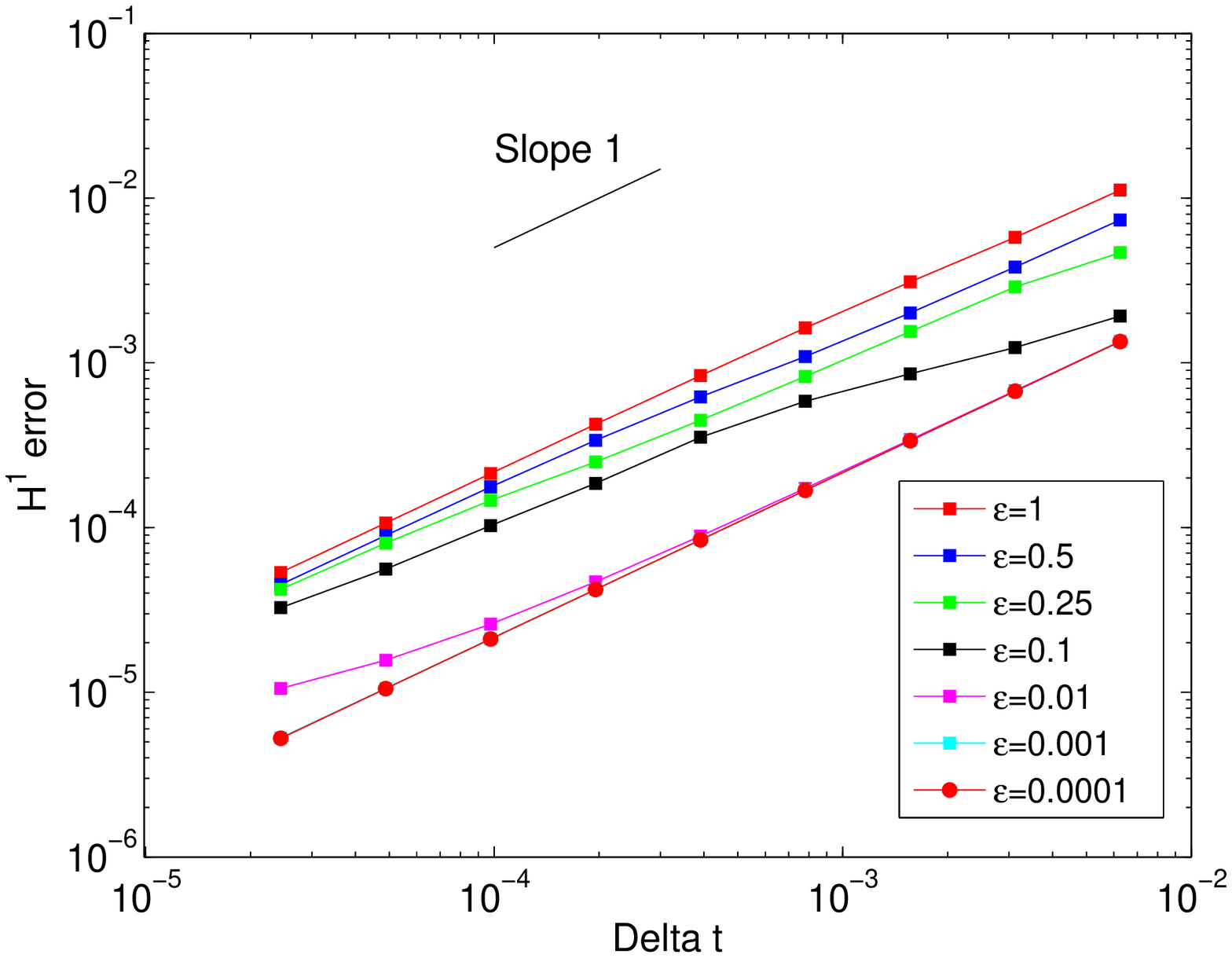}}\hspace*{-7mm}
  \subfigure[Error with respect to $\eps$]{\includegraphics[width=.717\textwidth]{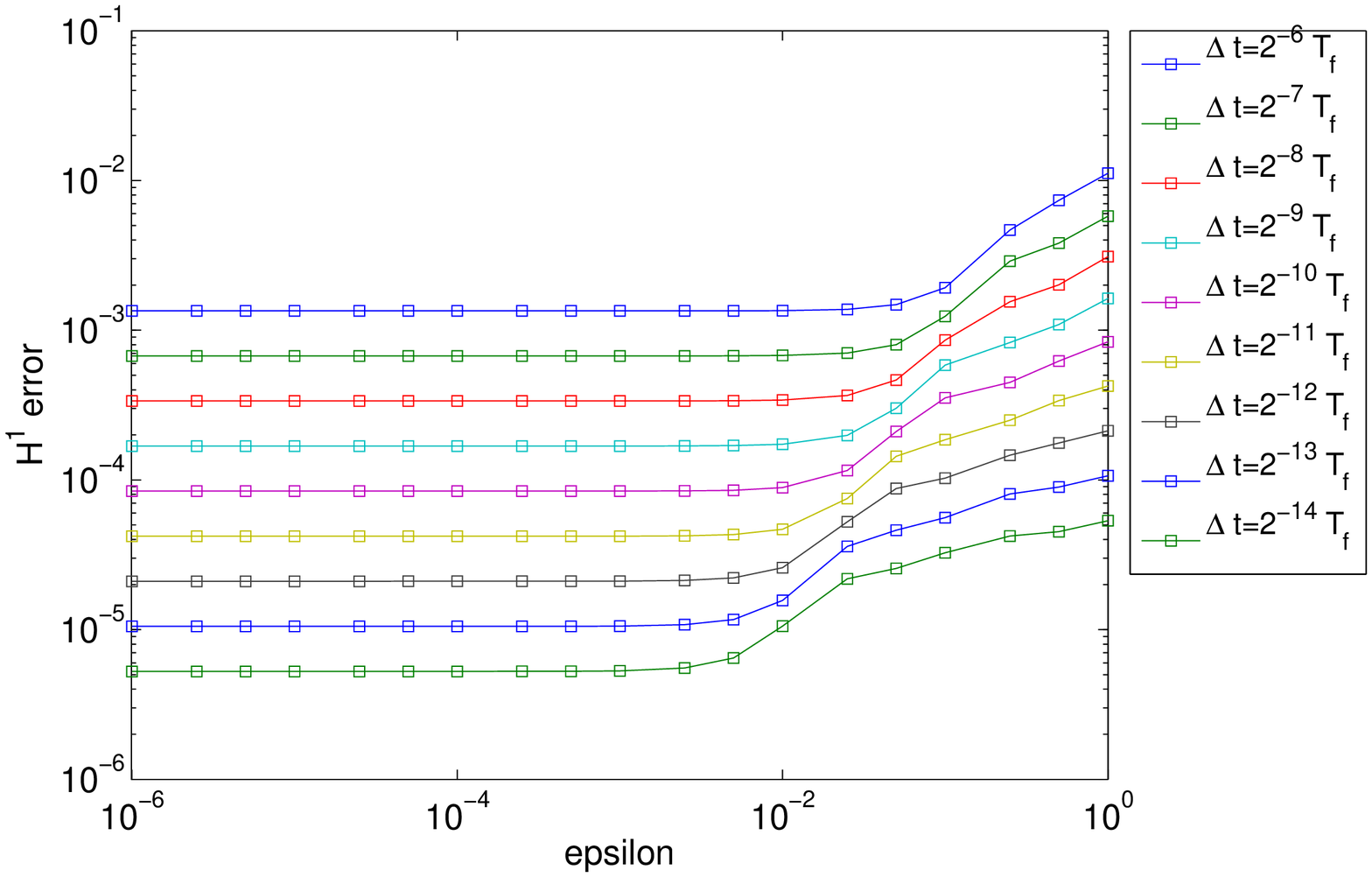}}}
    \caption{(NLS case) $H^1$ relative error for the first order UA scheme with the first order initial data.}
  \label{fig1006}
\end{figure}

\begin{figure}[!htbp]
  \centerline{
  \subfigure[Error with respect to $\Delta t$]{\includegraphics[width=.55\textwidth]{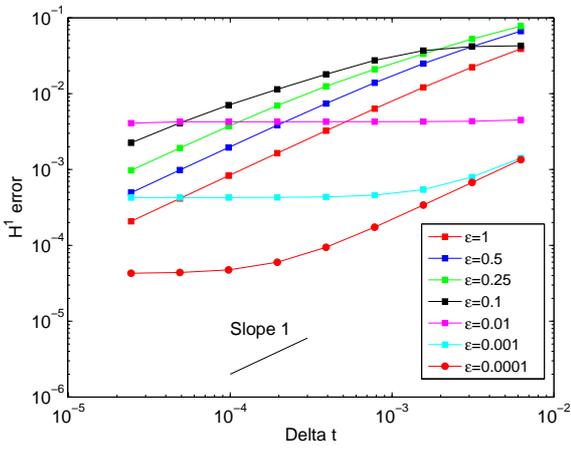}}\hspace*{-7mm}
  \subfigure[Error with respect to $\eps$]{\includegraphics[width=.717\textwidth]{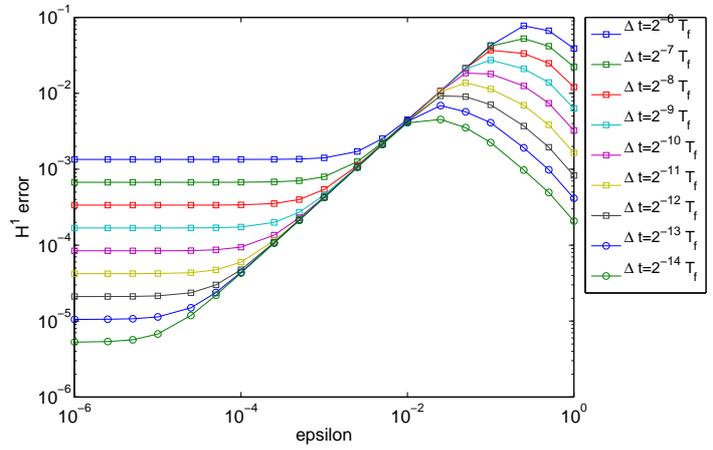}}}
    \caption{(NLS case) $H^1$ relative error for the first order UA scheme with the uncorrected initial data.}
  \label{fig1007}
\end{figure}

\bs
On Figures \ref{fig1005}, \ref{fig1006} and \ref{fig1007}, we plot the $H^1$ error for our first order numerical scheme, respectively in the three following cases: initial data with second order correction, first order correction, and with no correction. As expected, the error is uniform with respect to $\eps$ in the first  two  cases, and loses its uniformity if the initial data is taken with no correction.

\clearpage

\begin{figure}[!htbp]
  \centerline{
  \subfigure[With the third order initial data]{\includegraphics[width=.717\textwidth]{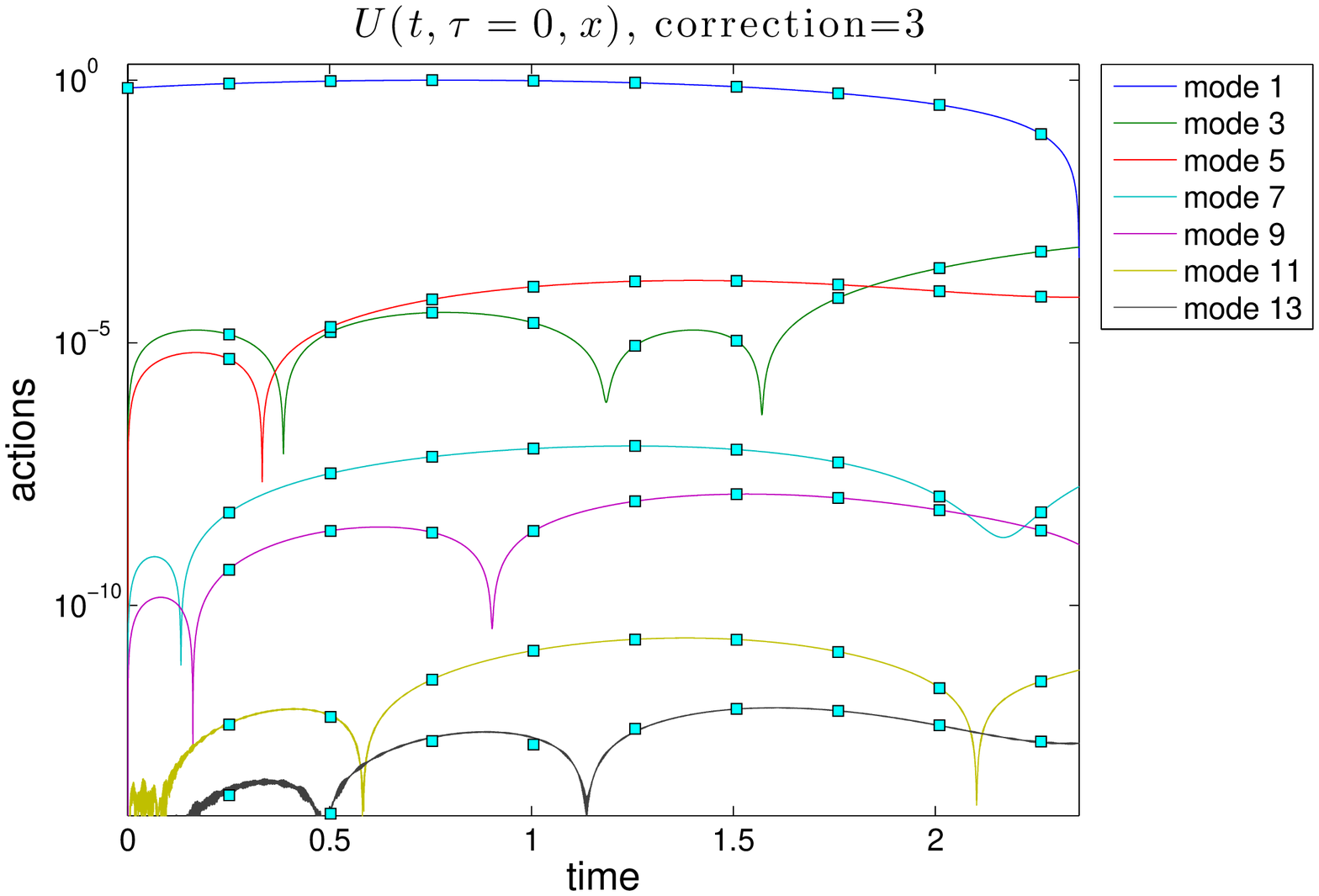}}\hspace*{-10mm}
  \subfigure[With the second order initial data]{\includegraphics[width=.717\textwidth]{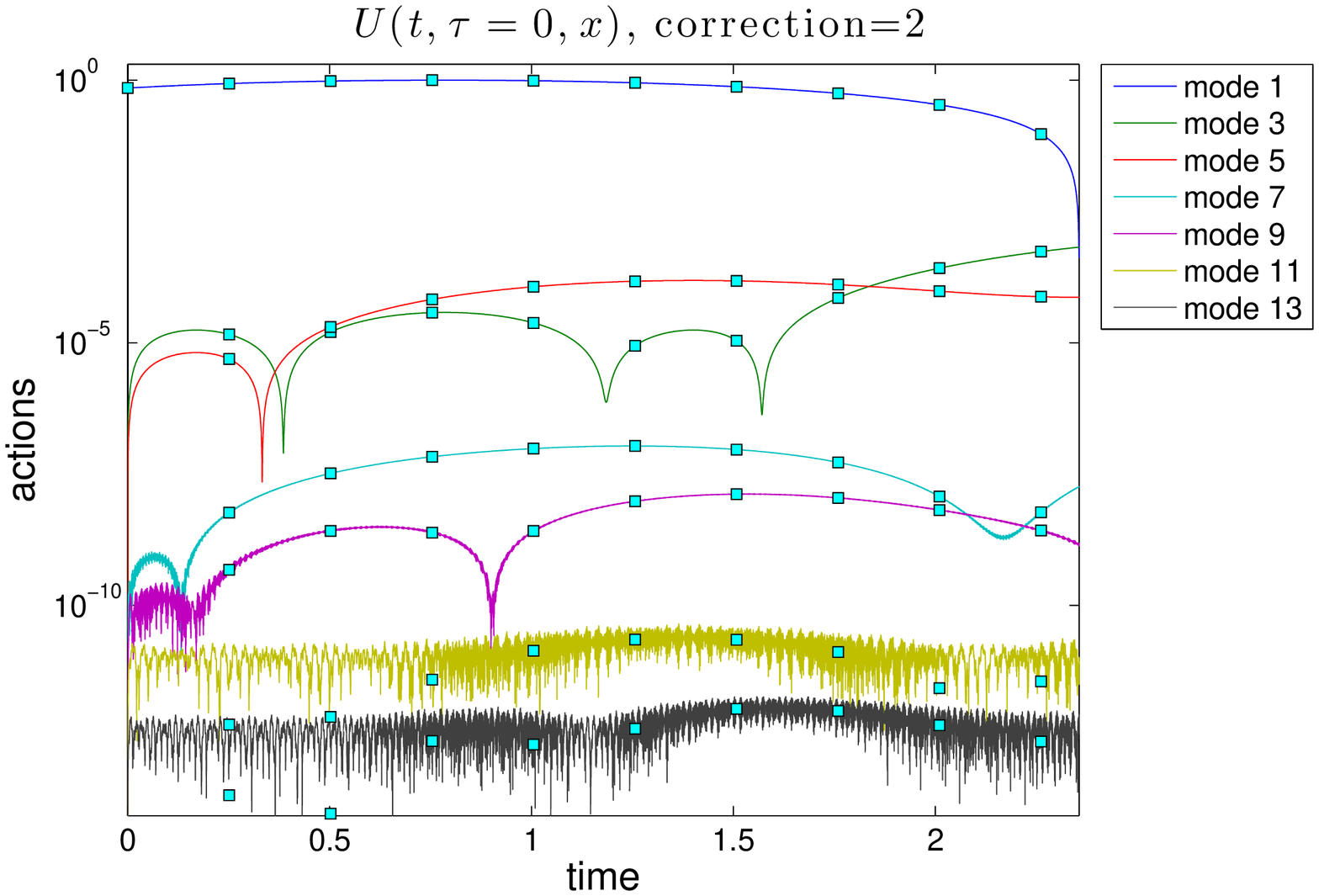}}}
    \caption{(NLS) Time evolution of the first Fourier modes in $x$ of the function $U(t, \tau=0,x)$, in the log-scale. At blue squares is plotted the reference solution at some stroboscopic points.}
  \label{fig1010}
\end{figure}

\begin{figure}[!htbp]
  \centerline{
  \subfigure[With the first order initial data]{\includegraphics[width=.717\textwidth]{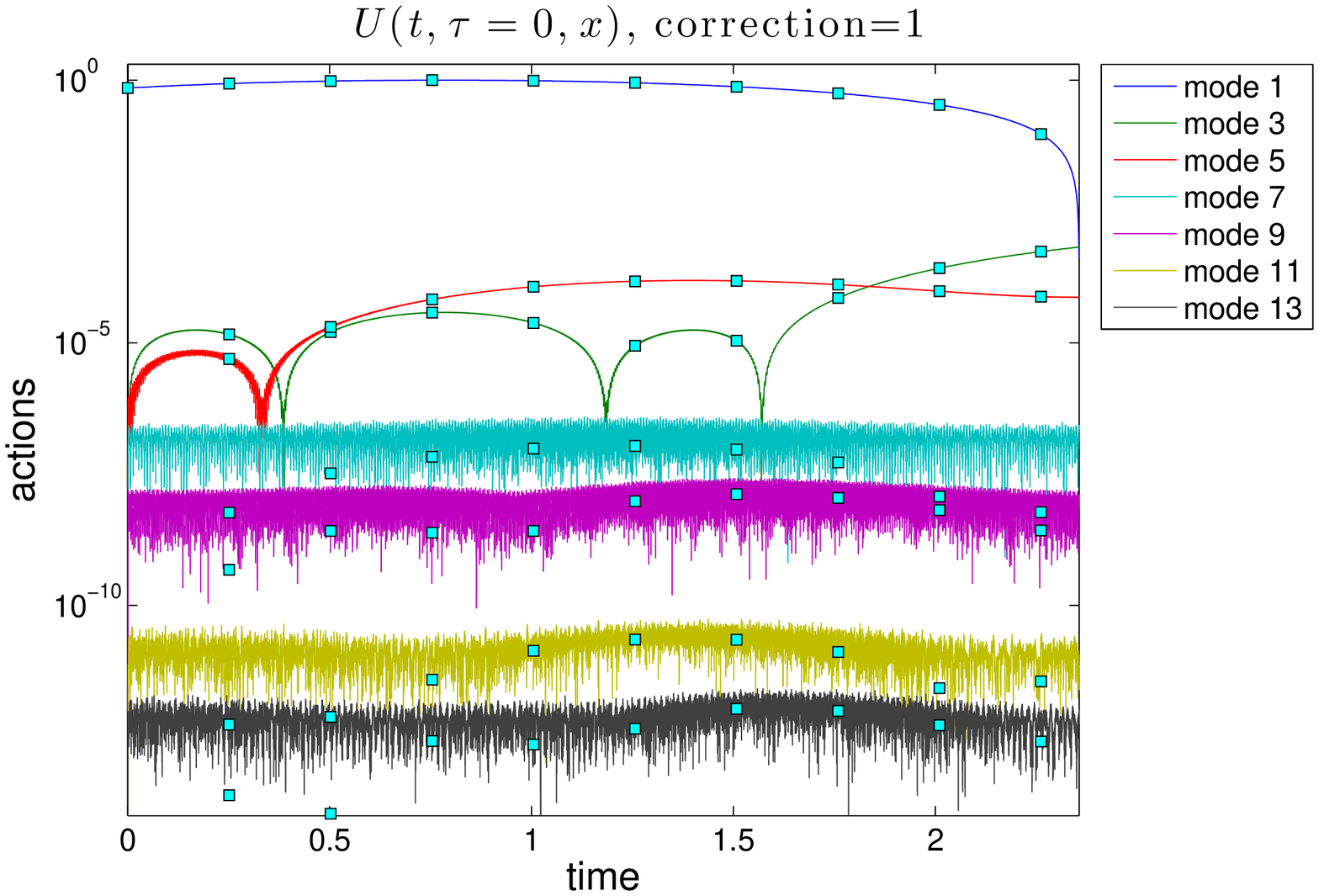}}\hspace*{-10mm}
  \subfigure[With the uncorrected initial data]{\includegraphics[width=.717\textwidth]{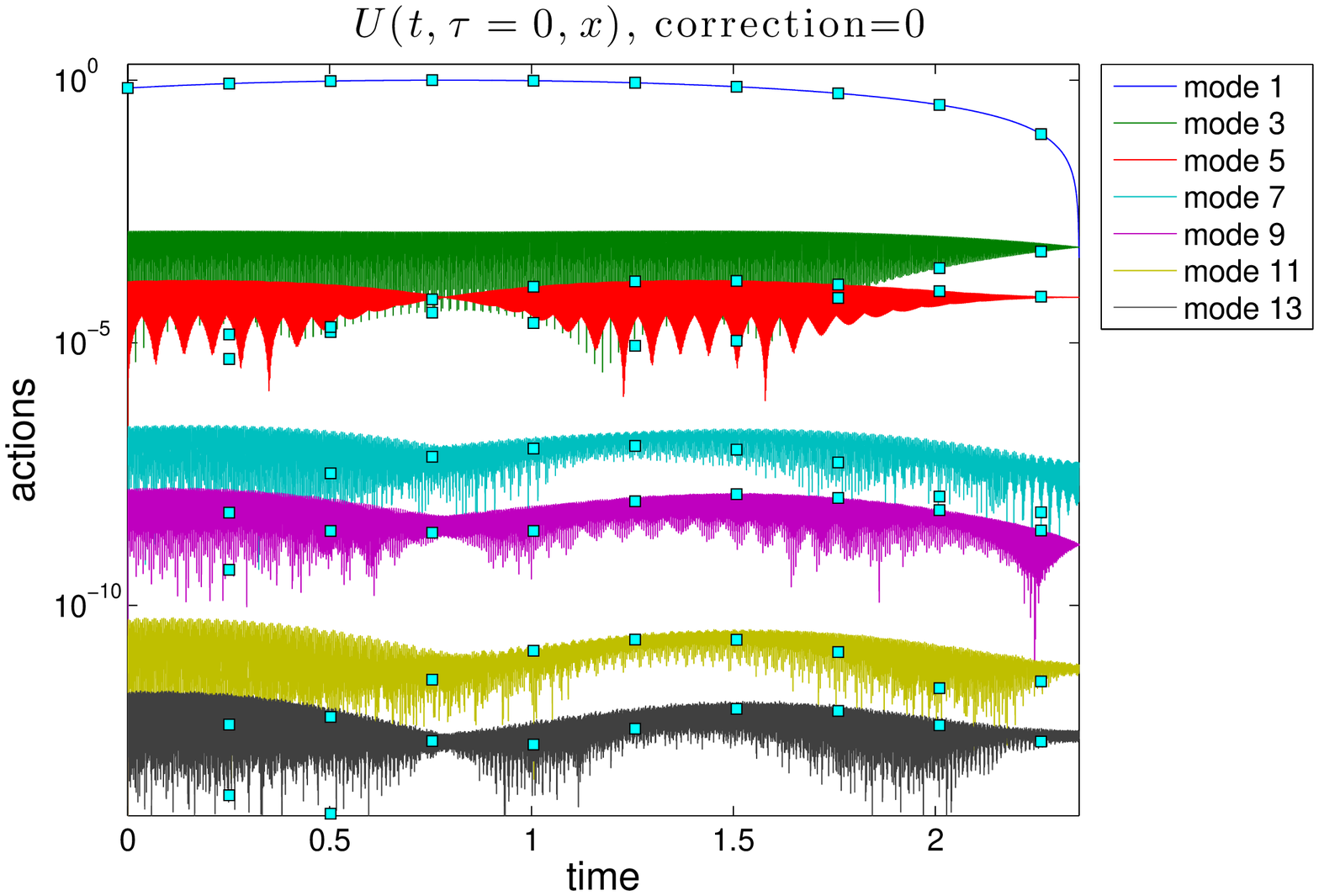}}}
    \caption{(NLS case) Time evolution of the first Fourier modes in $x$ of the function $U(t, \tau=0,x)$, in the log-scale.  At blue squares is plotted the reference solution at some stroboscopic points.}
  \label{fig1011}
\end{figure}

\begin{figure}[!htbp]
  \centerline{
  \subfigure[Reference solution obtained with the Strang splitting scheme]{\includegraphics[width=.717\textwidth]{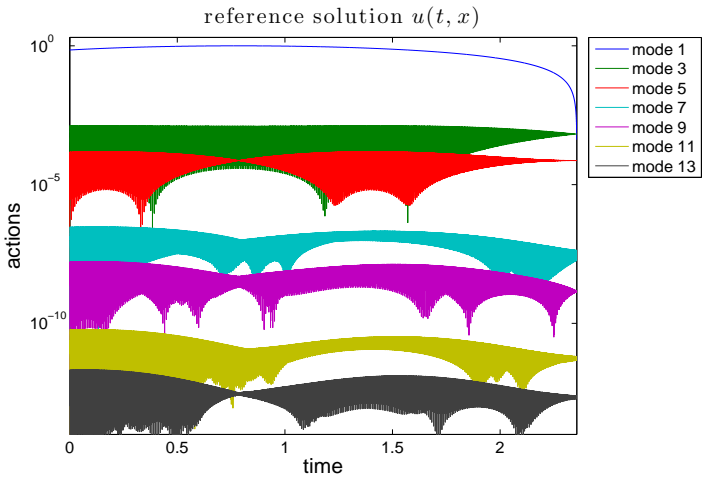}}\hspace*{-10mm}
  \subfigure[Numerical solution obtained with our scheme]{\includegraphics[width=.717\textwidth]{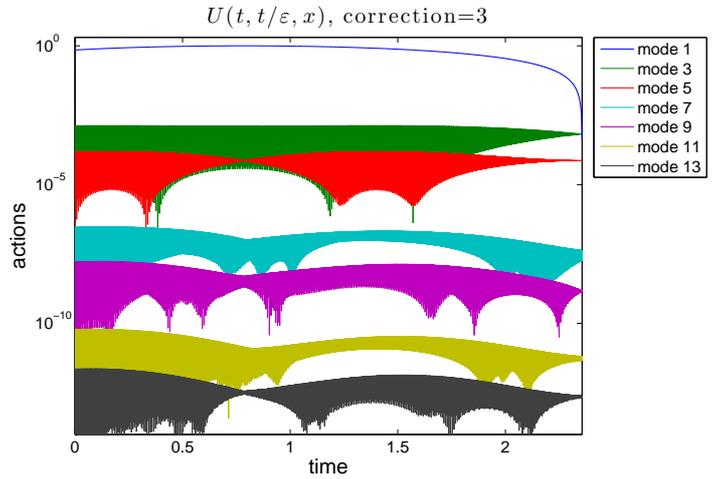}}}
    \caption{(NLS case) Time evolution of the first Fourier modes in $x$ of the solution $u(t,x)$ (in the log-scale).}
  \label{fig1012}
\end{figure}

Finally, in order to emphasize once again the importance of the choice of the initial data on the augmented problem in $U(t,\tau,x)$, we plot on Figures \ref{fig1010} and \ref{fig1011} the time evolution of the modulus of the first odd Fourier modes in $x$: $|u_1(t,\tau=0)|$, $|u_3(t,\tau=0)|$, \ldots, $|u_{13}(t,\tau=0)|$, with different initial data. Here we take $\eps=0.005$ and the steps in $t$, $\tau$ and $x$ are chosen small enough (we assume that the numerical schemes have reached their convergence). The NLS equation \eqref{nls} with the above choice of functions $\gamma$ and $u_0$ has a particular interesting property: as $\eps\to 0$, we have
$$u_1=\mathcal O(1),\quad (u_3,\, u_5)=\mathcal O(\eps),\quad (u_7,\, u_9)=\mathcal O(\eps^2),\quad (u_{11},\, u_{13})=\mathcal O(\eps^3).$$
This property allows to observe more easily the influence of the choice of the initial data. With uncorrected initial data (Figure \ref{fig1011}, right), all the terms of order $\mathcal O(\eps^k)$ with $k\geq 1$ are highly oscillatory. With the first order corrected initial data (Figure \ref{fig1011}, left), only the terms of order $\mathcal O(\eps^k)$ with $k\geq 2$ are rapidly oscillatory.  With the second order corrected initial data (Figure \ref{fig1010}, right), only the terms of order $\mathcal O(\eps^k)$ with $k\geq 3$ are rapidly oscillatory.  Finally, with the third order corrected initial data (Figure \ref{fig1010}, left), all the observed modes have smooth behaviors. Recall that, by construction, the solution of the augmented problem always satisfies $U(t,t/\eps,x)=u(t,x)$, so in particular we have the coincidence $U(t_k,\tau=0,x)=u(t_k,x)$ at the 'stroboscopic points' $t_k=2\pi k\varepsilon$, $k\in \N$. On Figures \ref{fig1010} and \ref{fig1011}, we plot in blue squares the modes of the solution $u$ of \eqref{nls} at the stroboscopic points $t_k$ for $k\in \{0,8,16,24,32,40,48,56,84,72\}$. We observe the coincidence between $U$ and $u$ at these times. As a comparison, the modes of the solution $u(t,x)$, which are all highly oscillatory (except for $|u_1|$ and $|u_{-1}|$), are finally represented  for all times on Figure \ref{fig1012} (on the left, computed with the Strang splitting scheme and on the right, computed with our UA scheme: both solutions coincide).

\end{document}